\documentclass{amsart}

\usepackage{amssymb}

\newtheorem{thm}{Theorem}%[section]

\newtheorem{lem}[thm]{Lemma}
\newtheorem{cor}[thm]{Corollary}

\newtheorem{prop}[thm]{Proposition}

\theoremstyle{definition}
\newtheorem{defn}[thm]{Definition}

\newtheorem{say}[thm]{}
\newtheorem{exmp}[thm]{Example}

\newtheorem{rem}[thm]{Remark}          
\newtheorem{aside}[thm]{Aside}          

\newtheorem{notation}[thm]{Notation}   
\newtheorem{warning}[thm]{Warning}  
\newtheorem{defn-thm}[thm]{Definition--Theorem}  %!!!!!!!!!!!!!!!!!!!!!!!!
\newtheorem{defn-lem}[thm]{Definition--Lemma}  %!!!!!!!!!!!!!!!!!!!!!!!!
\newtheorem{keyidea}[thm]{Key idea}

\theoremstyle{remark}

\renewcommand{\c}[0]{{\mathbb C}}  

\renewcommand{\o}[0]{{\mathcal O}} 
\newcommand{\z}[0]{{\mathbb Z}}

\renewcommand{\a}[0]{{\mathbb A}}

\newcommand{\p}[0]{{\mathbb P}}

\newcommand{\map}[0]{\dasharrow}
\newcommand{\qtq}[1]{\quad\mbox{#1}\quad}
\newcommand{\spec}[0]{\operatorname{Spec}}
\newcommand{\pic}[0]{\operatorname{Pic}}

\newcommand{\mult}[0]{\operatorname{mult}}

\newcommand{\supp}[0]{\operatorname{Supp}}    
    
\newcommand{\codim}[0]{\operatorname{codim}}    
\newcommand{\im}[0]{\operatorname{im}}

\newcommand{\aut}[0]{\operatorname{Aut}}

\newcommand{\sing}[0]{\operatorname{Sing}}    
\newcommand{\ex}[0]{\operatorname{Ex}}

\newcommand{\res}[0]{\operatorname{\mathcal R}}

\newcommand{\onto}[0]{\twoheadrightarrow}

\newcommand{\lcm}[0]{\operatorname{lcm}}

\newcommand{\ord}[0]{\operatorname{ord}}
\newcommand{\mord}[0]{\operatorname{max-ord}}
\newcommand{\der}[0]{\operatorname{Der}}
\newcommand{\cosupp}[0]{\operatorname{cosupp}}

\newcommand{\bres}[0]{\operatorname{\mathcal{BR}}} 
\newcommand{\bpri}[0]{\operatorname{\mathcal{BP}}} 
\newcommand{\bb}[0]{\operatorname{\mathcal{B}}} 
\newcommand{\bord}[0]{\operatorname{\mathcal{BO}}} 
\newcommand{\bdd}[0]{\operatorname{\mathcal{BD}}} 
\newcommand{\bmord}[0]{\operatorname{\mathcal{BMO}}} 
\newcommand{\extot}[0]{\operatorname{Ex_{\rm tot}}}

\def\into{\DOTSB\lhook\joinrel\to}

\begin{document}
\bibliographystyle{amsalpha}

\title{Resolution of  Singularities -- Seattle Lecture }
\author{J\'anos Koll\'ar}

\maketitle

\today

\tableofcontents

%\chapter{Strong Resolution in Characteristic Zero}

The most influential paper on resolution of singularities
is Hironaka's magnum opus \cite{hir-ann}\index{Hironaka!magnum opus}.
 Its starting point is
a profound shift in emphasis from resolving singularities of
varieties to resolving ``singularities of ideal sheaves.''
Ideal sheaves of smooth or  simple normal crossing divisors
 are the simplest ones. Locally, in a suitable coordinate system, these
 ideal sheaves are generated by a  single monomial.
The aim 
is to transform an arbitrary   ideal sheaf
into such a ``locally monomial''  one by a sequence
of blow-ups.
Ideal sheaves are much more flexible than varieties, and this 
opens up new ways of running induction.

Since then, resolution of singularities emerged as a very unusual
subject whose main object has been a deeper understanding
of the proof, rather than the search for new theorems.
A better grasp of the proof leads to improved theorems, with
the ultimate aim of extending the method to positive characteristic.
Two seemingly contradictory aspects make it very interesting
to study and develop Hironaka's approach.

  First, the method is very robust, in that many variants of the
proof work. One can even change basic definitions  and be
rather confident that the other parts can be modified to fit.

 Second, the complexity of the proof is very sensitive to details.
Small changes in definitions and presentation
may result in major simplifications.

This duality also makes it difficult to write a reasonable historical
presentation and to correctly appreciate the contributions of
various researchers. Each step ahead can be viewed as small or large,
depending on whether we focus on the change in the ideas or on their effect.
In some sense, all the results of the past forty years have their seeds
in \cite{hir-ann}, nevertheless, the improvement in the methods has been
enormous. 
Thus, instead of historical notes, here is a list of the
most important contributions to the development of the
Hironaka method, more or less in historical order:
Hironaka \cite{hir-ann, hir-ide};
Giraud \cite{gir};
Villamayor  \cite{vil-con, vil-pat, vil-int}  with his coworkers
Bravo \cite{bra-vil} and Encinas \cite{enc-vil, enc-vil-2};
Bierstone and Milman \cite{bm-uni, bm-sim, bm-can, bm-des};
Encinas and Hauser \cite{enc-hau} and 
W{\l}odar\-czyk \cite{wlo}.
The following proof 
relies mostly on the works of
Villamayor  and W{\l}odar\-czyk.

The methods of Bierstone and Milman and of Encinas and Hauser
differ from ours (and from each other) in key technical aspects,
though the actual resolution procedures end up very similar.

I have also benefited from the
surveys and books \cite{gir-bou, lip-arc, ahv, cgo-bk, 
res-bas, hau, cut}.

Abhyankar's  book  \cite{abh-book} shows some of
 the additional formidable difficulties that appear in 
positive characteristic.

A very elegant approach to resolution 
following 
de Jong's results on alterations \cite{deJ}
is developed in the papers \cite{bog-pan, abr-deJ, ab-wa}.
This method produces a resolution
as in (\ref{res.res.say}), which is however neither  strong
 (\ref{strongres.res.say}) nor functorial (\ref{functres.res.say}).
The version given in  \cite{par} is especially simple.

Another feature of the study of resolutions is that everyone seems
to use different terminology, so I also felt free to
introduce my own.

It is very instructive to compare the current methods with
Hironaka's  ``idealistic'' paper  
\cite{hir-ide}\index{Hironaka!idealistic exponents}.
The main theme is that resolution becomes simpler if we
do not try to control the process very tightly, as
 illustrated by the following three examples.

(1)  The original method of  \cite{hir-ann} worked with
the Hilbert-Samuel function  of an ideal sheaf
at  a point.
It was gradually realized that
the process simplifies if one considers only the
vanishing order of an ideal sheaf---a much cruder invariant.

(2) Two ideals $I$ and $J$ belong to the same {\it idealistic exponent} 
if they behave  ``similarly'' with respect to
any birational map $g$. 
(That is, $g^*I$ and $g^*J$ agree at the generic point of every
divisor for every $g$.)
 Now we see that 
it is easier to work with an equivalence relation
that requires the ``similar'' behavior only with respect to 
some birational maps
(namely, composites of smooth blow-ups
along subvarieties, where the vanishing order is maximal).

(3) The concept of a {\it distinguished presentation} 
\index{Hironaka!distinguished presentation} attempts to
pick a local coordinate system that is optimally adjusted to the
resolution of a variety or ideal sheaf. A key result of
W{\l}odarczyk \cite{wlo}
 says that for a suitably modified ideal,
all reasonable choices are equivalent, and thus we do not have to
be very careful. Local coordinate systems are not needed at all.
\medskip

The arguments given here differ from their predecessors in
two additional aspects. The first of these is  a matter of choice,
but the second one makes the structure of the proof
 patent.

(4) The inductive proof gives resolutions only locally, and
patching the local resolutions has been quite difficult.
The best way would be to define  an
 invariant on points of varieties 
$\operatorname{inv}(x, X)$ 
with values in an ordered set 
such that 

(i) $x\mapsto \operatorname{inv}(x, X)$ is an
 upper semi continuous function, and

(ii)  at each step of the resolution we 
blow up the locus where the invariant is maximal.

With some modification, this is 
 accomplished
in \cite{vil-con, bm-can, enc-hau}. 
All known  invariants are, however,
rather complicated.
W{\l}odar\-czyk
 suggested in \cite{wlo} that with his methods it should not be necessary to
define such an invariant.
We show that, by a slight change in the definitions, the resolution algorithm
automatically globalizes, obviating the need for the invariant.

(5) Traditionally, the results of Sections 9--12  constituted one
intertwined package, which had to be carried
through the whole induction together.
The introduction of the notions of
{\it $D$-balanced} and {\it MC-invariant} ideal sheaves
makes it possible to disentangle these to obtain
four independent parts.

\section{What is a good resolution algorithm?}

Before we consider the resolution of singularities in general,
it is worthwhile to contemplate what   the
properties of a good resolution algorithm should be.

Here I concentrate on the case of resolving  singularities of 
varieties only.
In practice, one may want to keep track and improve additional
objects, for instance,
  subvarieties or sheaves  as well, but
for now these variants would only obscure the general picture.

\begin{say}[Weakest resolution]
\index{Resolution!weak}{\it  Given  a variety $X$,
find a projective variety $X'$
such that $X'$ is smooth and birational to $X$.}

This is what the Albanese method  gives for curves and surfaces.
In these cases one can then use this variant
to get better resolutions, so we do not lose anything at the end.
These stronger forms are, however, not automatic, and it is
not at all clear that such a ``weakest resolution''
would be  powerful enough in higher dimensions.

(Note that even if $X$ is not proper, we have to insist
on $X'$ being proper; otherwise, one could  take
the open subset of smooth points of $X$ for $X'$.)

In practice it is
useful, sometimes crucial,  to have  additional properties.
\end{say}

\begin{say}[Resolution]\label{res.res.say}\index{Resolution}
 {\it  Given  a variety $X$,
find a  variety $X'$ and a projective morphism
$f:X'\to X$ 
such that $X'$ is smooth and $f$ is birational.}

This is the usual definition of resolution of singularities.

For many applications this is all one needs, but there are
plenty of situations when additional properties would be very useful.
Here are some of these.

\ref{res.res.say}.1 (Singularity theory).
\index{Resolution!and singularity theory} Let us start with an isolated
singularity $x\in X$. One frequently would like to study it
by taking a resolution $f:X'\to X$ and connecting the
properties of $x\in X$ with properties of the exceptional divisor
$E=\ex(f)$. Here everything works best if $E$ is projective,
that is, when $E=f^{-1}(x)$. 

It is reasonable to hope that we can achieve this. Indeed, by
assumption, $X\setminus \{x\}$ is smooth, so it
should be possible to resolve without changing $X\setminus \{x\}$. 

\ref{res.res.say}.2 (Open varieties).
\index{Resolution!and open varieties} It is  natural to study
a noncompact variety $X^0$ via a compactification $X\supset X^0$.
Even if $X^0$ is smooth, the compactifications that are easy to
obtain are usually singular. Then one would like to resolve
 the singularities of $X$ and get a smooth compactification
$X'$. If we take  any resolution  $f:X'\to X$, 
the embedding $X^0\into X$ does not lift to an embedding
$X^0\into X'$.  Thus we would like to find a
resolution   $f:X'\to X$ such that $f$ is an isomorphism over $X^0$.

In both of the above examples, we would like
the exceptional set $E$ or the boundary $X'\setminus X^0$
to be ``simple.'' Ideally we would like them to be smooth,
but this is rarely possible. The next best situation is when
$E$ or $X'\setminus X^0$ are simple normal crossing divisors.

These considerations lead to the  following variant.
\end{say}

\begin{say}[Strong resolution]\label{strongres.res.say}
\index{Resolution!strong}
 {\it  Given  a variety $X$,
find a  variety $X'$ and a projective morphism
$f:X'\to X$ 
such that we have the following:
\begin{enumerate}
\item $X'$ is smooth and $f$ is birational,
\item $f: f^{-1}(X^{ns})\to X^{ns}$ is an isomorphism, and
\item $f^{-1}(\sing X)$ is a divisor with simple normal crossings.
\end{enumerate}
}
Here $\sing X$ denotes the set of {\it singular points} of
$X$\index{Singular points, $\sing X$}\index{SingX@$\sing X$, singular points}
and $X^{ns}:=X\setminus \sing X$ the  set of
{\it smooth points}.
\index{Smooth points, $X^{ns}$}\index{Xn@$X^{ns}$, smooth points}

Strong resolution seems to be the variant that is most frequently
used in applications, but 
sometimes other versions are needed. For instance, one might  
need  condition (\ref{strongres.res.say}.3) scheme theoretically.

A more important question arises when one has several varieties
$X_i$ to work with simultaneously. In this case we may need to know
that certain morphisms $\phi_{ij}:X_i\to X_j$ lift to the
resolutions $\phi'_{ij}:X'_i\to X'_j$.

It would be nice to have this for all morphisms, which would give a
``resolution functor'' from the category of all varieties 
and morphisms to the
category of smooth varieties. This is, however, impossible.
\medskip
 
{\it Example} \ref{strongres.res.say}.4.
Let $S:=(uv-w^2=0)\subset \a^3$ be the quadric cone,
and consider the morphism
$$
\phi:\a^2_{x,y}\to S\qtq{given by $(x,y)\mapsto (x^2,y^2,xy)$.}
$$
The only sensible resolution of
$\a^2$ is itself, and any resolution of $S$ dominates the
minimal resolution $S'\to S$ obtained by blowing up the origin.

The morphism $\phi$  lifts to a rational map
$\phi':\a^2\map S'$, but $\phi'$  is not a morphism.
\end{say}

It seems that the best one can hope for is that the
resolution commutes with smooth morphisms.

\begin{say}[Functorial resolution]\label{functres.res.say}
\index{Resolution!functorial}\index{Functoriality!of resolution}
 {\it  For every variety $X$
find a   resolution
$f_X:X'\to X$ that is functorial with respect to smooth morphisms.
That is,  any smooth morphism $\phi:X\to Y$
lifts to a smooth morphism $\phi':X'\to Y'$,
which gives a fiber product square
 $$
\begin{array}{ccc}
X' & \stackrel{\phi'}{\longrightarrow} & Y'\\
f_X\downarrow \hphantom{f_X} &\square & \hphantom{f_Y}\downarrow f_Y\\
 X & \stackrel{\phi}{\longrightarrow} & Y.
\end{array}
$$
}

Note that if $\phi'$ exists, it is unique, and so we indeed get a functor
 from the category of all varieties 
and smooth morphisms to the
category of smooth varieties and smooth morphisms.

This is   quite a strong property
with many useful implications.

\ref{functres.res.say}.1 (Group actions).
\index{Resolution!and group actions}Functoriality of resolutions  implies that
any group action on $X$ lifts to $X'$.
For discrete groups this is just functoriality
plus the observation that the only lifting of the identity map 
on $X$ is the identity map of $X'$. 
For an algebraic group $G$ a few more
steps are needed; see (\ref{functres.res.aside}.1).

\ref{functres.res.say}.2  (Localization). 
\index{Resolution!and localization}
Let $f_X:X'\to X$ be a functorial resolution.
The embedding of any open subset $U\into X$ is smooth, and so
the  functorial resolution of $U$ is the restriction of
the  functorial resolution of $X$. That is,
$$
(f_U:U'\to U)\cong \bigl(f_X|_{f_X^{-1}(U)}: f_X^{-1}(U)\to U\bigr).
$$
Equivalently, 
a  functorial resolution is Zariski local. 
 More generally, a functorial resolution is \'etale local
since  \'etale morphisms are smooth.

Conversely, we show in
(\ref{functres.res.aside}.2) that any resolution that is
functorial with respect to \'etale morphisms is also 
functorial with respect to smooth morphisms.

Since any resolution $f:X'\to X$ is birational, it is
an isomorphism over some smooth points of $X$.
Any two smooth points of $X$ are \'etale equivalent, and thus
a resolution that is functorial with respect to \'etale morphisms
is an isomorphism over smooth points.
Thus any functorial resolution satisfies
(\ref{strongres.res.say}.2).

\ref{functres.res.say}.3 (Formal localization).
\index{Resolution!and formal localization}
Any sensible \'etale local construction in algebraic geometry
is also formal local. In our case this means that
the behavior of the resolution $f_X:X'\to X$ near a point
$x\in X$ should depend only on the completion $\widehat{\o}_{x,X}$.
(Technically speaking, $\spec \widehat{\o}_{x,X}$ is not a variety and the map
$\spec \widehat{\o}_{x,X}\to \spec \o_{x,X}$ is only formally smooth,
so this is a stronger condition than functoriality.)

\ref{functres.res.say}.4 (Resolution of products).
\index{Resolution!of products} 
It may appear surprising, but a strong and functorial  resolution
should {\em not} commute with products.

For instance, consider the quadric cone
$0\in S=(x^2+y^2+z^2=0)\subset \a^3$.
This is resolved by blowing up the origin
$f:S'\to S$ with exceptional curve $C\cong \p^1$. On the other hand, 
$$
f\times f:S'\times S'\to S\times S
$$
cannot be the outcome of an \'etale local strong resolution.
The singular locus of $S\times S$
has two components, $Z_1=\{0\}\times S$ and $Z_2=S \times \{0\}$,
and correspondingly, the exceptional divisor has two components,
 $E_1=C\times S'$ 
and $E_2=S'\times C$,  which intersect along  $C\times C$.

If we work \'etale locally at $(0,0)$, we cannot tell whether the
two branches of the singular locus $Z_1\cup Z_2$
are on different irreducible components of $\sing S$ or on one
non-normal  irreducible component. Correspondingly,
the germs of $E_1$ and $E_2$ could be on the same
irreducible exceptional divisor, and on a strong
resolution self-intersections of exceptional divisors are not allowed.

\end{say}

So far we concentrated on the end result $f_X:X'\to X$ of the
resolution. Next we look at some properties of the
resolution algorithm itself.

\begin{say}[Resolution by blowing up smooth centers]
\index{Resolution!by smooth blow-ups}
 {\it  For every variety $X$
find a   resolution
$f_X:X'\to X$
such that $f_X$ is a composite of morphisms
$$
f_X:X'=X_n\stackrel{p_{n-1}}{\longrightarrow}
X_{n-1}\stackrel{p_{n-2}}{\longrightarrow}\cdots
\stackrel{p_{1}}{\longrightarrow}
X_{1}\stackrel{p_{0}}{\longrightarrow}X_0=X,
$$
where each $p_i:X_{i+1}\to X_i$ is obtained by
blowing up a smooth subvariety $Z_i\subset X_i$.
}

If we want $f_X:X'\to X$ to be a strong resolution, then
the condition $Z_i\subset \sing X_i$ may also be required,
though we  need only that
$p_0\cdots p_{i-1}(Z_i)\subset \sing X$.

Let us note first that in low dimensions some of the best
resolution algorithms do not have this property.
\begin{enumerate}
\item The quickest way to resolve a curve is to
normalize it. The normalization usually cannot be obtained by
blowing up points (though it is a composite of blow-ups of points).
\item  A normal surface can be resolved
by repeating the procedure:
``blow up the singular points and  normalize'' \cite{zar-40}.
\item  A toric variety is best resolved by toric blow-ups. These are
 rarely given by blow-ups of subvarieties (cf.\ \cite[2.6]{fult-toric}).
\item Many of the best-studied singularities  
are easier to resolve by doing a weighted blow-up first.
\item 
The theory of Nash blow-ups offers a---so far mostly hypothetical---approach
to resolution that does not rely on blowing up smooth centers;
cf.\ \cite{hir-nash}.

\end{enumerate}

On the positive side, resolution by blowing up smooth centers
has the great advantage that we do not mess up what is already nice.
For instance, if we want to resolve $X$ and $Y\supset X$ is a smooth
variety containing $X$, then a resolution by blowing up smooth centers
automatically carries along the smooth variety. Thus we get a
sequence of smooth varieties $Y_i$ fitting in a diagram
$$
\begin{array}{ccccccl}
X_n&\stackrel{p_{n-1}}{\longrightarrow}&
X_{n-1}&\cdots
&X_{1}&\stackrel{p_{0}}{\longrightarrow}&X_0=X\\
\downarrow &&\downarrow &&\downarrow &&\ \downarrow \\
Y_n&\stackrel{q_{n-1}}{\longrightarrow}&
Y_{n-1}&\cdots
&Y_{1}&\stackrel{q_{0}}{\longrightarrow}&Y_0=Y,
\end{array}
$$
where the vertical arrows are closed embeddings.
\end{say}

Once we settle on resolution by successive blowing ups,
the main question is how to find the centers that we need to blow up.
From the algorithmic point of view, 
the best outcome would be the following.

\begin{say}[Iterative resolution, one blow-up at a time]
\label{it.res.1bu.say}
\index{Resolution!iterative}\index{Resolution!one blow-up at a time}
{\it 
For any variety $X$, 
 identify a  subvariety $W(X)\subset X$ consisting of the
``worst'' singularities. Set $R(X):=B_{W(X)}X$
and $R^m(X):=R(R^{m-1}(X))$ for $m\geq 2$.
Then we get  resolution by iterating this procedure.
That is, $R^m(X)$ is smooth for $m\gg 1$.
}

Such an algorithm exists for curves with $W(X)=\sing X$.

The situation is not so simple in higher dimensions.
\medskip

{\it Example} \ref{it.res.1bu.say}.1.
 Consider the pinch point, or Whitney umbrella,
\index{Resolution!of $x^2-y^2z=0$}
  $S:=(x^2-y^2z=0)\subset \a^3$. 
$S$ is singular along the line $(x=y=0)$. It has a normal crossing point
if $z\neq 0$ but a more complicated singularity at $(0,0,0)$.

If we blow up the ``worst'' singular point  $(0,0,0)$ of the surface $S$,
then in the chart  with coordinates $x_1=x/z,y_1=y/z, z_1=z$ we
get the birational transform
$S_1=(x_1^2-y_1^2z_1=0)$.
This is isomorphic to the original surface.

Thus we conclude that one cannot resolve surfaces  by blowing up the
``worst'' singular point all the time.

We can, however, resolve the pinch point by blowing up the whole
singular line.  In this case, using the 
multiplicity (which is a rough invariant) gives  the right
blow-up, whereas distinguishing the pinch point from a normal crossing
point (using some finer invariants) gives the wrong blow-up.
 The message is that we should not look at the
singularities too carefully.

The situation gets even worse for normal 3-folds.

\medskip

 {\it Example} \ref{it.res.1bu.say}.2. Consider the 3-fold
\index{Resolution!of $x^2+y^2+z^mt^m=0$}
$$
X:=(x^2+y^2+z^mt^m=0)\subset \a^4.
$$
The singular locus is the union of the two lines
$$
L_1:=(x=y=z=0)\qtq{and} L_2:=(x=y=t=0).
$$
There are two reasons why  no sensible resolution procedure
should  start by blowing up either of the lines.
\begin{enumerate}
\item[(i)] 
The two lines are interchanged by the involution
$\tau: (x,y,z,t)\mapsto(x,y,t,z)$, and thus they should be
blown up in a $\tau$-invariant way.
\item[(ii)]  An \'etale local resolution procedure cannot tell
if $L_1\cup L_2$ is a union of two lines or just two local
branches of an irreducible curve. Thus picking one branch does not
make sense globally.
\end{enumerate}
Therefore, we must start by blowing up the intersection point
$(0,0,0,0)$ (or resort to blowing up a singular subscheme). 

Computing the $t$-chart $x=x_1t_1, y=y_1t_1, z=z_1t_1, t=t_1$, we get
$$
X_{1,t}=(x_1^2+y_1^2+z_1^mt_1^{2m-2}=0)
$$
and similarly in the $z$-chart. Thus on $B_0X$
the singular locus consists of three lines: $L'_1, L'_2$ and
an exceptional line $E$.

For $m=2$ we are thus back to the original situation,  
and for $m\geq 3$ we made the singularities worse by blowing up.
In 
the $m=2$ case there is nothing else one can do,
 and we get our first negative result.

\medskip

{\it Claim} \ref{it.res.1bu.say}.3. 
There is no iterative resolution algorithm that
works one smooth blow-up at a time.
\index{Resolution!failure of iterative method} \end{say}

\medskip

The way out is to notice that our two objections 
(\ref{it.res.1bu.say}.2.i--ii)  to
first blowing up one of the lines $L_1$ or $L_2$  
are not so strong when applied to the three lines $L_1$, $L_2$ and $E$
on the blow-up $B_0X$. Indeed, we know that
the new exceptional line $E$ is isomorphic to $\c\p^1$, and it is
invariant under every automorphism lifted from $X$.
Thus we can safely blow up $E\subset B_0X$
without the risk of running into problems with 
\'etale localization. 
(A key point is that we want to ensure  that the process is
\'etale local only on $X$, not on all the intermediate varieties.)
In the $m=2$ case we can then blow up the birational transforms of
the two lines $L_1$ and $L_2$ simultaneously, to achieve
resolution. (Additional steps are needed for $m\geq 3$.)

In general,  we have to  ensure that the
resolution process has some ``memory.''
That is, at each step 
the procedure is allowed to use information about the previous blow-ups.
For instance, it could  keep track of
  the exceptional divisors that
were created by earlier blow-ups of the resolution and  in which order
they were created.

\begin{say}[Other considerations]\label{res.other.say}

There are several other ways to judge how good a
resolution algorithm is.

\ref{res.other.say}.1 (Elementary methods).
A good resolution method should be part of
``elementary'' algebraic geometry. 
Both Newton's method of rotating rulers  and the Albanese projection
method  pass this criterion. On the other hand, several of the
methods for surfaces 
rely on more advanced machinery.

\ref{res.other.say}.2 
(Computability).
\index{Resolution!computability}  In concrete cases, one may wish to explicitly
determine
resolutions by hand or by a  computer.
As far as I can tell, the existing methods do rather poorly even
on the simplest singularities. In a more theoretical direction, one can ask
for the worst case or average complexity of the algorithms.
See \cite{bs1, bs2, fk-p} for  computer implementations.

\end{say}

\begin{say} Our resolution 
is strong and functorial with respect to smooth morphisms, but
it is very far from being 
iterative if we want to work one blow-up at a time.
Instead, at each step we specify a long sequence of
blow-ups to be performed. 

We shift our emphasis from resolution of singularities to
principalization of ideal sheaves. While principalization 
is achieved by a sequence of smooth blow-ups, the resolution 
of singularities may involve blow-ups of singular centers.
Furthermore,  at some stage we may blow
up a subvariety $Z_i\subset X_i$ along which the variety $X_i$ is smooth.
This  only happens for subvarieties that sit over
the original singular locus, so at the end we still get a strong
resolution.

The computability of the algorithm has not been investigated
much, but the early indications are not promising.
One issue is that starting with, say, a hypersurface $(f=0)\subset \a^n$ of
multiplicity $m$ the first step is to
replace the ideal $(f)$ with another ideal
$W(f)$, which has more than $e^{mn}$ generators, each of
multiplicity at least $e^m$; see (\ref{tuning.gen.say}.3).
Then we reduce to a resolution problem in $(n-1)$-dimensions,
and at the next reduction step we again may have an exponential
increase of the multiplicity and the number of generators.
For any reasonable computer implementation, some shortcuts
are essential.

\end{say}

\begin{aside}\label{functres.res.aside}  Here we prove the two claims made
 in (\ref{functres.res.say}).
These are not used in the rest of the chapter.
\medskip
 
 {\it Proposition} \ref{functres.res.aside}.1. 
The action of an algebraic group $G$ on a scheme $X$ 
lifts to an  action of  $G$ on its functorial resolution $X'$. 
\medskip

Proof. 
 The action of an algebraic group $G$ on a variety $X$ 
is given by
a smooth morphism
$m:G\times X\to X$. By functoriality,  the resolution $(G\times X)'$ of
$G\times X$ is given by the pull-back of $X'$ via $m$, that is, by
$f_X^*(m):(G\times X)'\to X'$.

On the other hand, the second projection
$\pi_2:G\times X\to X$ is also smooth, and so
 $(G\times X)'=G\times X'$. Thus we get a commutative diagram
$$
\begin{array}{ccccc}
m':G\times X'\hphantom{m':} & \cong & (G\times X)'& 
\stackrel{f_X^*(m)}{\longrightarrow}& X'\\
\hphantom{id_G\times f_X} \downarrow id_G\times f_X &&&& 
\hphantom{f_X}\downarrow  f_X\\
G\times X & = & G\times X & \stackrel{m}{\longrightarrow}& X.\\
\end{array}
$$
We claim that the composite in the top row
$m':G\times X'\to X'$ defines a group action.
This means that the following diagram is commutative,
where $m_G:G\times G\to G$ is the group multiplication:
$$
\begin{array}{ccc}
G\times G\times X' & \stackrel{id_G\times m'}{\longrightarrow}& G\times X'\\
m_G\times id_{X'} \downarrow \hphantom{m_G\times id_{X'}} &&
\hphantom{m'}\downarrow  m'\\
G\times X' &  \stackrel{m'}{\longrightarrow}& X'.\\
\end{array}
$$
Since $m:G\times X\to X$ defines a group action, we know that
the diagram is commutative over a dense open set. 
Since all schemes in the diagram are separated and reduced,
this implies commutativity.\qed
\medskip

  {\it Proposition} 
\ref{functres.res.aside}.2.  Any resolution that is
functorial with respect to \'etale morphisms is also 
functorial with respect to smooth morphisms.
\medskip

Proof. As we noted in 
(\ref{functres.res.say}.2), 
a resolution that is functorial with respect to \'etale morphisms
is an isomorphism over smooth points.

\'Etale locally, a smooth morphism is a direct product,
and so it is sufficient to prove that $(X\times A)'\cong X'\times A$
for any abelian variety $A$. Such an isomorphism is unique;
thus it is enough to prove existence for $X$ proper.

Since $(X\times A)'$ is proper, the connected component of
its automorphism group is an algebraic group $G$ 
(see, for instance, \cite[I.1.10]{rcbook}).
Let $G_1\subset G$ denote the subgroup  whose elements commute with
the projection  $\pi:(X\times A)'\to X$.

Let $Z\subset X^{ns}$ be a finite subset. Then
$\pi^{-1}(Z)\cong Z\times A$, and the action of $A$ on itself
gives a subgroup $j_Z:A\into \aut(\pi^{-1}(Z))$.
There is a natural restriction map
$\sigma_Z:G_1\to \aut(\pi^{-1}(Z))$; set
$G_Z:=\sigma_Z^{-1}(j_Z A)$.

As we increase $Z$, the subgroups $G_Z$ form a decreasing
sequence, which eventually stabilizes at a 
subgroup $G_X\subset G$ 
such that for every finite set $Z\subset X^{ns}$ the action
of $G_X$ on $\pi^{-1}(Z)$ is through the action of $A$ on itself.
This gives an injective   homomorphism of algebraic groups
$G_Z\into A$.

On the other hand, $A$ acts on $X\times A$ by isomorphisms, and
by assumption this action lifts to an action of the {\em discrete} group
$A$ on $(X\times A)'$. Thus the injection  
$G_Z\into A$ has a set-theoretic inverse, so it is an
isomorphism of algebraic groups.\qed
\end{aside}

\section{Examples of resolutions}

We start the study of resolutions with some examples.
First, we  describe how the resolution method  deals with
two particular surface singularities $S\subset \a^3$.
 While these are  relatively simple cases,
they allow us to isolate six  problems facing the method.
Four of these we solve later, and we can live with the other two.

Then we see how the problems can be tackled for
Weierstrass polynomials and what this solution
tells us about the general case.
For curves and surfaces, this method was already used
in Sections 1.10 and 2.7.

\begin{keyidea} We look at the trace of $S\subset \a^3$
 on a suitable smooth surface $H\subset \a^3$
and reconstruct the whole resolution of $S$ from $S\cap H$.
\end{keyidea}

More precisely, starting with a surface singularity
$0\in S\subset \a^3$ of multiplicity $m$, we will be guided by
$S\cap H$  until the multiplicity of
the birational transform %%(\ref{bir.tr.var.defn})
 of $S$ drops below $m$.
Then we need to 
repeat the method to achieve further multiplicity reduction.

\begin{exmp}[Resolving $S:=(x^2+y^3-z^6=0)\subset \a^3$]
\label{deg1Econe.w.meth} 
\index{Resolution!of $x^2+y^3-z^6=0$}
(We already know %%from (\ref{x^2+y^3+z^6.exmp}) 
that the minimal resolution
has a single exceptional curve $E\cong (x^2z+y^3-z^3=0)\subset \p^2$ and it has
self-intersection $(E^2)=-1$ but let us forget  it for now.)

Set $H:=(x=0)\subset \a^3$, and  work with $S\cap H$.
\medskip

{\it Step 1.}  Although the trace $S\cap H=(y^3-z^6=0)\subset \a^2$
has multiplicity 3, we came from a multiplicity 2 situation,
and we blow up until the multiplicity drops below 2.
\medskip

Here it takes two blow-ups to achieve this.
The crucial local  charts and equations are
$$
\begin{array}{ll}
x^2+y^3-z^6=0,& \qquad \\
x_1^2+(y_1^3-z_1^3)z_1=0,& \qquad  x_1=x/z, y_1=y/z, z_1=z,\\
x_2^2+(y_2^3-1)z_2^2=0,&\qquad 
 x_2=x_1/z_1, y_2=y_1/z_1, z_2=z_1.
\end{array}
$$
At this stage the trace of the  dual graph %%(\ref{dual.graphs.defn}) 
 of the 
birational transform of $S$ on the
birational transform of $H$ 
is the following, where the numbers indicate the multiplicity
(and not minus the self-intersection number  as usual) 
and $\bullet$ indicates the birational transform of the
original curve $S\cap H$:
$$
1\ -\ 2
\begin{array}{c}
\diagup\\
-\\
\diagdown
\end{array}
\begin{array}{c}
\bullet\\
{\ }\\
\bullet\\
{\ }\\
\bullet
\end{array}
$$
\medskip

{\it Step 2.} The birational transform of $S\cap H$ intersects
some of the new exceptional curves that appear with
positive coefficient. We blow up until these intersections are
removed.
\medskip

In our case each intersection point needs to be blown up twice.
At this stage the trace of the birational transform of $S$ on the
birational transform of $H$ 
looks like
$$
1\ -\ 2\  
\begin{array}{c}
\diagup\\
-\\
\diagdown
\end{array}
\begin{array}{ccccc}
 1 & - & 0 & - &\bullet\\
&&&&\\
 1 & - & 0 & - &\bullet\\
&&&&\\
 1 & - & 0 & - &\bullet
\end{array}
$$
where multiplicity $0$ indicates that the curve is no
longer contained in the birational transform of $H$
(so strictly speaking, we should not draw it at all). 
\medskip

{\it Step 3.} The trace now has multiplicity $< 2$ along the 
birational transform of $S\cap H$, but it still has some points
of multiplicity $\geq 2$. We remove these by blowing up the
exceptional curves with multiplicity $\geq 2$.
\medskip

In our case there is only one such curve. After blowing it up,
we get the final picture
$$
1\ -\ \boxed{0}\  
\begin{array}{c}
\diagup\\
-\\
\diagdown
\end{array}
\begin{array}{ccccc}
 1 & - & 0 & - &\bullet\\
&&&&\\
 1 & - & 0 & - &\bullet\\
&&&&\\
 1 & - & 0 & - &\bullet
\end{array}
$$
where the boxed curve is elliptic.
\end{exmp}

More details of the
resolution method appear in the following example.

\begin{exmp}[Resolving $S:=(x^3+(y^2-z^6)^2+z^{21})=0)\subset \a^3$]
\label{messy.exmp.w.meth} 
\index{Resolution!of $x^3+(y^2-z^6)^2+z^{21}=0$}
As before, we look at the trace of $S$ on the plane $H:=(x=0)$
and reconstruct the whole resolution of $S$ from $S\cap H$.
\medskip

{\it Step 1.}  Although the trace $S\cap H=((y^2-z^6)^2+z^{21}=0)\subset \a^2$
has multiplicity 4, we came from a multiplicity 3 situation,
and we blow up until the multiplicity drops below 3.
\medskip

Here it takes three blow-ups to achieve this.
The crucial local charts and equations  are
$$
\begin{array}{ll}
x^3+(y^2-z^6)^2+z^{21}=0,& \qquad \\
x_1^3+z_1(y_1^2-z_1^4)^2+z_1^{18}=0,&
 \qquad  x_1=x/z, y_1=y/z, z_1=z,\\
x_2^3+z_2^2(y_2^2-z_2^2)^2+z_2^{15}=0&\qquad 
 x_2=x_1/z_1, y_2=y_1/z_1, z_2=z_1,\\
x_3^3+z_3^3(y_3^2-1)^2+z_3^{12}=0,&\qquad 
 x_3=x_2/z_2, y_3=y_2/z_2, z_3=z_2.\end{array}
$$
The birational transform of $S\cap H$ has equation
$$
(y_3^2-1)^2+z_3^{9}=0
$$
and has two higher cusps at $y_3=\pm 1$ on the last exceptional curve.
The trace of the birational transform of $S$ on the
birational transform of $H$ 
looks like
$$
1\ -\ 2\ -\ 3\ 
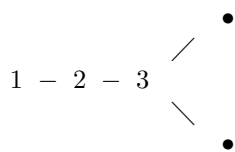
\begin{array}{c}
\diagup\\
{\ }\\
\diagdown
\end{array}
\begin{array}{c}
\bullet\\
{\ }\\{\ }\\{\ }\\
\bullet
\end{array}
$$
(As before, the numbers indicate the multiplicity,
and $\bullet$ indicates the birational transform of the
original curve $S\cap H$.
Also note that here the curves marked $\bullet$ have multiplicity 2 at 
their intersection point with the curve marked $3$.)
\medskip

{\it Step 2.} The birational transform of $S\cap H$ intersects
some of the new exceptional curves that appear with
positive coefficient. We blow up until these intersections are
removed.
\medskip

In our case each intersection point needs to be blown up three times,
and we get the following picture:
$$
1\ -\ 2\ -\ 3\ 
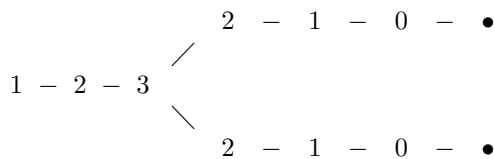
\begin{array}{c}
\diagup\\
{\ }\\
\diagdown
\end{array}
\begin{array}{ccccccc}
2 & - & 1 & - & 0 & - &\bullet\\
&&&&&&\\&&&&&&\\&&&&&&\\
2 & - & 1 & - & 0 & - &\bullet
\end{array}
$$

\medskip

{\it Step 3.} The trace now has multiplicity $< 3$ along the 
birational transform of $S\cap H$, but it still has some points
of multiplicity $\geq 3$. 
There is one 
exceptional curve with multiplicity $\geq 3$; we blow that up.
This drops its coefficient from $3$ to $0$. 
There are four more points of multiplicity $3$,
where a curve with multiplicity 2 intersects a curve 
with multiplicity 1. After blowing these up
we get the final picture
$$
1\ -\ 0\ -\ 2\ -\ 0\ 
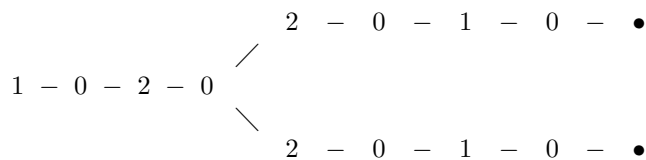
\begin{array}{c}
\diagup\\
{\ }\\
\diagdown
\end{array}
\begin{array}{ccccccccc}
2 & - &  0 & - &1 & - & 0 & - &\bullet\\
&&&&&&&&\\&&&&&&&&\\&&&&&&&&\\
2 & - & 0 & - &1  & - & 0 & - &\bullet
\end{array}
$$

\end{exmp}

\begin{say}[Problems with the method]\label{probls.with.meth.say}
\index{Resolution!problems of inductive approach}
\index{Problems of the inductive approach}
There are at least six different problems with the method.
Some are clearly visible from the examples, while some are hidden by the
presentation.
\medskip

{\it Problem}
\ref{probls.with.meth.say}.1.  In (\ref{deg1Econe.w.meth}) 
we end up with eight exceptional curves, when we
need only one to resolve $S$. %%; see  (\ref{x^2+y^3+z^6.exmp}).
 In general, for many surfaces the
method gives  a resolution that is much bigger than the
minimal one. However, in higher dimensions there is no
minimal resolution, and it is not clear how to measure
the ``wastefulness'' of a resolution. 

We will not be able to
deal with this issue.

\medskip

{\it Problem}
\ref{probls.with.meth.say}.2. 
 The resolution problem for surfaces in $\a^3$ was reduced
not to the resolution problem for curves in $\a^2$ but
to a related problem that also takes into account
 exceptional curves and their multiplicities in some way.

We  have to set up a somewhat artificial-looking
resolution problem that allows true induction on the dimension.

\medskip

{\it Problem}
\ref{probls.with.meth.say}.3.
 The end result of the resolution process guarantees that
the birational transform of $S$ has multiplicity $<2$ along the
birational transform of $H=(x=0)$, but we have said nothing
about the 
singularities that occur outside the
birational transform of $H$. 

There are indeed such singularities  if we do not choose $H$ carefully.
For instance, if we take $H':=(x-z^2=0)$,
then at the end of Step.1 of (\ref{deg1Econe.w.meth}),
 that is, after two blow-ups, the birational transform of $H'$
is $(x_2-1=0)$, which does not contain the
singularity 
that is at the origin $(x_2=y_2=z_2=0)$.

Thus a careful choice of $H$ is needed.
This is solved by the theory of {\it maximal contact},
developed by Hironaka and  Giraud \cite{gir, ahv-mc}.

\medskip

 {\it Problem}
\ref{probls.with.meth.say}.4.
 In some cases, the opposite problem happens.
All the singularities end up on the birational transforms
of $H$, but we also pick up extra tangencies, so
we see too many singularities.

For instance, take $H'':=(x-z^3=0)$. Since
$$
x^2+y^3-z^6=(x-z^3)(x+z^3)+y^3,
$$
the trace of $S$ on $H''$ is a triple line.
The trace shows a 1-dimensional singular set when
we have only an isolated singular point.

In other cases, these problems  may appear only after many blow-ups.

At  first glance, this may not  be a problem at all.
This simply means that we make some  unnecessary blow-ups as well.
Indeed, if our aim is to resolve surfaces only, then this problem can be
mostly ignored. However, for the general inductive procedure
this is a serious difficulty since unnecessary blow-ups
can increase the multiplicity.
For instance, 
$$
S=(x^4+y^2+yz^2=0)\subset \a^3
$$
is an isolated double point. If we blow up the line
$(x=y=0)$, in the $x$-chart we get a triple point
$$
x_1^3+x_1y_1^2+y_1z^2=0,\qtq{where $x=x_1, y=y_1x_1$.}
$$
One way to solve this problem is to switch from resolving
varieties to ``resolving'' ideal sheaves by  introducing a
{\it coefficient ideal}  $C(S)$ such that 
\begin{enumerate}
\item[(i)] resolving $S$ is equivalent to ``resolving'' $C(S)$, and
\item[(ii)] ``resolving'' the traces $C(S)|_H$ does not generate
extra blow-ups for $S$. 
\end{enumerate}
This  change of emphasis is crucial for our approach.
\medskip

{\it Problem}
\ref{probls.with.meth.say}.5.
 No matter how carefully we choose $H$, we can never
end up with a unique choice.
For instance, the analytic automorphism of $S=(x^2+y^3-z^6=0)$,
$$
(x,y,z)\mapsto (x+y^3, y\sqrt[3]{1-2x-y^3}, z),
$$
shows that no internal property
distinguishes the choice $x=0$
from the choice $x+y^3=0$.

Even with the careful ``maximal contact'' choice of $H$,
we end up with cases where the traces $S\cap H$ are not 
isomorphic. Thus our resolution process seems to depend on the choice of $H$.

This is again only a minor inconvenience for surfaces, 
but in higher dimensions we have to deal with patching together
the local resolution processes into a global one.
(We cannot  even avoid  this issue
by pretending to care only about  isolated singularities, since
 blowing up frequently leads  to
nonisolated singularities.)

An efficient solution of this  problem 
developed in  \cite{wlo}   replaces $S$ with an ideal
$W(S)$ such that
\begin{enumerate}
\item[(i)] resolving $S$ is equivalent to resolving $W(S)$, 
\item[(ii)] the traces $W(S)|_H$ are locally 
  isomorphic for all
hypersurfaces of maximal contact through $s\in S$
(here ``locally'' is meant in the analytic or \'etale topology), and 
\item[(iii)] the resolution of  $W(S)|_H$
tells us how to resolve $W(S)$.
\end{enumerate}
The local ambiguity is thus removed from the process,
and there is no longer a patching problem.

\medskip

 {\it Problem}
\ref{probls.with.meth.say}.6.
 At Steps 2 and 3 in (\ref{deg1Econe.w.meth}),
 the choices we make are not
 canonical. For instance, in Step 2 we could have blown up the central curve  
with multiplicity 2 first,
to complete the resolution in just one step. Even if we do Step 2 as above,
in general there are many curves to blow-up in Step 3, and the
order of blow-ups matters. (In $\a^3$, one can blow up
two intersecting smooth curves
 in either order, and the resulting 3-folds are
not isomorphic.)

This problem, too, remains unsolved. We
 make a choice, and it is good
enough that the resolutions we get  commute with any smooth
morphism.  Thus we get a resolution that one can 
call functorial. I would not call it a
canonical resolution, since even in the framework of this proof
other equally functorial choices  are possible.

This is very much connected with the lack of minimal
resolutions.
\end{say}

Next we see how Problems (\ref{probls.with.meth.say}.2--5) can be approached 
for hypersurfaces using Weierstrass polynomials.
As was the case with curves and surfaces,
this example motivates the whole proof.
(To be fair, this example provides much better guidance
with hindsight. One might argue that the whole history
of resolution by smooth blow-ups is but an ever-improving
 understanding of this single example.
It has taken a long time to sort out how to generalize
various aspects of it, and it is by no means certain that
we have learned all the right lessons.)

\begin{exmp}\label{maxcont.key.exmp}
\index{Resolution!of hypersurfaces using Weierstrass  equation}
\index{Weierstrass!equation and resolution} 
Let $X\subset \c^{n+1}$ be a hypersurface.
Pick a point $0\in X$, where $\mult_0X=m$.
Choose suitable local coordinates $x_1,\dots,x_{n},z$, and apply
the Weierstrass preparation theorem %%(\ref{gen.weierstrass.thm})
 to get 
(in an analytic neighborhood) an 
 equation  of the form
$$
z^m+a_1({\mathbf x})z^{m-1}+\cdots+a_m({\mathbf x})=0
$$
for $X$.
We can kill the $z^{m-1}$ term by a substitution
$z= y-\frac1{m}a_{1}({\mathbf x})$ to get
another local equation
$$
f:=y^m+b_2({\mathbf x})y^{m-2}+\cdots+b_m({\mathbf x})=0.
\eqno{(\ref{maxcont.key.exmp}.1)}
$$
Here $\mult_0b_i\geq i$ since $\mult_0X=m$.

Let us blow up the point $0$ to get $\pi:B_0X\to X$, and consider the chart 
$x'_i=x_i/x_n, x'_n=x_n, y'=y/x_n$. We get an equation
for $B_0X$
$$
F:=
(y')^m+(x'_n)^{-2}b_2(x'_1x'_n,\dots,x'_n)(y')^{m-2}+\cdots+
(x'_n)^{-m}b_m(x'_1x'_n,\dots,x'_n).
\eqno{(\ref{maxcont.key.exmp}.2)}
$$

Where are the  points of multiplicity $\geq m$ on $B_0X$?
Locally we can  view $B_0X$ as a hypersurface
in $\c^{n+1}$ given by the equation
$F({\mathbf x'},y')=0$, and a point $p$ has multiplicity
 $\geq m$
iff all the $(m-1)$st partials of $F$ vanish.  
First of all, we get that
$$
\frac{\partial^{m-1}F}{\partial {y'}^{m-1}}=m!\cdot y'
\qtq{vanishes at $p$.}
\eqno{(\ref{maxcont.key.exmp}.3)}
$$

This means
that all points of
multiplicity $\geq m$ on $B_0X$ are on the birational transform
of the hyperplane $(y=0)$. Since the new equation (\ref{maxcont.key.exmp}.2)
has the same form as the original  (\ref{maxcont.key.exmp}.1),
 the conclusion continues to hold after further blow-ups,
solving   (\ref{probls.with.meth.say}.3):

\medskip

 {\it Claim}
\ref{maxcont.key.exmp}.4.  After a sequence of blow-ups
at points of multiplicity $\geq m$
$$
\Pi: X_r=B_{p_{r-1}}X_{r-1}\to X_{r-1}=B_{p_{r-2}}X_{r-2}\to\cdots
\to X_1=B_{p_0}X\to X,
$$
all points of
multiplicity $\geq m$ on $X_r$ are on the birational transform
of the hyperplane $H:=(y=0)$,
and all points of $X_r$ have  multiplicity $\leq m$.
\medskip

This property of the hyperplane $(y=0)$ will be 
encapsulated by the concept of
{\it   hypersurface of maximal contact}.
\index{Maximal contact!using Weierstrass equation}
\medskip

In order to determine the location of points of
multiplicity $m$, we need to look at 
all the other $(m-1)$st partials of $F$
restricted to $(y'=0)$. These
can be written as
$$
\frac{\partial^{m-1}F}{\partial {\mathbf x'}^{i-1} \partial {y'}^{m-i}}
|_{(y'=0)}=
(m-i)!\cdot 
\frac{\partial^{i-1}\bigl((x'_n)^{-i}b_i(x'_1x'_n,\dots,x'_n)\bigr)}
{\partial {\mathbf x'}^{i-1}}.
\eqno{(\ref{maxcont.key.exmp}.5)}
$$

Thus  we can actually read off from $H=(y=0)$
which points of $B_0X$ have multiplicity $m$.
 For this, however,  we need
not only the restriction $f|_H=b_m({\mathbf x})$
but   all the  
other coefficients  $b_i({\mathbf x})$ as well.

There is one further twist. The usual rule
for transforming a polynomial under a blow-up is
$$
b(x_1,\dots,x_n)\mapsto
(x'_n)^{-\mult_0b}b(x'_1x'_n,\dots,x'_n),
$$
but instead we use the rule
$$
b_i(x_1,\dots,x_n)\mapsto
(x'_n)^{-i}b_i(x'_1x'_n,\dots,x'_n).
$$
That is, we 
``pretend'' that $b_i$ has multiplicity $i$ at the origin.
To handle this, we introduce the notion of a marked
function $(g,m)$ and define 
 the birational transform of
a marked function $(g,m)$ to be
$$
\pi^{-1}_*\bigl(g(x_1,\dots,x_n),m\bigr):=
\bigl((x'_n)^{-m}g(x'_1x'_n,\dots,x'_n),m\bigr).
\eqno{(\ref{maxcont.key.exmp}.6)}
$$
{\it Warning}. If we change coordinates, the right-hand side of 
(\ref{maxcont.key.exmp}.6)  
changes by a unit. Thus the ideal $(\pi^{-1}_*(g,m))$
is  well defined but not $\pi^{-1}_*(g,m)$ itself.
  Fortunately, this does not lead to
any problems.

By induction we define
$\Pi^{-1}_*(g,m)$, where $\Pi$ is a sequence of blow-ups as
in (\ref{maxcont.key.exmp}.4).

This leads to a solution of Problems (\ref{probls.with.meth.say}.2)
and (\ref{probls.with.meth.say}.4).
\medskip

 {\it Claim}
\ref{maxcont.key.exmp}.7. After a sequence of blow-ups
at points of multiplicity $\geq m$,
$$
\Pi: X_r=B_{p_{r-1}}X_{r-1}\to X_{r-1}=B_{p_{r-2}}X_{r-2}\to\cdots
\to X_1=B_{p_0}X\to X,
$$
a  point $p\in X_r$ has 
multiplicity $< m$ on $X_r$ iff 
\begin{enumerate}
\item[(i)] either $p\not\in H_r$, the birational transform of $H$,
\item[(ii)] or 
there is an index $i=i(p)$ such that
$$
\mult_p(\Pi|_{H_r})^{-1}_*\bigl(b_i({\mathbf x}),i\bigr)<i.
$$
\end{enumerate}
\medskip

A further observation is  that we can obtain the $b_i({\mathbf x})$
 from the derivatives of $f$:
$$
 b_i({\mathbf x})=\frac{1}{(m-i)!}\cdot 
\frac{\partial^{m-i}f}{\partial y^{m-i}}({\mathbf x},y)|_H.
$$
Thus (\ref{maxcont.key.exmp}.7)
can be restated in a more invariant-looking but
also vaguer form.
\medskip

 {\it Principle}
\ref{maxcont.key.exmp}.8. Multiplicity reduction for
the $n+1$-variable function
$f({\mathbf x},y)$
 is equivalent to  multiplicity  reduction 
 for certain $n$-variable functions
constructed from  the partial derivatives of $f$
with suitable  markings.
\medskip

\ref{maxcont.key.exmp}.9. Until now  we have completely ignored 
that everything we do  depends on the
initial choice of the coordinate system
$(x_1,\dots,x_n,z)$.
The fact that in (\ref{maxcont.key.exmp}.7--8)
we get  equivalences suggests that the
 choice of the coordinate system
should not matter much. The problem, however, remains:
 in globalizing the local resolutions
constructed above, we have to choose local resolutions
out of the many possibilities and hope that the different
local choices patch together.

This has been a surprisingly serious obstacle.
\end{exmp}

\section{Statement of the main theorems}\label{st.of.main.res}

So far we have been concentrating on resolution of singularities,
but now we switch our focus, and instead of dealing with
singular varieties, we consider ideal sheaves on {\em smooth} varieties.
Given an ideal sheaf $I$ on a smooth variety $X$, our first aim
is to write down a birational morphism $g:X'\to X$
such that $X'$ is smooth and the pulled-back ideal sheaf $g^*I$ is
locally principal. This is called the
principalization of $I$.

\begin{notation} Let $g:Y\to X$ be a morphism of schemes
and $I\subset \o_X$  an ideal sheaf. 
I will be sloppy and use $g^*I$ to denote
the {\it inverse image ideal sheaf} of $I$. 
\index{Ideal sheaf!inverse image of}
\index{Inverse image ideal sheaf, $g^*I$}
This is the ideal sheaf
generated by the pull-backs of local sections of $I$.
(It is denoted by $g^{-1}I\cdot \o_Y$ or by $I\cdot \o_Y$
in \cite[Sec.II.7]{hartsh}.) 
 
We should be mindful that $g^*I$ (as an  inverse image ideal sheaf)
 may differ from  the usual sheaf-theoretic
pull-back, also commonly 
denoted by  $g^*I$; see \cite[II.7.12.2]{hartsh}.
This can happen even if $X,Y$ are both smooth.

For the rest of the chapter, we use only
 inverse image ideal sheaves, so hopefully this should not lead to  any
confusion.
\end{notation}

It is easy to see that resolution of singularities
implies principalization. Indeed, let $X_1:=B_IX$ be
the blow-up of $I$ with projection $\pi:X_1\to X$.
Then $\pi^*I$ is locally principal (cf.\ \cite[II.7.13]{hartsh}). Thus if
$h:X'\to X_1$ is any resolution, then
$\pi\circ h:X'\to X$ is a principalization of $I$.

Our aim, however, is to derive
resolution theorems from principalization results.
Given a singular variety $Z$, choose an embedding
of $Z$ into a smooth variety $X$, and let $I_Z\subset \o_X$
be its ideal sheaf.
(For $Z$  quasi-projective, we can just take 
any embedding $Z\into \p^N$ 
into a projective space, but in general
such an embedding may not exist; see (\ref{no.smooth.emb.exmp}).)  Then
we turn a principalization  of the ideal
sheaf $I_Z$ into a resolution of $Z$.

In this section we state four, increasingly stronger
versions of principalization and derive
from them various resolution theorems.
The rest of the chapter is then devoted to
proving these principalization theorems.

\begin{say}[Note on terminology] 
Principal ideals are much simpler than arbitrary ideals,
but they can still be rather complicated since they capture
all the intricacies of hypersurface singularities. 

An ideal sheaf $I$ on a smooth scheme $X$ is called (locally) 
{\it  monomial}\index{Ideal sheaf!monomial}
if the following equivalent conditions hold.
\begin{enumerate}
\item For every $x\in X$ there are local coordinates
$z_i$ and natural numbers $c_i$ such that
$I\cdot \o_{x,X}=\prod_i z_i^{c_i}\cdot \o_{x,X}$.
\item $I$ is the ideal sheaf of a simple normal
crossing divisor (\ref{nc.divs.defn}).
\end{enumerate}

I would like to call a
 birational morphism $g:X'\to X$,
such that $X'$ is smooth and  $g^*I$ is
 monomial,   a  resolution of $I$.

However, for many people, the phrase
``resolution of an ideal sheaf'' brings to mind a long  exact
sequence
$$
\cdots \to E_2\to E_1\to I\to 0,
$$
where the $E_i$ are locally free sheaves. 
This has nothing to do with resolution of singularities.
Thus, rather reluctantly,
I follow convention and  talk about 
{\it principalization}  or {\it monomialization} of an ideal sheaf $I$.
\index{Principalization!of ideal sheaves}
\index{Ideal sheaf!principalization}
\index{Monomialization!of ideal sheaves}
\index{Ideal sheaf!monomialization}
\end{say}

We start with the simplest version of  principalization
(\ref{princ.thm.1})
and its first consequence, the resolution of indeterminacies
of rational maps  
(\ref{res.ind.thm.1}). Then we consider a stronger version
of principalization
(\ref{princ.thm.2}),  which
implies resolution of singularities
(\ref{res.sings.thm.1}). 
Monomialization of ideal sheaves is given in (\ref{funct.princ.thm}),
which implies strong, functorial resolution for
quasi-projective varieties  (\ref{res.sings.main.thm}). 
The proof of the strongest
variant of monomialization (\ref{funct.princ.thm.IV})
occupies the rest of the chapter. 
At the end of the section we observe that the
 functorial properties proved in
(\ref{funct.princ.thm.IV}) imply that 
the  monomialization and resolution theorems automatically extend
to algebraic and analytic spaces; see
(\ref{algsp.res.say}) and  (\ref{analsp.res.say}).

\begin{thm}[Principalization, I]\label{princ.thm.1}
\index{Principalization!of ideal sheaves}
\index{Ideal sheaf!principalization}
  Let $X$ be a smooth variety
over a field of characteristic zero
and $I\subset \o_X$ a nonzero  ideal sheaf.
Then there is a smooth variety $X'$ and a
birational and projective morphism $f:X'\to X$ such that
$f^*I\subset\o_{X'}$ is a locally principal ideal sheaf.
\end{thm}

\begin{cor}[Elimination of indeterminacies]\label{res.ind.thm.1}
\index{Elimination!of  indeterminacies}
\index{Indeterminacy!elimination of}
  Let $X$ be a smooth variety over a field of characteristic zero
and $g:X\map \p$ a rational map
to some projective space.
Then there is a smooth variety $X'$ and a
birational and projective morphism $f:X'\to X$ such that
the composite $g\circ f:X'\to \p$ is a morphism.
\end{cor}

Proof. Since $\p$ is projective and $X$ is normal,
there is a subset $Z\subset X$ of codimension $\geq 2$
such that $g:X\setminus Z\to \p$ is a morphism.
Thus $g^*\o_{\p}(1)$ is a line bundle on
$X\setminus Z$. 
Since $X$ is smooth, it extends to a line bundle on
$X$; denote it by $L$. Let $J\subset L$ be the subsheaf generated
by  $g^*H^0(\p,\o_{\p}(1))$.  Then $I:=J\otimes L^{-1}$ is an ideal sheaf,
and so by (\ref{princ.thm.1})
there is  a
projective morphism $f:X'\to X$ such that
$f^*I\subset\o_{X'}$ is a locally principal ideal sheaf.

Thus the  global sections 
$$
(g\circ f)^*H^0\bigl(\p,\o_{\p}(1)\bigr)\subset H^0(X', f^*L)
$$
generate the locally free sheaf
$L':=f^*I\otimes f^*L$. Therefore,
$g\circ f:X'\to \p$ is a morphism given by
the nowhere-vanishing subspace of global sections
 $$
(g\circ f)^*H^0\bigl(\p,\o_{\p}(1)\bigr)\subset H^0(X', L').\qed
$$

\begin{notation}[Blow-ups]\label{bu.defn}
Let $X$ be a scheme and $Z\subset X$ a closed subscheme.
Let $\pi=\pi_{Z,X}:B_ZX\to X$ denote the {\it blow-up}\index{Blow-up}
of $Z$ in $X$; see \cite[II.4]{shaf} or  \cite[Sec.II.7]{hartsh}.
Although resolution by definition involves singular schemes $X$,
we will almost always study the case 
 where $X$ and $Z$ are both smooth,
called a {\it smooth blow-up}.\index{Blow-up!smooth}
The {\it exceptional divisor}\index{Exceptional!divisor}
\index{Divisor!exceptional} of a  blow-up
is  $F:=\pi_{Z,X}^{-1}(Z)\subset B_ZX$. If $\pi_{Z,X}$ is a smooth blow-up,
then $F$ and $B_ZX$ are both smooth.
\end{notation}

\begin{warning}[Trivial and empty blow-ups] \label{triv.bu.warn}
A  blow-up is called {\it trivial}\index{Blow-up!trivial}%
\index{Trivial blow-up}%
if $Z$ is a  Cartier divisor in $X$. In  these cases $\pi_{Z,X}:B_ZX\to X$
is an isomorphism.
We also allow   the possibility
$Z=\emptyset$, called the {\it empty}\index{Blow-up!empty}%
\index{Empty!blow-up} blow-up.

We have to deal with trivial blow-ups to make induction work
since the blow-up of a codimension 2 smooth subvariety
$Z^{n-2}\subset X^n$ corresponds to a trivial blow-up
on a smooth hypersurface $Z^{n-2}\subset H^{n-1}\subset X^n$.

Two peculiarities of trivial blow-ups cause trouble.
\begin{enumerate}
\item For a  nontrivial smooth blow-up $\pi:B_ZX\to X$,
the morphism $\pi$  determines
the center $Z$, but this fails for a  trivial blow-up.
One usually thinks of $\pi$ as {\em the} blow-up, 
hiding the dependence on $Z$. By contrast, we always think of
a smooth blow-up as having a  specified center.
\item The   exceptional divisor
of a trivial blow-up $\pi_{Z,X}:B_ZX\to X$
is $F=Z\subset X$. \index{Exceptional!divisor of a trivial blow-up}
This is, unfortunately, at variance with the usual definition
of exceptional set/divisor 
(see \cite[Sec.II.4.4]{shaf} or (\ref{bir.tr.divs.def})), but it is the right
concept for blow-ups.
\end{enumerate}
These are both minor inconveniences, but they could lead
to confusion.

Empty blow-ups naturally occur when we  restrict a blow-up sequence to an
open subset  $U\subset X$ and the center of the blow-up is disjoint from
$U$. We will exclude empty blow-ups from the final blow-up sequences,
but  we have to keep them in mind since they mess up the numbering
of the blow-up sequences.
\end{warning}

\begin{thm}[Principalization, II]\label{princ.thm.2}
\index{Principalization!strong form}
\index{Ideal sheaf!principalization, strong form}
  Let $X$ be a smooth variety over a field of characteristic zero
and $I\subset \o_X$ a nonzero ideal sheaf.
Then there is a smooth variety $X'$ and a
birational and projective morphism $f:X'\to X$ such that
\begin{enumerate}
\item 
$f^*I\subset\o_{X'}$ is a locally principal ideal sheaf, 
\item $f:X'\to X$ is an isomorphism over $X\setminus \cosupp I$, 
where $\cosupp I$ (or $\supp (\o_X/I)$) is the cosupport of $I$,
\index{Ideal sheaf!cosupport of}\index{Cosu@$\cosupp I$, cosupport}
 and
\item $f$ is a composite of smooth blow-ups
$$
f:X'=X_r\stackrel{\pi_{r-1}}{\longrightarrow} X_{r-1}
\stackrel{\pi_{r-2}}{\longrightarrow}\cdots
\stackrel{\pi_{1}}{\longrightarrow} X_1
\stackrel{\pi_{0}}{\longrightarrow} X_0=X.
$$
\end{enumerate}
\end{thm}

This form of principalization implies resolution of
singularities, seemingly by accident.
(In practice, one can follow the steps of a
principalization method and see how resolution happens,
though this is not always easy.)

\begin{cor}[Resolution of singularities, I]\label{res.sings.thm.1}
\index{Resolution!standard form} 
 Let $X$ be a  quasi-projective variety.
Then there is a smooth variety $X'$ and a
birational and projective morphism $g:X'\to X$.
\end{cor}

Proof. Choose an embedding of  $X$ into 
a smooth variety $P$ such that $N\geq \dim X+2$.
(For instance, $P=
\p^N$ works for all $N\gg \dim X$.)
Let $\bar X\subset P$ denote the closure and $I\subset \o_{P}$
its ideal sheaf. Let $\eta_X\in X\subset P$ be the generic point of $X$.

By (\ref{princ.thm.2}), there is a sequence of smooth blow-ups
$$
\Pi:P'=P_r\stackrel{\pi_{r-1}}{\longrightarrow}
P_{r-1}\stackrel{\pi_{r-2}}{\longrightarrow}\cdots 
P_1\stackrel{\pi_0}{\longrightarrow}
P_0=P
$$
such that $\Pi^*I$ is locally principal.

Since $X$ has codimension $\geq 2$, its ideal sheaf 
 $I$ is not locally principal at $\eta_X$, and therefore, 
 some blow-up center must contain $\eta_X$.
Thus, there is a unique $j$ such that
$\pi_{0}\cdots \pi_{j-1}:P_j\to P$
is a local isomorphism around $\eta_X$ but
$\pi_j:P_{j+1}\to P_j$ is a blow-up with center $Z_j\subset P_j$
such that $\eta_X\in Z_j$.

 By (\ref{princ.thm.2}.2), 
$\pi_{0}\cdots \pi_{j-1} (Z_j)\subset \bar X$, and this implies that
$\eta_X$ is the generic point of $Z_j$. 
Thus
$$
g:=\pi_{0}\cdots \pi_{j-1}:Z_j\to \bar X
$$
is birational.

 $Z_j$ is smooth since we blow it up, 
and by (\ref{princ.thm.2}.3) we only blow up smooth subvarieties.
Therefore $g$ is 
 a resolution of singularities of $\bar X$.
 Set $X':=g^{-1}(X)\subset Z_j$.
Then $g:X'\to X$ is a resolution of singularities
of $X$.
\qed

\begin{warning} The resolution
    $g:X'\to X$ constructed in (\ref{res.sings.thm.1})
need not be  a composite of {\em smooth} blow-ups.
Indeed, the process exhibits $g$ as the composite of
blow-ups whose centers are obtained by intersecting
the smooth centers $Z_i$ with the birational transforms of $X$.
Such intersections may be singular.
See (\ref{nsbu.exmp}) for a concrete example.
\end{warning}

We also need a form of  resolution that keeps track of 
a suitable simple normal crossing divisor.
This feature is very useful  in applications and in the inductive proof.

\begin{defn}\label{nc.divs.defn}
Let $X$ be a smooth variety and $E=\sum E^i$
 a {\it simple normal crossing divisor} on $X$.
\index{Simple normal crossing divisor}%
\index{Normal crossing!simple}%
\index{Divisor!simple normal crossing}%
\index{Snc, simple normal crossing}%
This means that each $E^i$ is smooth, and for each point $x\in X$
one can choose local  coordinates 
$z_1,\dots,z_n\in  m_x$ in the maximal ideal of the
 local ring $\o_{x,X}$
 such that for each $i$
\begin{enumerate}
\item either $x\not\in E^i$, or
\item $E^i= (z_{c(i)}=0)$ in a neighborhood of $x$ for some $c(i)$, and
\item $c(i)\neq c(i')$ if $i\neq i'$.
\end{enumerate}

A subvariety $Z\subset X$ has {\it simple normal crossings} with $E$
if
one can choose  $z_1,\dots,z_n$ as above 
such that in addition
\begin{enumerate}\setcounter{enumi}{3}
\item $Z= (z_{j_1}=\cdots=z_{j_s}=0)$ for some 
$j_1,\dots,j_s$, again in some open neighborhood of $x$.
\end{enumerate}
In particular, $Z$ is smooth,
and some of the $E^i$ are allowed to contain $Z$.

If  $E$ does not contain $Z$, then
$E|_Z$ is again a simple normal crossing divisor on $Z$.
\end{defn}

\begin{defn}\label{bir.tr.divs.def}
Let $g:X'\to X$ be a birational morphism.
Its {\it exceptional set} is 
the set of points $x'\in X'$ such that $g$ is not a local isomorphism at $x'$.
It is denoted by $\ex(g)$.
\index{Exceptional!set $\ex(f)$}\index{Exf@$\ex(f)$, exceptional set}%
If $X$ is smooth, then $\ex(g)$ is a divisor \cite[II.4.4]{shaf}.
Let $$
\Pi:X'=X_r\stackrel{\pi_{r-1}}{\longrightarrow} X_{r-1}
\stackrel{\pi_{r-2}}{\longrightarrow}\cdots
\stackrel{\pi_{1}}{\longrightarrow} X_1
\stackrel{\pi_{0}}{\longrightarrow} X_0=X
$$
be a sequence of smooth blow-ups with centers $Z_i\subset X_i$.
Define the {\it total exceptional set} 
\index{Exceptional!set, total $\extot(f)$}%
\index{Extotf@$\extot(f)$, exceptional set}%
\index{Total!exceptional set}%
to be
$$
\extot(\Pi):=
\bigcup_{i=0}^{r-1}\bigl(\pi_i\circ\cdots\circ \pi_{r-1})^{-1}(Z_i).
$$
If all the blow-ups are nontrivial, then $\ex(\Pi)=\extot(\Pi)$.

 Let $E$
 be a simple normal crossing divisor on $X$. We say that the centers 
$Z_i$
have {\it simple normal crossings with $E$} if
\index{Simple normal crossing!with a divisor}
 each
 blow-up  center
$Z_i\subset X_i$ 
has simple normal crossings (\ref{nc.divs.defn}) with 
$$
(\pi_0\cdots \pi_{i-1})^{-1}_*(E)+\extot(\pi_0\cdots \pi_{i-1}).
$$
If this holds, then   
$$
\Pi^{-1}_{\rm tot}(E):=\Pi^{-1}_*(E)+\extot(\Pi)
$$
is a  simple normal crossing divisor, called the 
{\it total transform}\index{Total!transform of a divisor}%
\index{Transform!total, for a divisor}
 of $E$. (A refinement for divisors with ordered index set will
be introduced in (\ref{test.bu.defn}).)
\end{defn}

We can now strengthen the 
theorem on principalization of  ideal sheaves.

\begin{thm}[Principalization, III]\label{funct.princ.thm}
\index{Principalization!strong form}
\index{Ideal sheaf!principalization, strong form}
Let  $X$  be a smooth variety over
  a field of characteristic zero,
 $I\subset \o_X$ a nonzero ideal sheaf
and $E$ a simple normal crossing divisor on $X$.
Then there is a 
sequence of  smooth blow-ups
$$
\Pi:R_{I,E}(X):=X_r\stackrel{\pi_{r-1}}{\longrightarrow} X_{r-1}
\stackrel{\pi_{r-2}}{\longrightarrow}\cdots
\stackrel{\pi_{1}}{\longrightarrow} X_1
\stackrel{\pi_{0}}{\longrightarrow} X_0=X
$$
whose centers have simple normal crossing with $E$
such that
\begin{enumerate}
\item 
$\Pi^*I\subset\o_{R_{I,E}(X)}$ is the ideal sheaf of a simple normal
crossing divisor, and 
\item $\Pi:R_{I,E}(X)\to X$ is functorial on smooth morphisms
(\ref{functres.res.say}).
\end{enumerate}
\end{thm}

Note that since $\Pi$ is a composite of smooth blow-ups, 
 $R_{I,E}(X)$ is smooth and $\Pi:R_{I,E}(X)\to X$
is birational and projective.

As a consequence we get  strong 
 resolution of singularities
for quasi-projective  schemes 
over  a field of characteristic zero.

\begin{thm}[Resolution of singularities, II]\label{res.sings.main.thm}
\index{Resolution!strong form}
Let $X$ be a quasi-projective  variety
over a   field of characteristic zero. 
Then there is a birational and projective morphism
$\Pi:R(X)\to X$ such that
\begin{enumerate}
\item $R(X)$ is smooth,
\item $\Pi:R(X)\to X$ is an isomorphism over 
the smooth locus $X^{ns}$, and
\item $\Pi^{-1}(\sing X)$ is a
divisor with simple normal crossing.
\end{enumerate}
\end{thm}

Proof. We have already seen in
(\ref{res.sings.thm.1}) that given a (locally closed)  embedding
$i:X\into P$ we get a resolution  $R(X)\to X$ from the
principalization of the ideal sheaf $I$ of the closure of $i(X)$.
We need to  check that applying (\ref{funct.princ.thm})
to $(P,I, \emptyset)$
gives a strong resolution of $X$.
(We do not claim that $R(X)\to X$ is independent of the
embedding $i:X\into P$. This will have to wait until after the stronger
principalization theorem (\ref{funct.princ.thm.IV}).)

As in 
 the proof of (\ref{res.sings.thm.1}), 
 there is a sequence of smooth blow-ups
$$
\Pi:P'=P_r\stackrel{\pi_{r-1}}{\longrightarrow}
P_{r-1}\stackrel{\pi_{r-2}}{\longrightarrow}\cdots 
P_1\stackrel{\pi_0}{\longrightarrow}
P_0=P
$$
such that $\Pi^*I$ is locally principal. Moreover, there is a first blow-up
$$
\pi_j:P_{j+1}\to P_j\qtq{with center $Z_j\subset P_j$}
$$
such that   
$g:=\pi_{0}\cdots \pi_{j-1}|_{Z_j}: Z_j\to \bar X$
is birational. 
We claim that $g:Z_j\to \bar X$ is a strong resolution
of $\bar X$, and  hence
$g:g^{-1}(X)\to X$ is a strong resolution
of $ X$.

First, we prove that $g$ is an isomorphism over $\bar X^{ns}$.
As in (\ref{functres.res.say}.2), this follows from the functoriality condition
(\ref{funct.princ.thm}.2). Note, however, that 
(\ref{funct.princ.thm}.2) asserts functoriality for
$\Pi:P_r\to P$ but not for the intermediate
maps $P_j\to P$. Thus a little extra work is needed.
(This seems like a small technical point, but actually
it has been the source of serious troubles.
The notion of blow-up sequence functors (\ref{busf.defnition}) is 
designed
to deal with it.)

Let $F'_{j+1}\subset P_r$ denote the birational transform
of  $F_{j+1}\subset P_{j+1}$, the exceptional divisor of
$\pi_j$. Since $F'_{j+1}\subset \extot(\Pi)$, it is a smooth
divisor and so $\Pi|_{F'_{j+1}}:F'_{j+1}\to \bar X$
is generically smooth. Thus there is a  smooth
point $x\in X$ such that $\Pi|_{F'_{j+1}}$ is smooth over $x$.

For any other smooth point $x'\in X$, the embeddings
$$
(x\in X\into P)\qtq{and} (x'\in X\into P)
$$
have isomorphic \'etale neighborhoods. Thus by
(\ref{funct.princ.thm}.2), $\Pi|_{F'_{j+1}}$ is also smooth over $x'$.
We can factor
$$
\Pi|_{F'_{j+1}}:F'_{j+1}\to Z_j\stackrel{g}{\to} \bar X.
$$
Thus $g:Z_j\to \bar X$ is smooth over every smooth point of $\bar X$.
It is also birational, and thus $g$ is an isomorphism
over $\bar X^{ns}$.

Since  $Z_j$ is smooth, $g$ is not an isomorphism
over any point of $\sing \bar X$, and thus
$$
g^{-1}(\sing \bar X)= Z_j\cap \extot (\pi_{0}\cdots \pi_{j-1}),
$$
where $\extot$ denotes the total exceptional divisor
(\ref{bir.tr.divs.def}).
Observe that in  (\ref{funct.princ.thm}) we can only blow-up
$Z_j$ if it has simple normal crossings with
$\extot(\pi_{0}\cdots \pi_{j-1})$; hence
$$
g^{-1}(\sing \bar X)=Z_j\cap \extot(\pi_{0}\cdots \pi_{j-1})
$$
is a simple normal crossing divisor on $Z_j$. \qed

\begin{rem}\label{nonproj.res.loc.glob.say}
The proof of   the implication 
(\ref{funct.princ.thm}) $\Rightarrow$ (\ref{res.sings.main.thm})  
also  works for any scheme that can be
embedded into a smooth variety. We see in (\ref{no.smooth.emb.exmp})
 that 
not all schemes  can be
embedded into a smooth scheme, so in general one has to
proceed differently. 
It is worthwhile to contemplate further the
local nature of resolutions and its consequences.

Let $X$ be a  scheme of finite type and
$X=\cup U_i$ an affine cover.  For each $U_i$ 
(\ref{res.sings.main.thm}) gives   a  resolution $R(U_i)\to U_i$,
and we would like to patch these together 
to $R(X)\to X$.

First, we need to show that $R(U_i)$ is well defined;
that is, it does not depend on the
embedding $i:U_i\into P$ chosen in the proof of 
(\ref{res.sings.main.thm}).

Second, we need to show that $R(U_i)$ and $R(U_j)$ agree
over the intersection $U_i\cap U_j$.

If these hold, then 
 the $R(U_i)$ patch together into a 
resolution $R(X)\to X$, but
there is one problem. 
$R(X)\to X$ is {\em locally} projective, but it
may not be {\em globally} projective.
The following is an example of this type.
\medskip

{\it Example} \ref{nonproj.res.loc.glob.say}.1. 
\index{Resolution!nonprojective example}%
Let $X$ be a smooth 3-fold and  $C_1,C_2$ a pair of ir\-re\-ducible curves,
intersecting  at two points $p_1,p_2$.
Assume, furthermore, that $C_i$ is smooth away from $p_i$,
where it has a cusp whose tangent plane is transversal
to the other curve. Let $I\subset \o_X$ be the ideal sheaf
of $C_1\cup C_2$.

On $U_1=X\setminus \{p_1\}$, the curve $C_1$ is smooth;
we can  blow it up first.
The birational transform of $C_2$ becomes smooth,
and we can blow it up next to get $Y_1\to U_1$.
Over $U_2=X\setminus \{p_2\}$ 
we would work in the other order.
Over $U_1\cap U_2$ we get the same thing, and 
thus $Y_1$ and $Y_2$ glue together
to a variety $Y$ such that $Y\to X$ is proper and  locally projective
but not globally projective.
\medskip

We see that the gluing problem comes from the circumstance that
the birational map $Y_1\cap Y_2\to U_1\cap U_2$ is the
blow-up of two disjoint curves, and we do not know which one to blow up first.

For a sensible resolution algorithm there is only one choice:
we have to blow them up at the same time. Thus in the above example,
the ``correct'' method is to blow up the points $p_1,p_2$ first.
The curves $C_1,C_2$ become smooth and disjoint, and then both can
be blown up. (More blow-ups are needed if we want to have only
simple normal crossings.)
\end{rem}

These problems can be avoided if we make
(\ref{funct.princ.thm})  sharper. A key point is to prove
 functoriality conditions
not only for the end result $R_{I,E}(X)$
 but for all intermediate steps, including the center
of each blow-up.

\begin{defn}[Blow-up sequences]\label{bu.seq.notation.defn}
Let $X$ be a scheme.
A {\it  blow-up sequence}\index{Blow-up!sequence}
 of length $r$ starting with $X$ is a
chain of morphisms
$$
\begin{array}{ccccccl}
\Pi:X_r&\stackrel{\pi_{r-1}}{\longrightarrow} &X_{r-1}&
\stackrel{\pi_{r-2}}{\longrightarrow}\ \cdots\  
\stackrel{\pi_{1}}{\longrightarrow}& X_1&
\stackrel{\pi_{0}}{\longrightarrow} &X_0=X,\\
&& \cup && \cup && \, \cup \\
&& Z_{r-1} &\cdots & Z_1 && Z_0
\end{array}
\eqno{(\ref{bu.seq.notation.defn}.1)}
$$
where each $\pi_i=\pi_{Z_i,X_i}:X_{i+1}\to X_i$ is a  blow-up with center
$Z_i\subset X_i$ and exceptional divisor $F_{i+1}\subset X_{i+1}$.
Set
$$
\Pi_{ij}:=\pi_{j}\circ\cdots\circ\pi_{i-1}:X_i\to X_j
\qtq{and} \Pi_i:=\Pi_{i0}:X_i\to X_0.
$$
We say that (\ref{bu.seq.notation.defn}.1) is a
{\it smooth blow-up sequence}\index{Blow-up!sequence, smooth}
if each $\pi_i:X_{i+1}\to X_i$ is a  smooth blow-up.

We allow trivial and  empty blow ups (\ref{triv.bu.warn}).

For the rest of the chapter, $\pi$ always denotes a  blow-up,
$\Pi_{ij}$ a composite of blow-ups and $\Pi$ the composite
of all blow-ups in a  blow-up sequence
(whose length we frequently leave unspecified).
We usually drop the centers $Z_i$ from the notation,
 to avoid cluttering up the diagrams.
\end{defn}

\begin{defn}[Transforming blow-up sequences]\label{bu.seq.functs}
There are three basic ways to transform blow-up sequences from one
scheme to another. Let ${\mathbf B}:=$ 
$$
\begin{array}{ccccccl}
\Pi:X_r&\stackrel{\pi_{r-1}}{\longrightarrow} &X_{r-1}&
\stackrel{\pi_{r-2}}{\longrightarrow}\ \cdots\  
\stackrel{\pi_{1}}{\longrightarrow}& X_1&
\stackrel{\pi_{0}}{\longrightarrow} &X_0=X\\
&& \cup && \cup && \, \cup \\
&& Z_{r-1} &\cdots & Z_1 && Z_0
\end{array}
$$
be a blow-up sequence starting with $X$.
\medskip

\ref{bu.seq.functs}.1. 
For a  smooth morphism $h:Y\to X$ define the
{\it pull-back}\index{Pull-back!of a blow-up sequence}%
\index{Blow-up!sequence, pull-back}
$h^*{\mathbf B}$ to be the  blow-up sequence
$$
\begin{array}{ccccccl}
h^*\Pi:X_r\times_XY
&\!\!\!\stackrel{h^*\pi_{r-1}}{\longrightarrow}\!\!\! &X_{r-1}\times_XY&
\cdots & X_1\times_XY&
\!\!\!\stackrel{h^*\pi_{0}}{\longrightarrow}\!\!\! &X_0\times_XY=Y.\\
&& \cup && \cup && \ \ \ \ \ \cup \\
&& Z_{r-1}\times_XY &\cdots & Z_1\times_XY && Z_0\times_XY
\end{array}
$$
If ${\mathbf B}$ is a smooth blow-up sequence then so is $h^*{\mathbf B}$.
If $h$ is surjective then $h^*{\mathbf B}$ determines 
${\mathbf B}$ uniquely. However, if $h$ is not surjective, then
$h^*{\mathbf B}$ may contain some empty blow-ups, and
we lose information about the centers living above $X\setminus h(Y)$.
\medskip

\ref{bu.seq.functs}.2. 
Let $X$ be a  scheme and $j:S\into X$ a  closed subscheme.
Given a   blow-up sequence ${\mathbf B}$
 starting with $X$ as above, define its {\it restriction} to  $S$
\index{Restriction!of blow-up sequence}%
\index{Blow-up!sequence, restriction}%
as the sequence 
$$
\begin{array}{ccccccl}
\Pi^S:S_r&\stackrel{\pi^S_{r-1}}{\longrightarrow} &S_{r-1}&
\stackrel{\pi^S_{r-2}}{\longrightarrow}\ \cdots\  
\stackrel{\pi^S_{1}}{\longrightarrow}& S_1&
\stackrel{\pi^S_{0}}{\longrightarrow} &S_0=S.\\
&& \cup && \cup && \, \cup \\
&& Z_{r-1}\cap S_{r-1} &\cdots & Z_1\cap S_1 && Z_0\cap S_0
\end{array}
$$
It is denoted by $j^*{\mathbf B}$ or  ${\mathbf B}|_S$.

Note that  $S_{i+1}:=B_{Z_i\cap S_i}S_i$ is naturally
identified with the birational transform $(\pi_i)^{-1}_*S_i\subset X_{i+1}$
(cf.\ \cite[II.7.15]{hartsh}),
thus there are natural embeddings $S_i\into X_i$ for every $i$.

The restriction of a 
smooth  blow-up sequence
need not be a smooth  blow-up sequence.
\medskip

\ref{bu.seq.functs}.3. 
Conversely, let  ${\mathbf B}(S):=$
$$
\Pi:=S_r\stackrel{\pi_{r-1}}{\longrightarrow} S_{r-1}
\stackrel{\pi_{r-2}}{\longrightarrow}\cdots
\stackrel{\pi_{1}}{\longrightarrow} S_1
\stackrel{\pi_{0}}{\longrightarrow} S_0=S
$$
be a  blow-up sequence with centers $Z^S_i\subset S_i$.
Define its {\it push-forward}
\index{Push-forward!of a blow-up sequence}%
\index{Blow-up!sequence, push-forward}%
as the sequence
$j_*{\mathbf B}:=$
$$
\Pi^X:X_r\stackrel{\pi^X_{r-1}}{\longrightarrow} X_{r-1}
\stackrel{\pi^X_{r-2}}{\longrightarrow}\cdots
\stackrel{\pi^X_{1}}{\longrightarrow} X_1
\stackrel{\pi^X_{0}}{\longrightarrow} X_0=X,
$$
whose centers $Z^X_i\subset X_i$ are defined inductively 
as $Z^X_i:=(j_i)_*Z^S_i$, where the $j_i:S_i\into X_i$
are the natural 
inclusions. Thus, for all practical purposes, $Z^X_i=Z^S_i$.

If  ${\mathbf B}$ is a smooth  blow-up sequence,
then so is $j_*{\mathbf B}$.
\end{defn}

\begin{defn}[Blow-up sequence functors]\label{busf.defnition}
A {\it blow-up sequence functor}%
\index{Blow-up!sequence functor}%
\index{Functoriality!of blow-up sequence}
is a functor $\bb$\index{BB@$\bb(X)$, blow-up sequence functor}
 whose
\begin{enumerate} 
\item inputs are triples
$(X,I,E)$, where $X$ is a scheme, $I\subset \o_X$ an ideal sheaf
that is nonzero on every irreducible component and
 $E$ a divisor on $X$ with ordered index set, and
\item  outputs are  blow-up sequences
$$
\begin{array}{ccccccl}
\Pi:X_r&\stackrel{\pi_{r-1}}{\longrightarrow} &X_{r-1}&
\stackrel{\pi_{r-2}}{\longrightarrow}\ \cdots\  
\stackrel{\pi_{1}}{\longrightarrow}& X_1&
\stackrel{\pi_{0}}{\longrightarrow} &X_0=X\\
&& \cup && \cup && \, \cup \\
&& Z_{r-1} &\cdots & Z_1 && Z_0
\end{array}
$$
with specified centers. Here the length of the sequence $r$,
the schemes $X_i$ and the centers $Z_i$ all depend on
$(X,I,E)$. (Later we will add ideal sheaves $I_i$ and divisors $E_i$
to the notation.)
\end{enumerate}
If each $Z_i$ is smooth, then a nontrivial 
blow-up $\pi_i:X_{i+1}\to X_i$ uniquely determines $Z_i$,
so we can drop $Z_i$ from the notation. However,
 in general  many different centers give the same
birational map.

The  (partial) {\it resolution  functor}\index{Resolution!functor}
\index{Functoriality!of resolution}
 $\res$ {\it associated}\index{Associated resolution functor}
\index{Resolution!associated to a blow-up sequence} to a
blow-up sequence functor $\bb$ is the functor that sends
$(X,I,E)$ to the end result of the
blow-up sequence
$$
\res: (X,I,E)\mapsto (\Pi:X_r\to X).
$$
Sometimes we write simply $\res_{(I,E)}(X)=X_r$.
\end{defn}

\begin{say}[Empty blow-up convention]\label{ebu.conv}
\index{Empty!blow-up convention} We basically try to avoid
empty blow-ups, but we are forced to deal with them because a pull-back 
or a restriction of a
nonempty blow-up may be an empty blow-up.

Instead of saying repeatedly that we perform a certain blow-up
unless its center  is empty, we adopt the convention
 that  the final outputs of the named blow-up sequence functors
$\bdd, \bmord, \bord, \bpri$ do not contain empty blow-ups.

The process of their construction may contain
blow-ups that are empty in certain cases.
(For instance, we may be told to blow up
$E^1\cap E^2$ and the intersection may be empty.)
These steps are then ignored without explicit mention
whenever they happen to lead to  empty blow-ups.
\end{say}

\begin{rem} \label{bu.dndet.ord.rem} The end result of a
sequence of blow-ups $\Pi:X_r\to X$ often determines the whole
sequence, but this is not always the case.

First, there are some genuine counterexamples.
Let $p\in C$ be a smooth pointed curve in
a smooth 3-fold $X_0$.
We can first blow up $p$ and then the birational transform of $C$
to get
$$
\Pi:X_2\stackrel{\pi_{1}}{\longrightarrow} X_{1}=B_pX_0
\stackrel{\pi_{0}}{\longrightarrow} X_0,
$$
with exceptional divisors $E_0,E_1\subset X_2$,
or we can blow up first $C$ and then the preimage  $D=\sigma_0^{-1}(p)$
to get
$$
\Sigma:X'_2\stackrel{\sigma_{1}}{\longrightarrow} X'_{1}=B_CX_0
\stackrel{\sigma_{0}}{\longrightarrow} X_0
$$
with exceptional divisors $E'_0,E'_1\subset X'_2$.

It is easy to see that $X_2\cong X'_2$, and under this isomorphism
$E_1$ corresponds to $E'_0$ and $E_0$ corresponds to $E'_1$.

Second, there are some ``silly'' counterexamples.
If $Z_1, Z_2\subset X$ are two disjoint smooth
subvarieties, then we get the same result whether
we blow up first $Z_1$ and then $Z_2$, or 
first $Z_2$ and then $Z_1$, or in one step we blow up
$Z_1\cup Z_2$. 

While it seems downright stupid to distinguish
between these three processes, it is precisely this
ambiguity that caused the difficulties in
(\ref{nonproj.res.loc.glob.say}.1).
\end{rem}

It is also convenient to have a 
unified way to look at the functoriality
properties of various resolutions.

\begin{say}[Functoriality package]\label{funct.pack.proc.say}
\index{Functoriality!package}
There are three  functoriality properties of blow-up sequence functors
 $\bb$ that we are
interested in. Note that in all three cases the claimed isomorphism
is unique, and hence the existence is a local question.

Functoriality for \'etale morphisms is an essential ingredient
of the proof. As noted in (\ref{functres.res.aside}.2), this is equivalent
to functoriality for smooth morphisms (\ref{funct.pack.proc.say}.1).
Independence of the base field (\ref{funct.pack.proc.say}.2) is  very useful
in applications, but it is not needed for the proofs.

Functoriality for closed embeddings (\ref{funct.pack.proc.say}.3)
is used
for resolution of singularities, but it is not
needed for the principalization theorems.
This property is quite delicate, and we are not able
to prove it in full generality, see (\ref{pr.to.res.say}).
\medskip

\ref{funct.pack.proc.say}.1 (Smooth morphisms).
\index{Functoriality!for smooth morphisms}
We would like our resolutions to commute with
smooth morphisms, and it is best to build this
into the blow-up sequence functors.

We say that a blow-up sequence functor
$\bb$ commutes with  $h$ 
if 
$$
\bb\bigl(Y,h^*I,h^{-1}(E)\bigr) =h^*\bb(X,I,E).
$$
This sounds quite reasonable until one notices that
even when $Y\to X$ is an open immersion it can happen that
$Z_0\times_XY$ is empty. It is, however, reasonable to expect that
a good blow-up sequence functor commutes with
smooth {\em surjections}.

Therefore, we say that $\bb$ commutes with smooth morphisms if
\index{Blow-up!sequence functor, commutes with smooth morphisms}%
\index{Commute!with smooth morphisms}%
\begin{enumerate}%\setcounter{enumi}{2}
\item[$\bullet$]
 $\bb$ commutes with every smooth surjection $h$, and
\item[$\bullet$]
 for every smooth  morphism $h$,  $\bb(Y,h^*I,h^{-1}(E))$ is 
obtained from the pull-back $h^*\bb(X,I,E)$ by deleting every 
blow-up $h^*\pi_i$ whose center is empty
and reindexing the resulting blow-up sequence.
\end{enumerate}

\ref{funct.pack.proc.say}.2 (Change of fields).
\index{Functoriality!for field change}%
We also would like the resolution
to be independent of the field we work with.

Let $\sigma:K\into L$ be a field extension.
Given a $K$-scheme of finite type $X_K\to \spec K$, we can view
$\spec L$ as a scheme over $\spec K$
(possibly not of finite type) and take the fiber product
$$
X_{L,\sigma}:=X_K\times_{\spec K}\spec L,
$$
which is an $L$-scheme of finite type.
If $I$ is an ideal sheaf and $E$ a divisor on $X$,
then similarly we get $I_{L,\sigma}$ and $E_{L,\sigma}$.

We say that
$\bb$ {\it commutes}%
\index{Blow-up!sequence functor, commutes with change of fields}%
\index{Commute!with change of fields} with  $\sigma$ if
 $\bb(X_{L,\sigma},I_{L,\sigma},E_{L,\sigma})$ is the
blow-up sequence
$$
\begin{array}{ccccccl}
\Pi_{L,\sigma}:(X_r)_{L,\sigma}&\stackrel{(\pi_{r-1})_{L,\sigma}}
{\longrightarrow} &(X_{r-1})_{L,\sigma}&
\cdots & (X_1)_{L,\sigma}&
\stackrel{(\pi_{0})_{L,\sigma}}{\longrightarrow} 
&(X_0)_{L,\sigma}.\\
&& \cup && \cup && \ \ \ \ \ \cup \\
&& (Z_{r-1})_{L,\sigma}&\cdots & (Z_1)_{L,\sigma} && (Z_0)_{L,\sigma}
\end{array}
$$
This property  will hold  automatically for all
blow-up sequence functors that we construct.

\ref{funct.pack.proc.say}.3 (Closed embeddings).
\index{Functoriality!for closed embeddings}%
In the proof of (\ref{res.sings.thm.1}) we constructed a resolution
of a variety $Z$ by choosing an embedding
of $Z$ into a  smooth variety $Y$.
In order to get a well-defined resolution, we need to know
that our constructions do not depend on the embedding chosen.
The key step is to ensure independence from further
embeddings $Z\into Y\into X$.

We say that  $\bb$  {\it commutes with closed embeddings}
\index{Blow-up!sequence functor, commutes with closed embeddings}%
\index{Commute!with  closed embeddings}
 if
$$
\bb(X,I_X,E)=j_*\bb(Y,I_Y,E|_Y),
$$
 whenever 
\begin{enumerate}
\item[$\bullet$] $j:Y\into X$ is a closed  embedding 
of smooth schemes, 
\item[$\bullet$] $0\neq I_Y\subset \o_{Y}$  and $0\neq I_X\subset \o_{X}$
are  ideal sheaves such that $\o_X/I_X= j_*(\o_{Y}/I_Y)$, and
\item[$\bullet$] $E$ is a simple normal crossing divisor  on $X$ such that
$E|_Y$ is also a simple normal crossing divisor on $Y$.
\end{enumerate}
\medskip

\ref{funct.pack.proc.say}.4 (Closed embeddings, weak form).
Let the notation and assumptions be as in (\ref{funct.pack.proc.say}.3).
We say that $\bb$  {\it weakly commutes with closed embeddings}%
\index{Blow-up!sequence functor, weakly commutes with closed embeddings}%
\index{Commute!with  closed embeddings, weakly}
if
$$
j^*\bb(X,I_X,E)=\bb(Y,I_Y,E|_Y).
$$
The difference appears only in the proof of 
(\ref{funct.princ.thm.IV}) given in (\ref{pf.ored->unord}).
At the beginning of the proof we blow up various intersections
of the irreducible components of $E$. Since these
intersections are not contained in $Y$, this commutes
with restriction to $Y$ but it does not commute with push forward.

\end{say}

The strongest form of monomialization is the following.

\begin{thm}[Principalization, IV]\label{funct.princ.thm.IV}
\index{Principalization!strong form}%
\index{Ideal sheaf!principalization}%
There is a  blow-up sequence functor 
$\bpri$\index{BP@$\bpri(X)$, principalization functor}
defined on all triples $(X,I,E)$, where
 $X$ is a smooth scheme of finite type over  a field of characteristic zero,
 $I\subset \o_X$ is an ideal sheaf that is not zero on any
irreducible component of $X$
and $E$ is a simple normal crossing divisor on $X$. $\bpri$
satisfies the following conditions.
\begin{enumerate}
\item In the blow-up sequence $\bpri(X,I,E)=$
$$
\begin{array}{ccccccl}
\Pi:X_r&\stackrel{\pi_{r-1}}{\longrightarrow} &X_{r-1}&
\stackrel{\pi_{r-2}}{\longrightarrow}\ \cdots\  
\stackrel{\pi_{1}}{\longrightarrow}& X_1&
\stackrel{\pi_{0}}{\longrightarrow} &X_0=X,\\
&& \cup && \cup && \, \cup \\
&& Z_{r-1} &\cdots & Z_1 && Z_0
\end{array}
$$
all centers of blow-ups are smooth and have simple normal crossing
with   $E$ (\ref{bir.tr.divs.def}).
\item The pull-back 
$\Pi^*I\subset\o_{X_r}$ is the ideal sheaf of a simple normal
crossing divisor.
\item $\Pi:X_r\to X$ is an isomorphism over $X\setminus \cosupp I$.
\item $\bpri$  commutes with smooth morphisms
(\ref{funct.pack.proc.say}.1) 
and with change of fields (\ref{funct.pack.proc.say}.2).
\item $\bpri$  commutes with closed embeddings (\ref{funct.pack.proc.say}.3)
whenever $E=\emptyset$.
\end{enumerate}
\end{thm}

Putting together the proof of
(\ref{res.sings.main.thm}) with (\ref{aff.to.glob.prop}),
we obtain strong and functorial resolution.

\begin{thm}[Resolution of singularities, III]\label{res.sings.proc.main.thm}
\index{Resolution!strong form}%
There is a  blow-up sequence functor 
$\bres(X)=$\index{BR@$\bres(X)$, resolution functor}
$$
\begin{array}{ccccccl}
\Pi:X_r&\stackrel{\pi_{r-1}}{\longrightarrow} &X_{r-1}&
\stackrel{\pi_{r-2}}{\longrightarrow}\ \cdots\  
\stackrel{\pi_{1}}{\longrightarrow}& X_1&
\stackrel{\pi_{0}}{\longrightarrow} &X_0=X,\\
&& \cup && \cup && \, \cup \\
&& Z_{r-1} &\cdots & Z_1 && Z_0
\end{array}
$$
defined on all 
 schemes $X$  of finite type over  a field of characteristic zero,
 satisfying the following conditions.
\begin{enumerate}
\item $X_r$ is smooth.
\item $\Pi:X_r\to X$ is an isomorphism over the smooth locus
$X^{ns}$.
\item $\Pi^{-1}(\sing X)$ is a divisor with simple normal crossings.
\item $\bres$  commutes with smooth morphisms
(\ref{funct.pack.proc.say}.1) 
and with change of fields (\ref{funct.pack.proc.say}.2).
\end{enumerate}
\end{thm}

Proof. First we construct $\bres(X)$ for affine schemes.
Pick any embedding $X\into A$ into a smooth affine scheme
such that $\dim A\geq \dim X+2$.
As in the proof of  (\ref{res.sings.thm.1}),
the blow-up sequence for $\bpri(A,I_X,\emptyset)$
obtained in   (\ref{funct.princ.thm.IV})
gives a blow-up sequence $\bres(X)$.

Before we can even consider the
functoriality conditions, we need to prove that
$\bres(X)$ is independent of the
choice of the embedding $X\into A$.

 Thus assume that
$\Pi_1:R_1(X)\to \cdots \to X$
and $\Pi_2:R_2(X)\to \cdots \to X$ are two blow-up sequences constructed
this way.
Using that $\bpri$ weakly commutes with closed embeddings
(\ref{funct.pack.proc.say}.4),
it is enough to prove   uniqueness 
for  resolutions  constructed from embeddings 
into affine spaces $X\into \a^n$.
Moreover, we are allowed to increase $n$ anytime
by taking a further embedding 
$\a^n\into\a^{n+m}$. 

As (\ref{aff.emb.stab.lem}) shows, 
any two embeddings $i_1,i_2:X\into \a^n$ become equivalent
by an automorphism of  $\a^{2n}$, which gives the required uniqueness.

Thus (\ref{funct.pack.proc.say}.2) for $\bpri(A, I_X,\emptyset)$
implies (\ref{funct.pack.proc.say}.2) for $\bres(X)$
since  an embedding
$i:X\into A$ over $K$ and $\sigma:K\into L$ gives 
another embedding
$i_{\sigma,L}:X_{\sigma,L}\into A_{\sigma,L}$.

We can also reduce the condition
(\ref{funct.pack.proc.say}.1) for $\bres(X)$
to the same condition for  $\bpri(A, I_X,\emptyset)$.

To see this, 
let $h:Y\to X$ be a smooth morphism,
and choose any embedding $X\into A_X$ into a smooth affine variety. 
By (\ref{smooth.emb.stab.lem}), for
every $y\in Y$  there is an open neighborhood
$h(y)\in A_X^0\subset A_X$ 
and a smooth surjection $h_A:A_Y^0\onto A_X^0$ such that
$h_A^{-1}(X\cap A_Y^0)$ is isomorphic to 
an open neighborhood $y\in Y^0\subset Y$. 
Set  $X^0:=X\cap A_X^0$. Thus,
by  (\ref{funct.pack.proc.say}.1), 
$$
h_A^* \bpri(A_X^0,I_{X^0},\emptyset)=\bpri(A_Y^0,I_{Y^0},\emptyset),
$$
which shows that $h^*\bres(X^0)=\bres(Y^0)$.
As we noted earlier, (\ref{funct.pack.proc.say}.1)
is a local property, and thus $h^*\bres(X)=\bres(Y)$ as required.

We have now defined $\bres$ on (possibly reducible) affine schemes, and
it remains to prove that  one can glue together a global
resolution out of these local pieces.
This turns out to be a formal property of
blow-up sequence functors, which we treat next.
\qed

\begin{prop}\label{aff.to.glob.prop}
\index{Blow-up!sequence functor}
\index{Functoriality!of blow-up sequence}
 Let $\bb$ be a blow-up sequence functor
defined on affine schemes over a field $k$
that commutes with smooth surjections.

Then $\bb$ has a unique extension to
a blow-up sequence functor $\overline{\bb}$,
which is defined on all schemes of finite type over  $k$ and
which commutes with smooth surjections.
\end{prop}

Proof. For any $X$ choose an open affine cover
$X=\cup U_i$, and let $X':=\coprod_i U_i$ be the disjoint union.
Then $X'$ is affine, and there is a smooth surjection
$g:X'\to X$. We show that $\bb (X')$ descends to
a blow-up sequence of $X$.

Set $X'':=\coprod_{i\leq j}U_i\cap U_j$.
(We can also think of it as the fiber product $X'\times_XX'$.)
There are surjective open immersions 
$\tau_1,\tau_2:X''\to X'$, where
$\tau_1|_{U_i\cap U_j}:U_i\cap U_j\to U_i$ is the first inclusion
and $\tau_2|_{U_i\cap U_j}:U_i\cap U_j\to U_j$ is the second.

The blow-up sequence $\bb(X')$  starts with blowing up
$Z'_0\subset X'$,
and the blow-up sequence $\bb(X'')$ starts with blowing up
$Z''_0\subset X''$.
Since $\bb$ commutes with the $\tau_i$, we conclude
that
$$
\tau_1^*(Z'_0)=Z''_0=\tau_2^*(Z'_0).
\eqno{(\ref{aff.to.glob.prop}.1)}
$$
Since
$Z'_0\subset X'$ is a disjoint union of its pieces
$Z'_{0i}:=Z'_0\cap U_i$, 
(\ref{aff.to.glob.prop}.1) is equivalent to saying that for every $i,j$
$$
Z'_{0i}|_{U_i\cap U_j}=Z'_{0j}|_{U_i\cap U_j}.
\eqno{(\ref{aff.to.glob.prop}.2)}
$$
Thus the subschemes $Z'_{0i}\subset U_i$ glue together
to a subscheme $Z_0\subset X$.

This way we obtain $X_1:=B_{Z_0}X$
such that $X'_1=X'\times_XX_1$. We can repeat the above argument
to obtain the center $Z_1\subset X_1$ and eventually get the
whole blow-up sequence for $X$.\qed

\begin{warning} A key element of the above argument
is that we need to know $\bb$ for the
disconnected affine scheme $\coprod_i U_i$.

Any resolution functor defined on connected schemes
automatically extends to  disconnected schemes, but
for blow-up sequence functors this is not at all the case.
Although the blow-ups on different connected components
do not affect each other, in a resolution process we need to know
in which order we perform them, see (\ref{nonproj.res.loc.glob.say}.1).

Besides proving resolution for
nonprojective schemes and for algebraic spaces, 
the method of (\ref{aff.to.glob.prop})
is used in the  proof of
the principalization theorems. The inductive proof
naturally produces resolution processes only locally,
and this method shows that they automatically globalize.
\end{warning}

The following lemma shows that an affine scheme has a
unique embedding into affine spaces, if we stabilize
the dimension.

\begin{lem} \label{aff.emb.stab.lem}
Let $X$ be an affine scheme and
$i_1:X\into \a^n$ and $i_2:X\into \a^m$ two closed embeddings.
Then the two embeddings into the  coordinate subspaces
$$
i'_1:X\into \a^n\into \a^{n+m}\qtq{and}
i'_2:X\into \a^m\into \a^{n+m}
$$
are equivalent under a (nonlinear) automorphism of $\a^{n+m}$.
\end{lem}

Proof. We can extend $i_1$ to a morphism $j_1:\a^m\to \a^n$ and
 $i_2$ to a morphism $j_2:\a^n\to \a^m$.

Let ${\mathbf x}$ be coordinates on $\a^n$ and
${\mathbf y}$  coordinates on $\a^m$.
Then
$$
({\mathbf x},{\mathbf y})\mapsto 
({\mathbf x}, {\mathbf y}+j_2({\mathbf x}))
$$
is an automorphism of $\a^{n+m}$, which sends
the image of $i'_1$ to 
$$
\im\bigl[ i_1\times i_2: X\to \a^n\times \a^m\bigr].
$$
Similarly, 
$$
({\mathbf x},{\mathbf y})\mapsto 
({\mathbf x}+j_1({\mathbf y}), {\mathbf y})
$$
is an automorphism of $\a^{n+m}$, which sends
the image of $i'_2$ to 
$$
\im\bigl[ i_1\times i_2: X\to \a^n\times \a^m\bigr].\qed
$$

\begin{aside} It is worthwhile to mention a
local variant of (\ref{aff.emb.stab.lem}).
Let $X$ be a scheme and $x\in X$  a point 
whose Zariski tangent space has dimension $d$.
Then, for  $m\geq 2d$, $x\in X$ has a unique embedding into a  smooth
scheme of dimension $m$, up to \'etale coordinate changes.

See \cite{jelon, kalim} for affine versions.
\end{aside}

\begin{lem} \label{smooth.emb.stab.lem}
Let $h:Y\to X$ be a smooth morphism, $y\in Y$ a point 
and  $i:X\into A$ a closed embedding. 
Then  there are open neighborhoods
$y\in Y^0\subset Y$,  $f(y)\in A_X^0\subset A_X$, $X^0=X\cap A^0_X$;
 a smooth morphism $h_A:A_Y^0\to A_X^0$; and a closed embedding
$j: Y^0\into A_Y^0$ such that 
 the following diagram is a fiber product square:
 $$
\begin{array}{ccc}
Y^0 & \stackrel{j}{\hookrightarrow} & A_Y^0\\
h\downarrow \hphantom{h} &\square & \hphantom{h_A}\downarrow h_A\\
 X^0 & \stackrel{i}{\hookrightarrow} & A_X^0.
\end{array}
$$
\end{lem}

Proof. We prove this over infinite fields, which is the only case that we use.

 The problem is local, and thus we may assume that
$X,Y, A_X$ are affine and  $Y\subset X\times \a^N$.
If $h$ has relative dimension $d$, choose a general projection
$\sigma:\a_x^N\to \a_x^{d+1}$ such that
$\sigma:h^{-1}(x)\to  \a_x^{d+1}$ is finite and 
an embedding in a neighborhood of $y$.
(Here we need that the residue field of $x$ is infinite.)
Thus, by shrinking $Y$, we may assume that
$Y$ is an open subset of a hypersurface
$H\subset X\times \a^{d+1}$ and the first projection is smooth at
$y\in H$. $H$ is defined by an equation
$\sum_I \phi_I z^I$, where the $\phi_I$ are regular functions on
$X$ and  $z$ denotes the coordinates on $\a^{d+1}$.
Since $X\into A_X$ is a closed embedding, the $\phi_I$
extend  to regular functions  $\Phi_I$ on $A_X$.
Set
$$
A_Y:=(\sum_I \Phi_Iz^I=0)\subset A_X\times \a^{d+1}.
$$
Thus $Y\subset A_Y$ and the projection $A_Y\to A_X$ is smooth at $y$.
Let $y\in A_Y^0\subset A_Y$  and $A_X^0\subset A_X$
be  open sets 
such that the projection $h_A:A_Y^0\to A_X^0$ is smooth and surjective.
Set $Y^0:=Y\cap A_Y^0$.
\qed 
\medskip

The following comments on resolution for
algebraic and analytic spaces are not used
elsewhere in these notes.

\begin{say}[Algebraic spaces]\label{algsp.res.say}
\index{Algebraic space}

All we need to know about algebraic spaces is that
\'etale locally they are like schemes. That is,
there is a (usually nonconnected) scheme of
finite type $U$ and an \'etale surjection
$\sigma:U\to X$. We can even assume that $U$ is affine.

The fiber product
$V:=U\times_XU$ is again a 
scheme of
finite type with two surjective, \'etale projection morphisms
$\rho_i:V\to U$, and for all purposes one can identify
the algebraic space with the diagram
of schemes
$$
X=\bigl[\rho_1,\rho_2:V \rightrightarrows U\bigr].
\eqno{(\ref{algsp.res.say}.1)}
$$
The argument of (\ref{aff.to.glob.prop}) applies
to show that any  blow-up sequence functor $\bb$ that is
defined on affine schemes over a field $k$ 
and commutes with \'etale surjections,
 has a unique extension to
a blow-up sequence functor $\overline{\bb}$,
which is defined on all algebraic spaces over  $k$. 
(See (\ref{glob.bus.prop}) for details.)
Thus we obtain the following.
\end{say}

\begin{cor}\label{funct.princ.algsp.thm}
\index{Principalization!for algebraic spaces}
\index{Ideal sheaf!principalization for algebraic spaces}
\index{Resolution!for algebraic spaces}
The theorems (\ref{funct.princ.thm.IV}) and (\ref{res.sings.proc.main.thm})
also hold for algebraic spaces of finite type over
a field of characteristic zero.\qed
\end{cor}

\begin{say}[Analytic spaces]\label{analsp.res.say}
\index{Analytic space}

It was always understood that a good resolution method
should also work for complex, real  or $p$-adic analytic spaces.
 (See \cite{gr-re.an.st} for an introduction to
analytic spaces.)

The traditional methods almost all worked well locally,
but globalization sometimes presented  technical 
difficulties.
We leave it to the reader to follow the proofs
in this chapter and see that they all extend to
analytic spaces over locally compact fields,
at least locally. Once, however, we have a locally defined
blow-up sequence functor that  commutes with smooth surjections,
the argument of (\ref{aff.to.glob.prop}) 
 shows that
we get a globally defined
  blow-up sequence functor 
for small neighborhoods of compact sets
on all analytic spaces. Once we have a resolution
functor on neighborhoods of compact sets that commutes
with open embeddings, we get resolution for any
analytic space that is an increasing union of
its compact subsets.
Thus we obtain the following.
\end{say}

\begin{thm}\label{funct.princ.analsp.thm}
\index{Resolution!for analytic spaces}
Let $K$ be a locally compact field of characteristic zero.
There is a resolution functor
$\res:X\to (\Pi_X:R(X)\to X)$ defined on all separable
$K$-analytic spaces with the following properties.
\begin{enumerate}
\item $R(X)$ is smooth.
\item $\Pi:R(X)\to X$ is an isomorphism over 
the smooth locus $X^{ns}$.
\item $\Pi^{-1}(\sing X)$ is a
divisor with simple normal crossing.
\item  $\Pi_X$ is  projective over
any compact subset of $X$.
\item $\res$ commutes with smooth $K$-morphisms.\qed
\end{enumerate}
\end{thm}

\begin{aside}\label{no.smooth.emb.exmp}
We give an example of a normal, proper surface $S$ 
over $\c$ that cannot be embedded
into a smooth scheme. 

Start with $\p^1\times C$, where $C$ is any smooth curve
of genus $\geq 1$. Take two points $c_1,c_2\in C$.
Blow up $(0,c_1)$ and $(\infty,c_2)$
to get $f:T\to \p^1\times C$. We claim the following.
\begin{enumerate}
\item
 The birational transforms
$C_1\subset T$ of  $\{0\}\times C$ and  
$C_2\subset T$ of $\{\infty\}\times C$
can be contracted, and we get a normal, proper surface $g:T\to S$.
\item  If $\o_C(c_1)$ and $\o_C(c_2)$ are
independent in $\pic(C)$, then $S$ can not be embedded
into a smooth scheme. 
\end{enumerate}

To get the first part, it is easy to check that  a  multiple  of
the birational transform
 of  $\{1\}\times C+\p^1\times\{c_i\}$
on $T$ is base point free and contracts $C_i$ only, giving $g_i:T\to S_i$.
Now $S_1\setminus C_2$ and $S_2\setminus C_1$ can be glued together
to get $g:T\to S$.

If  $D$  is a Cartier divisor  on $S$, then 
$\o_T(g^*D)$ is trivial on both
$C_1$ and $C_2$. Therefore, $f_*(g^*D)$ is a  Cartier divisor  on 
$\p^1\times C$
such that its restriction to $\{0\}\times C$
is linearly equivalent to a multiple of $c_1$
and its restriction to $\{\infty\}\times C$
is linearly equivalent to a multiple of $c_2$.

Since $\pic(\p^1\times C)=\pic(C)\times \z$  
and $\o_C(c_1)$ and $\o_C(c_2)$ are
independent in $\pic(C)$, 
every Cartier divisor on $S$ is linearly equivalent to a multiple
of $\{1\}\times C$.
Thus the points of $\{1\}\times C\subset S$ cannot be separated from 
each other 
by Cartier divisors on $S$.
 
Assume now that $S\into Y$ is an  embedding
into a smooth scheme. Pick a point $p\in \{1\}\times C\subset Y$,
and let $p\in U\subset Y$ be an affine neighborhood.
Any two points of $U$ can be separated from 
each other 
by Cartier divisors on $U$.  Since $Y$ is smooth,
the  closure of a Cartier divisor on $U$ is 
automatically Cartier on $Y$.
Thus any two points of $U\cap S$ can be separated from 
each other 
by Cartier divisors on $S$, a contradiction.\qed

An example of a toric  variety with no Cartier divisors
is given in \cite[p.65]{fult-toric}. This again has no smooth embeddings.
\end{aside}

\section{Plan of the proof}

This section contains a still somewhat informal review
of the main steps of the proof. For simplicity, the role
of the divisor $E$ is ignored for now. All the definitions and
theorems will be made precise later.

We need some way to measure
how complicated an ideal sheaf is at a point.
For the present proof a very crude measure---the
order of vanishing or, simply, order---is enough. 

\begin{defn} \label{order.defn}
Let $X$ be a smooth variety and $0\neq I\subset \o_X$ an ideal sheaf.
For a point $x\in X$ with  ideal sheaf $m_x$, we define the
{\it order of vanishing} or {\it order of} $I$ at $x$ to be
\index{Ord@$\ord_xI$, order of vanishing}
\index{Order!of vanishing}
$$
\ord_xI:=\max\{r: m_x^r\o_{x,X}\supset I\o_{x,X}\}.
$$
It is easy to see that $x\mapsto \ord_xI$ is a
constructible and upper-semi-continuous function on $X$.

For an irreducible subvariety $Z\subset X$, we define
the {\it order of} $I$ along $Z\subset X$ as
$$
\ord_ZI:=\ord_{\eta}I,\qtq{where $\eta\in Z$ is the generic point.}
$$
 Frequently we also use the notation
$\ord_ZI=m$ 
(resp., $\ord_ZI\geq m$) when $Z$ is not irreducible. In this case
we always assume that the order of $I$ at every generic point
of $Z$ is $m$ (resp., $\geq m$).

The {\it maximal order} of $I$ along $Z\subset X$ is
\index{Maxo@$\mord_ZI$, maximal order of vanishing}
\index{Order!of vanishing, maximal}
$$
\mord_ZI:=\max\{\ord_zI: z\in Z\}.
$$
We frequently use $\mord I$ to denote $\mord_XI$.
\end{defn}

If $I=(f)$ is a principal ideal, then the order of $I$ at a point $x$
is the same as the multiplicity of the hypersurface $(f=0)$ at $x$.
This is a simple but quite strong invariant.

In general, however, the order is a very stupid invariant.
For resolution of singularities we always start with
an embedding $X\into \p^N$, where $N$ is larger than the
embedding dimension of $X$ at any point. Thus the ideal
sheaf $I_X$ of $X$ contains an order 1 element 
 at every point (the local equation of a smooth hypersurface containing $X$),
so the order of $I_X$ is 1 at every point of $X$. 
Hence the order of $I_X$ does not ``see'' the singularities
of $X$ at all.
(In the proof given in Section 3.12, trivial 
steps reduce the principalization of the ideal sheaf of
$X\subset  \p^N$ near a point $x\in X$ to the 
principalization of the ideal sheaf of
$X\subset P$, where $P\subset  \p^N$ is smooth and
has the smallest possible dimension locally near $x$. 
Thus we start actual work only when $\ord I\geq 2$.)

There  is one useful property of  $\ord_ZI$, which
is exactly what we need:
the number $\ord_ZI$ equals the 
multiplicity of  $\pi^*I$ along the exceptional divisor 
of  the blow-up
$\pi:B_ZX\to X$.

\begin{defn}[Birational transform of ideals]\label{bit.tr.intro.defss}

Let $X$ be a smooth variety and $I\subset \o_X$ an ideal sheaf.
For $\dim X\geq 2$ an ideal cannot be written as
the product of prime ideals, but the codimension 1 primes can
be separated from the rest. That is, there is a unique
largest effective divisor $\operatorname{Div}(I)$ such that
$I\subset \o_X(-\operatorname{Div}(I))$, and we can write
$$
I=\o_X\bigl(-\operatorname{Div}(I)\bigr)\cdot I_{{\rm cod}\geq 2},\qtq{where
 $\codim\supp (\o_X/I_{{\rm cod}\geq 2})\geq 2$.}
$$
We call $\o_X\bigl(-\operatorname{Div}(I)\bigr)$ the {\it divisorial part}
\index{Divisorial part}\index{Ideal sheaf!divisorial part}%
\index{Divi@$\operatorname{Div}(I)$, divisorial part}%
 of $I$ and
$I_{{\rm cod}\geq 2}=\o_X\bigl(\operatorname{Div}(I)\bigr)\cdot I$
 the {\it codimension $\geq 2$ part} of $I$.

Let $f:X'\to X$ be a birational morphism between smooth varieties.
Assume for simplicity that $I$ has no divisorial part, that is,
$I=I_{{\rm cod}\geq 2}$. 
We are interested in the codimension $\geq 2$ part of $f^*I$,
called the  {\it birational transform} of $I$
\index{Birational!transform of an ideal sheaf}
\index{Ideal sheaf!birational transform}
\index{Transform!birational, of ideal sheaf}
and
denoted by $f^{-1}_*I$.
(It is also frequently called the weak transform in the literature.)
Thus
 $$
f^{-1}_*I=\o_{X'}\bigl(\operatorname{Div}(f^*I)\bigr)\cdot f^*I.
$$
We have achieved principalization iff the
codimension $\geq 2$ part of $f^*I$ is not there, that is, when
$f^{-1}_*I=\o_{X'}$.

For reasons connected with 
 (\ref{probls.with.meth.say}.2),
 we also need another version, where
we ``pretend'' that $I$ has order $m$.

A {\it marked ideal sheaf}  on $X$ is a pair $(I,m)$
\index{Ideal sheaf!marked}\index{Im@$(I,m)$, marked ideal sheaf}%
where $I\subset \o_X$ is an ideal sheaf on  $X$ and $m$ is
a natural number.

Let $\pi:B_ZX\to X$ be the blow-up of a smooth  subvariety $Z$ 
 and $E\subset B_ZX$  the exceptional divisor.
Assume that $\ord_ZI\geq m$. Set
$$
\pi^{-1}_*(I,m):=\bigl(\o_{B_ZX}(mE)\cdot \pi^*I,m\bigr),
$$
and call it the 
 {\it birational transform} of $(I,m)$.
\index{Birational!transform of a marked ideal}

If $\ord_ZI= m$, then this coincides with $f^{-1}_*I$, but
for $\ord_ZI> m$ the cosupport of  $f^{-1}_*(I,m)$ 
also contains $E$. 
(We never use the case where $\ord_ZI< m$, since then 
$f^{-1}_*(I,m)$ is not an ideal sheaf.)
One can iterate this procedure to
define $f^{-1}_*(I,m)$  whenever 
$f:X'\to X$ is the composite of blow-ups of  smooth irreducible subvarieties
as above, but one has to be quite careful with this; see
 (\ref{bir.tr.marked.warning}).
\end{defn}

\begin{say}[Order reduction theorems]\label{ord.red.intro.thms}

The technical core of the proof
consists of two order reduction theorems
using smooth blow-ups that match the order that we work with.

Let $I$ be an ideal sheaf with $\mord I\leq m$.
A {\it smooth blow-up sequence}
\index{Blow-up!sequence of order $m$}
\index{Order!blow-up sequence of order $m$}
 of order $ m$  starting with $(X,I)$ is a
smooth  blow-up sequence (\ref{bu.seq.notation.defn})
$$
\Pi:(X_r,I_r)\stackrel{\pi_{r-1}}{\longrightarrow} 
(X_{r-1},I_{r-1})
\stackrel{\pi_{r-2}}{\longrightarrow}\cdots
\stackrel{\pi_{1}}{\longrightarrow} (X_1,I_1)
\stackrel{\pi_{0}}{\longrightarrow} (X_0,I_0)=(X,I),
$$
where  each $\pi_i:X_{i+1}\to X_i$ is a  smooth blow-up with center
$Z_i\subset X_i$, 
 the $I_i$ are defined recursively by the
formula
$I_{i+1}:=(\pi_i)^{-1}_*I_i$ and $\ord_{Z_i}I_i=m$ for every $i<r$.

A {\it  blow-up sequence}  of order $ \geq m$ 
\index{Blow-up!sequence of order $\geq m$}
\index{Order!blow-up sequence of order $\geq m$}
 starting with a marked ideal $(X,I,m)$ is 
defined analogously, except  we use the recursion
formula
$(I_{i+1},m):=(\pi_i)^{-1}_*(I_i,m)$ and 
we require $\ord_{Z_i}I_i\geq m$ for every $i<r$.

Using these notions, the inductive versions of the main results 
are the following.
\medskip

\ref{ord.red.intro.thms}.1 (Order reduction for ideals).
\index{Order reduction!for ideal sheaves} {\it 
Let $X$ be a smooth variety,  $0\neq I\subset \o_X$  an ideal sheaf
and $m=\mord I$. 
By  a suitable blow-up sequence of order $m$
we  eventually get
$f:X'\to X$ such that $\mord f^{-1}_*I<m$.
}
\medskip

\ref{ord.red.intro.thms}.2 (Order reduction for marked ideals).
\index{Order reduction!for marked ideal sheaves}  {\it 
Let $X$ be a smooth variety, $0\neq I\subset \o_X$  an ideal sheaf
and  $m\leq \mord I$ a natural number. 
By a suitable blow-up sequence of  order $\geq m$,
we eventually get
$f:X'\to X$ such that $\mord f^{-1}_*(I,m)<m$.
}
\medskip

We prove these theorems
together in  a spiraling induction with 
two main reduction steps.
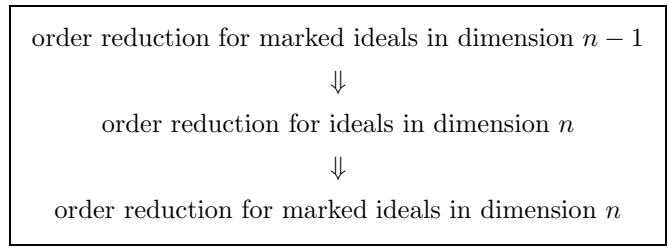
\index{Order reduction!spiraling induction} 
$$
\boxed{
\begin{array}{c}
\ \\*[-2ex]
\mbox{order reduction for marked ideals 
in dimension $n-1$}\\*[1ex]
\Downarrow\\*[1ex]
\mbox{order reduction for ideals
in dimension $n$}\\*[1ex]
 \Downarrow\\*[1ex]
\mbox{order reduction for marked ideals 
in dimension $n$}\\*[-2ex]
\ 
\end{array}
}
$$
\end{say}

The two steps are independent and use different methods.

The second implication is relatively easy and has been well understood for
a long time. We leave it  to Section 3.13.

Here we focus on the proof of the harder part, which is the
first implication.

\begin{say}[The heart of the proof]\label{heart.say}
\index{Order reduction!heart of the proof}
\index{Resolution!heart of the proof}  Methods to deal
with  Problems (\ref{probls.with.meth.say}.3--5)
 form the key steps of the proof.
My approach is to break apart the traditional inductive proof.
The problems can be solved  independently but only for certain ideals.
Then we need one more  step to show that order reduction for
 an arbitrary ideal
is equivalent to  order reduction for an ideal with all the required good
properties. 
\medskip

\ref{heart.say}.1 (Maximal contact). This deals with
 (\ref{probls.with.meth.say}.3) by showing that for suitable
hypersurfaces $H\subset X$ every step of an order reduction 
algorithm for $(X,I)$ with $m=\mord I$ is also a step of an
order reduction 
algorithm for $(H,I|_H,m)$.
This is explained in 
(\ref{max.cont.intr.say}) and completed in Section 3.8. 
\medskip

\ref{heart.say}.2 ($D$-balanced ideals). Problem 
(\ref{probls.with.meth.say}.4) 
has a solution for certain ideals only. For the so-called
$D$-balanced ideals, the converse of maximal contact
theory holds. That is,
for every 
hypersurface $S\subset X$, every  
order reduction 
step for $(S,I|_S,m)$ is also an
order reduction  step
 for $(X,I)$. This is outlined in (\ref{D-balance.say})
with all details in Section 3.9.

\medskip

\ref{heart.say}.3 (MC-invariant ideals). The solution of   
(\ref{probls.with.meth.say}.5) requires the
consideration of  MC-invariant ideals.
 For these, all hypersurfaces of maximal contact are
locally analytically isomorphic, with an
isomorphism preserving the ideal $I$.
See (\ref{mc-homog.say}),
with full proofs in Section 3.10.
\medskip

\ref{heart.say}.4 (Tuning of ideals).
It remains to show that order reduction for
an arbitrary ideal $I$ is
equivalent to order reduction 
for an ideal $W(I)$, which is both $D$-balanced
and MC-invariant. This  turns out to be
surprisingly easy; see (\ref{tuning.gen.say}) and Section 3.11.
\medskip

\ref{heart.say}.5 (Final assembly).
The main remaining problem is that a hypersurface of maximal contact
can be found only locally, not globally.
The local pieces are united in Section 3.12,
where we also  take
care of the divisor $E$, which we have ignored so far.
\end{say}

Let us now see these steps in more detail.

\begin{say}[Maximal contact]
\label{max.cont.intr.say}
\index{Maximal contact!outline}
Following the examples (\ref{deg1Econe.w.meth}) and  
(\ref{messy.exmp.w.meth}), given $X$ and $I$
with $m=\mord I$, we would like to
find a smooth hypersurface $H\subset X$ such that
order reduction for $I$  follows from
order reduction for $(I|_H, m)$.

As we noted in   (\ref{probls.with.meth.say}.3), first, we have to
ensure that the points where the birational transform of $I$ has
order $\geq m$ stay on the birational transform of $H$ all the time.
That is, we want to achieve the following.
\medskip

\ref{max.cont.intr.say}.1 (Going-down property of maximal contact).
\index{Maximal contact!going down}%
\index{Going!down}%
Restriction (\ref{bu.seq.functs}.2.) from $X$ to $H$ gives an injection
$$
\boxed{
\begin{array}{c}
\ \\*[-2ex]
\mbox{blow-up sequences of order $m$ for  $(X,I)$}\\*[1ex]
\bigcap\\*[1ex]
\mbox{blow-up sequences of order $\geq m$ for $(H,I|_H,m)$}\\*[-2ex]
\ 
\end{array}
}
$$
If this holds, then we say that $H$ is a hypersurface of {\it maximal contact}.
At least locally these are  easy to find  using derivative ideals.

Derivations of a smooth variety $X$ form a  sheaf $\der_X$, locally
generated by the usual partials  $\partial/\partial x_i$.
For an ideal sheaf $I$, let $D(I)$ denote the ideal sheaf
generated by all derivatives of local sections of $I$.
We can define higher-derivative ideals inductively by the rule
$D^{i+1}(I):=D\bigl(D^i(I)\bigr)$.

If $m=\mord I$, we are especially interested in 
the largest nontrivial derivative ideal. It is
also called the
{\it ideal of maximal contacts}
$$
MC(I):=D^{m-1}(I)=
\left(\ \frac{\partial^{m-1}f}{\partial x_1^{c_1}\cdots \partial x_n^{c_n}}:
f\in I,\ \sum c_i=m-1\right).
$$ 
\medskip

\ref{max.cont.intr.say}.2 (Local construction of maximal contact).
\index{Maximal contact!local existence}
 {\it For  a point $p\in X$ with $m=\ord_pI$, let
$h\in MC(I)$ be any local section with $\ord_ph=1$.
Then $H:=(h=0)$
 is a hypersurface of maximal contact in an open neighborhood
of $p$.
}

In general, hypersurfaces of maximal contact 
do not exist globally, and they are not unique locally. 
We deal with these problems later.

\end{say}

\begin{say}[$D$-balanced ideals]
\label{D-balance.say}
\index{D-bal@$D$-balanced ideal sheaf}
\index{Ideal sheaf!D-bal@$D$-balanced}
It is harder to deal with   (\ref{probls.with.meth.say}.4).
No matter how we choose the hypersurface of maximal contact $H$,
sometimes the restriction $(I|_H, m)$ is ``more singular'' than
$I$,
in the sense that order reduction for
$(I|_H, m)$ may involve blow-ups that are not needed
for any 
order reduction procedure of $I$; see (\ref{mc.badrestr.exmp}).

There are, however, some ideals for which this problem does not happen.
To define these, we again need to consider derivatives.

If $\ord_pf=m$, then typically $\ord_p(\partial f/\partial x_i)=m-1$,
so a nontrivial ideal is never $D$-closed. The best one can hope for is
that $I$ is $D$-closed, after we ``correct for
the lowering of the order.''

An ideal $I$ with $m=\mord I$ is called {\it $D$-balanced} if
$$
\bigl(D^i(I)\bigr)^m\subset I^{m-i}\quad \forall\ i<m.
$$
Such ideals behave very well with respect to   restriction
to  smooth subvarieties and smooth blow ups.
\medskip

\ref{D-balance.say}.1 (Going-up property of $D$-balanced ideals).
\index{Going!up}%
\index{D-bal@$D$-balanced ideal sheaf!and going up}%
Let $I$ be a $D$-balanced ideal with $m=\mord I$.
Then for any smooth hypersurface $S\subset X$
such that $S\not\subset \cosupp I$, 
push-forward (\ref{bu.seq.functs}.3.) from $S$ to $X$ gives an injection
$$
\boxed{
\begin{array}{c}
\ \\*[-2ex]
\mbox{blow-up sequences of order $m$ for  $(X,I)$}\\*[1ex]
\bigcup\\*[1ex]
\mbox{blow-up sequences of order $\geq m$ for $(S,I|_S,m)$}\\*[-2ex]
\ 
\end{array}
}
$$

{\it Example} \ref{D-balance.say}.2.  Start with the double point
ideal $I=(xy-z^n)$. Restricting to $S=(x=0)$ creates an $n$-fold line,
and blowing up this line is not an order 2 blow-up for  $I$.

We can see that
$$
I+D(I)^2=(xy, x^2, y^2,xz^{n-1}, yz^{n-1}, z^n)
$$
 is  $D$-balanced.
If we restrict $I+D(I)^2$ to $(x=0)$, we get the ideal
$(y^2, yz^{n-1}, z^n)$.
It is easy to check that the whole resolution of $S$
is correctly predicted by order reduction for  
$(y^2, yz^{n-1}, z^n)$.
\medskip

Putting (\ref{max.cont.intr.say}.1) and 
 (\ref{D-balance.say}.1) together, we
get the first  dimension reduction result.
\medskip

{\it Corollary} \ref{D-balance.say}.3. 
(Maximal contact for $D$-balanced ideals).
{\it Let $I$ be a $D$-balanced ideal with $m=\mord I$
and $H\subset X$ a smooth hypersurface of maximal
contact.  
Then we have an equivalence  
$$
\boxed{
\begin{array}{c}
\ \\*[-2ex]
\mbox{blow-up sequences of order $m$ for  $(X,I)$}\\*[1ex]
||\\*[1ex]
\mbox{blow-up sequences of order $\geq m$ for $(H,I|_H,m)$}\\*[-2ex]
\ 
\end{array}
}
$$
}
\end{say}

This equivalence suggests that the choice of $H$ should not be
important at all. However, in order to ensure functoriality we have to
choose a particular resolution. Thus we still need to show that our 
particular choices are independent of $H$.
A truly ``canonical'' resolution process would probably take care of such
problems automatically, but it seems that one has to make 
at least some artificial choices.

\begin{say}[MC-invariant ideals]\label{mc-homog.say}
\index{Maximal contact!ideal sheaf, invariant for}
\index{Ideal sheaf!maximal contact invariant}
\index{Ideal sheaf!MC-invariant}
\index{MC-invariant ideal sheaf}
\index{MCI@$MC(I)$, ideal of maximal contacts}
Dealing with  (\ref{probls.with.meth.say}.5)
is again possible only for certain ideals.

We say that  $I$ is {\it maximal contact invariant} or
{\it MC-invariant} if
$$
MC(I)\cdot D(I)\subset I.
\eqno{(\ref{mc-homog.say}.1)}
$$
Written in the equivalent form
$$
D^{m-1}(I)\cdot D(I)\subset I,
\eqno{(\ref{mc-homog.say}.2)}
$$
it is quite close in spirit to the
$D$-balanced condition.
The expected order of $D^{m-1}(I)\cdot D(I)$ is $m$,
so it is sensible to require  inclusion. There is no
need to correct for the change of order first.

For MC-invariant ideals the hypersurfaces of maximal contact
are still not unique, but different choices are equivalent
under local analytic isomorphisms (\ref{completions.def}).
\medskip

\ref{mc-homog.say}.3 (Formal uniqueness of maximal contact).
\label{anal.inv.say.4}%
\index{Maximal contact!formal uniqueness}%
 {\it
Let $I$ be an MC-in\-variant ideal sheaf on $X$ 
and    $H_1,H_2\subset X$  two
hypersurfaces of maximal contact through a point $x\in X$.
Then there is a local analytic automorphism
(\ref{completions.def}) 
$\phi: (x\in \hat X)\to (x\in \hat X)$ such that
\begin{enumerate}
\item[(i)] $\phi^{-1}(\hat H_1)=\hat H_2$, 
\item[(ii)]  $\phi^*\hat I=\hat I$, and
\item[(iii)] $\phi$ is the identity on $\cosupp \hat I$.
\end{enumerate}
}

\end{say}

\begin{say}[Tuning of ideals]\label{tuning.gen.say}
\index{Ideal sheaf!tuning}
\index{Tuning of ideal sheaves}
Order reduction using dimension induction is now in
quite good shape for ideals that are both $D$-balanced
and MC-invariant.

 The rest is taken care of by ``tuning'' the ideal $I$ first.
(I do not plan to give a precise meaning to the word ``tuning.''
The terminology follows
\cite{wlo}. The notion of tuning used in
\cite{enc-hau} is quite different.)
There are  many ways to tune an ideal; here is
one of the simplest ones.

To an ideal $I$ of order $m$, we would like to associate
the ideal generated by all products of derivatives of order at least $m$.
The problem with this is that if $f$ has order $m$,
then $\partial f/\partial x_i$ has order $m-1$, and so we 
are able to add $(\partial f/\partial x_i)^{2}$ 
(which has order $2m-2$), but we really
would like to add $(\partial f/\partial x_i)^{m/(m-1)}$
(which should have order $m$ in any reasonable definition).

We can avoid these fractional exponent problems by
working with all products of derivatives  whose order is
sufficiently divisible. For instance, the condition
$(\mbox{order})\geq m!$ works. 

Enriching an ideal with its derivatives was
used by Hironaka \cite{hir-ide}  and then developed by
 Villamayor \cite{vil-con}.
A larger ideal is introduced in \cite{wlo}.
The ideal $W(I)$ introduced below is even larger, and
this largest choice seems more
natural to me. 
That is, we set
$$
W(I):=
\left(
\prod_{j=0}^m  \bigl(D^{j}(I)\bigr)^{c_j} : \sum (m-j)c_j\geq m!
\right)\subset \o_X.
\eqno{(\ref{tuning.gen.say}.1)}
$$

The ideal $W(I)$ has all the properties that we need.

\medskip

 {\it Theorem} \ref{tuning.gen.say}.2. (Well-tuned  ideals). 
\index{Ideal sheaf!tuning}%
\index{Tuning of ideal sheaves}%
Let $X$ be a smooth variety, $0\neq I\subset \o_X$ an ideal sheaf and
$m=\mord I$. Then
\begin{enumerate}
\item[(i)] $\mord W(I)=m!$, 
\item[(ii)] $W(I)$ is   $D$-balanced,
\item[(iii)] $W(I)$ is    MC-invariant, and
\item[(iv)] there is an equivalence
$$
\boxed{
\begin{array}{c}
\ \\*[-2ex]
\mbox{blow-up sequences of order $m$ for  $(X,I)$}\\*[1ex]
||\\*[1ex]
\mbox{blow-up sequences of order $m!$ for $(X, W(I))$}\\*[-2ex]
\ 
\end{array}
}
$$\end{enumerate}
\medskip

\ref{tuning.gen.say}.3.   
It should be emphasized that there are many different ways to choose
an ideal with the properties of $W(I)$ as above, but all known choices have
rather high order.

I chose the order $m!$
for notational simplicity, but one could work with
any multiple of $\lcm(1,2,\dots,m)$ instead.
The smallest choice would be $\lcm(1,2,\dots,m)$, which is
roughly like $e^m$.  As discussed in (\ref{res.other.say}.2),
this is still too big for effective computations.
Even if we fix the order to be $m!$, many choices remain.

\ref{tuning.gen.say}.4. 
Similar constructions are also considered by Kawanoue
\cite{kaw-the} and by Villamayor \cite{vill-rees}.
\end{say}

\medskip

\begin{defn}[Completions]\label{completions.def}
\index{Completion!of a scheme}
This is the only piece of commutative algebra that we use.

For  a local ring $(R,m)$ its {\it completion} in the $m$-adic topology
is denoted by $\hat R$; cf.\ \cite[Chap.10]{at-ma}.
If $X$ is a $k$-variety and $x\in X$, then we denote by
$\hat X_x$  or by $\hat X$ the completion of $X$ at $x$,
which is $\spec_k \widehat{\o}_{x,X}$.

We say that $x\in X$ and $y\in Y$ are
{\it formally} 
\index{Formally isomorphic} isomorphic
if $\hat X_x$ is isomorphic to $\hat Y_y$.

We need Krull's intersection theorem (cf.\ \cite[10.17]{at-ma}),
which says that
for an ideal $I$ in a Noetherian local ring $(R,m)$ we have 
$$
I=\cap_{s=1}^{\infty} \bigl( I+m^s\bigr).
$$
In geometric language this implies that 
 if  $Z,W\subset X$ are two subschemes
such that $\hat Z_x=\hat W_x$, then there is an open neighborhood
$x\in U\subset X$ such that $Z\cap U=W\cap U$.

If $p\in X$ is  closed, 
then $\widehat{\o}_{p,X}\cong k(p)[[x_1,\dots,x_n]]$,
where $x_1,\dots,x_n$ are local coordinates.
If $k(p)=k$ or, more generally, when  there is
a {\it field of representatives}
\index{Field!of representatives} (that is,
 a subfield $k'\subset \widehat{\o}_{p,X}$
isomorphic to $k(p)$), this is proved
in \cite[II.2]{shaf}.
In characteristic zero one can find $k'$ as follows.
The finite field extension $k(p)/k$ is generated by a simple root
of a polynomial $f(y)\in k[y]\subset \widehat{\o}_{p,X}[y]$. 
Modulo the maximal ideal, $f(y)$ has a linear factor by
assumption, and thus by the general Hensel lemma %%(\ref{gen.hensel.lem}),
$f(y)$ has a linear factor and hence a root $\alpha\in \widehat{\o}_{p,X}$.
Then $k'=k(\alpha)$ is the required subfield.
(Note that usually one cannot find such
$k'\subset \o_{p,X}$, and so the  completion is necessary.)
\end{defn}

\begin{rem}\label{completions.rem}
By the approximation theorem of  \cite{art}, 
  $x\in X$ and $y\in Y$ are
 formally isomorphic iff there is a $z\in Z$ and
\'etale morphisms
$$
(x\in X)\leftarrow (z\in Z)\to (y\in Y).
$$
This implies that any resolution functor that commutes with
\'etale morphisms also commutes with formal
isomorphisms.

Our methods give resolution functors that
commute with  formal
isomorphisms by construction, so we do 
not need
to rely on \cite{art}.
\end{rem}

\begin{aside}[Maximal contact in  positive characteristic] 
\index{Maximal contact!fails in char.\ p}
Maximal contact, in the form presented above, 
works in
positive characteristic as long as the order of the ideal
is less than the  characteristic but 
fails in general. In some cases there is no smooth hypersurface
at all that contains the set of points where the order is maximal.
The following example is taken from \cite{na1-h}.
In characteristic 2 consider
$$
X:=(x^2+yz^3+zw^3+y^7w=0)\subset \a^4.
$$
The maximal multiplicity is 2, and the singular locus is 
given by
$$
x^2+yz^3+zw^3+y^7w=z^3+y^6w=yz^2+w^3=zw^2+y^7=0.
$$
It contains the monomial curve
$$
C:=\im[t\mapsto (t^{32}, t^7, t^{19}, t^{15})]
$$
(in fact, it is equal to it).
$C$ is not contained in any smooth hypersurface.
Indeed, assume that $(F=0)$ is a hypersurface containing
$C$ that is smooth at the origin. Then one of $x,y,z,w$
appears linearly in $F$ and $F(t^{32}, t^7, t^{19}, t^{15})\equiv 0$.
The linear term gives a nonzero
$t^m$ for some $m\in \{32,7,{19},15\}$, which must be canceled
by another term $t^m$. Thus we can write
 $m=32a+7b+{19}c+15d$, where $a+b+c+d\geq 2$
and $a,b,c,d\geq 0$. This is, however,  impossible
since none of the numbers $32,7,{19},15$ is a positive linear combination
of the other three.
\end{aside}

\section{Birational transforms and  marked ideals}

\begin{say}[Birational transform of ideals]\label{bir.transf.ideals.say}
 Let $X$ be a smooth scheme over a field $k$, $Z\subset X$ a
smooth subscheme and $\pi:B_ZX\to X$ the blow-up
with exceptional divisor $F\subset B_ZX$.
Let $Z=\cup Z_j$ and $F=\cup F_j$ be the irreducible components.

Let $I\subset \o_X$ be an ideal sheaf, and set $\ord_{Z_j}I=m_j$.
Then $\pi^*I\subset \o_{B_ZX}$ vanishes along $F_j$
with multiplicity $m_j$, and we aim to remove
the ideal sheaf $\o_{B_ZX}(-\sum m_jF_j)$ from $\pi^*I$.
That is, define the 
{\it birational transform}  (also called the 
controlled transform  or weak transform in the literature)
of  $I$ by the formula
\index{Birational!transform of an ideal sheaf}
\index{Ideal sheaf!birational transform of}
\index{Transform!birational, of an ideal sheaf}
\index{Pii@$\pi^{-1}_*I$, birational transform of an ideal sheaf}
$$
\pi^{-1}_*I:=
\o_{B_ZX}(\textstyle{\sum_j} m_jF_j)\cdot \pi^*I
\subset \o_{B_ZX}.
\eqno{(\ref{bir.transf.ideals.say}.1)}
$$
This is consistent with the definition given in
(\ref{ord.red.intro.thms}) for the case $I=I_{{\rm cod}\geq 2}$.

 Warning. If $Z\subset X$ is a smooth divisor, then
the blow-up is trivial. Hence
$\pi_{Z,X}:B_ZX\cong X$ is the identity map, and  
$$
\pi^{-1}_*I:=
\o_{X}(\textstyle{\sum_j} m_jZ_j)\cdot I
$$
depends not only on $\pi=\pi_{Z,X}$ but also on the center $Z$ of
the blow-up.
Unfortunately, I did not find any good way to fix this notational
inconsistency.
\medskip

One problem we have to deal with in resolutions is that if
$Z\subset H\subset X$ is a smooth hypersurface
with birational transform $B_ZH\subset B_ZX$ and projection $\pi_H:B_ZH\to H$,
then restriction to $H$
does not commute with taking birational transform.
That is,
$$
(\pi_H)^{-1}_*(I|_H)\supset (\pi^{-1}_*I)|_{B_ZH},
\eqno{(\ref{bir.transf.ideals.say}.2)}
$$
but  equality holds only if
$\ord_ZI=\ord_Z(I|_H)$.
\end{say}

The next definition is designed to remedy this problem.
We replace the ideal sheaf $I$ by a pair $(I,m)$,
where $m$  keeps track of the order of vanishing that we
pretend to have. The advantage is that we can redefine the
notion of birational transform to
achieve equality in (\ref{bir.transf.ideals.say}.2).

\begin{defn}\label{marked.id.defn} Let $X$ be a smooth scheme.
\index{Marked!function}\index{Marked!ideal sheaf}%
\index{Ideal sheaf!marked}\index{Im@$(I,m)$, marked ideal sheaf}%
A {\it marked function} on $X$ is a pair $(f,m)$,
where $f$ is a regular function on (some open set of) $X$ and $m$
a natural number.

A {\it marked ideal sheaf}  on $X$ is a pair $(I,m)$,
where $I\subset \o_X$ is an ideal sheaf on  $X$ and $m$
a natural number. 

The {\it cosupport} of $(I,m)$ is defined by
\index{Ideal sheaf!cosupport of}\index{Cosum@$\cosupp(I,m)$, cosupport}
$$
\cosupp(I,m):=\{x\in X: \ord_xI\geq m\}.
$$

The {\it product} of marked functions or marked ideal sheaves
\index{Product of marked ideals or functions}
is defined by
$$
(f_1,m_1)\cdot (f_2,m_2):=(f_1f_2, m_1+m_2),\quad
(I_1,m_1)\cdot (I_2,m_2):=(I_1I_2, m_1+m_2).
$$
The {\it sum} of marked functions or marked ideal sheaves
\index{Sum of marked ideals or functions}
is only sensible when the markings are the same:
$$
(f_1,m)+ (f_2,m):=(f_1+f_2, m)
\qtq{and}
(I_1,m)+ (I_2,m):=(I_1+I_2, m).
$$

The cosupport has the following elementary properties:
\begin{enumerate}
\item if $I\subset J$ then $\cosupp(I,m)\supset \cosupp(J,m)$,
\item $\cosupp(I_1I_2, m_1+m_2)\supset \cosupp(I_1,m_1)\cap \cosupp(I_2,m_2)$,
\item $\cosupp(I,m)=\cosupp(I^c,mc)$,
\item $\cosupp(I_1+I_2, m)=
 \cosupp(I_1,m)\cap \cosupp(I_2,m)$.
\end{enumerate}
\end{defn}

\begin{defn}\label{bir.tr.ids.defn}
  Let $X$ be a smooth variety, $Z\subset X$ a
smooth subvariety and $\pi:B_ZX\to X$ the blow-up
with exceptional divisor $F\subset B_ZX$.
Let $(I,m)$ be a marked ideal sheaf  on $X$ such that
 $m\leq \ord_ZI$. In analogy with
(\ref{bir.transf.ideals.say})  we define the
{\it birational transform} 
of  $(I,m)$ by the formula
\index{Birational!transform of a marked ideal sheaf}
\index{Ideal sheaf!marked, birational transform of}
\index{Transform!birational, of a marked ideal sheaf}
\index{Piim@$\pi^{-1}_*(I,m)$, birational transform of marked ideal sheaf}
$$
\pi^{-1}_*(I,m):=
\bigl(\o_{B_ZX}(mF)\cdot \pi^*I,m\bigr).
\eqno{(\ref{bir.tr.ids.defn}.1)}
$$
Informally speaking, we use the definition
(\ref{bir.transf.ideals.say}.2),
but we ``pretend that $\ord_ZI=m$.''

As in (\ref{bir.transf.ideals.say}),
it is worth calling special attention to the case
where $Z$ has codimension 1 in $X$. Then
$B_ZX\cong X$, and so scheme-theoretically there is no change.
However, the vanishing order of
$\pi^{-1}_*(I,m)$ along $Z$ is $m$ less than the
 vanishing order of
$I$ along $Z$.

In order to do  computations, choose
 local coordinates $(x_1,\dots,x_n)$ such that 
$Z=(x_1=\cdots=x_r=0)$. Then
$$
y_1=\tfrac{x_1}{x_r},\dots, y_{r-1}=
\tfrac{x_{r-1}}{x_r},\ y_r=x_r,\dots, y_n=x_n
\eqno{(\ref{bir.tr.ids.defn}.2)}
$$
give  local coordinates on a chart of $B_ZX$, and we define
$$
\pi^{-1}_*(f,m):=
\bigl(y_r^{-m}f(y_1y_r,\dots,y_{r-1}y_r,y_r,\dots,y_n),m\bigr).
\eqno{(\ref{bir.tr.ids.defn}.3)}
$$
This formula is the one we use to compute with blow-ups, but it
is coordinate system dependent. As we change coordinates, the result
of $\pi^{-1}_*$ changes by a unit. So we are free to use
$\pi^{-1}_*$ to compute the birational transform of ideal sheaves,
but one should not use it for individual functions, whose
birational transform cannot be defined (as a function).
\end{defn}

The following lemmas are easy.

\begin{lem}  \label{ored.leq.maxord.bu}
\index{Birational!transform and order of an ideal sheaf}
\index{Ideal sheaf!order after birational transform}
\index{Transform!birational, and order of an ideal sheaf}
Let $X$ be a smooth variety, $Z\subset X$ a
smooth subvariety, $\pi:B_ZX\to X$ the blow-up and
 $I\subset \o_X$  an ideal sheaf. Assume that
$\ord_ZI=\mord I$. Then
$$
\mord \pi^{-1}_*I\leq \mord I.
$$
\end{lem}

Proof. Choose local coordinates as above, and pick 
$f(x_1,\dots,x_n)\in I$ such that
$\ord_pf=\mord I=m$. Its birational transform is computed as
$$
\pi^{-1}_*f=
y_r^{-m}f(y_1y_r,\dots,y_{r-1}y_r,y_r,\dots,y_n).
$$
Since $f(x_1,\dots,x_n)$ contains a monomial  of degree $m$, the
corresponding monomial in $f(y_1y_r,\dots,y_{r-1}y_r,y_r,\dots,y_n)$
has degree $\leq 2m$, and thus in 
$\pi^{-1}_*f$ we get a monomial of degree $\leq 2m-m=m$.

This shows that $\ord_{p'}\pi^{-1}_*I\leq m$, where
$p'\in B_ZX$ denotes the origin of the chart we consider.
Performing a linear change of the $(x_1,\dots,x_r)$-coordinates
moves the origin  of the chart, and every preimage of $p$
appears as the origin after a suitable  linear change.
Thus our computation applies to all points
of the exceptional divisor of $B_ZX$.\qed

\begin{lem}  \label{bir.tr.ids.lem}
\index{Birational!transform, commutes with restriction}
\index{Ideal sheaf!marked, birational transform and restriction}
\index{Transform!birational, commutes with restriction}
Let the notation be as in (\ref{ored.leq.maxord.bu}).
Let $Z\subsetneq H\subset X$ be
 a smooth hypersurface
with birational transform $B_ZH\subset B_ZX$
and projection $\pi_H:B_ZH\to H$.
If $m\leq \ord_ZI$ and $I|_H\neq 0$, then
$$
(\pi_H)^{-1}_*(I|_H,m)= \bigl(\pi^{-1}_*(I,m)\bigr)|_{B_ZH}.
$$
\end{lem}

Proof. Again choose coordinates and assume that $H=(x_1=0)$.
Working with the chart as in (\ref{bir.tr.ids.defn}.2),
the birational transform of $H$ is  $(y_1=0)$,
and we see that it does not matter whether we set first $x_1=0$
and compute the transform or first compute the transform and
then set $y_1=0$.
We still need to contemplate  what happens
in the chart
$$
z_1={x_1}, z_{2}=\tfrac{x_{2}}{x_1},\dots, \ z_r=\tfrac{x_r}{x_1},\
z_{r+1}=x_{r+1}, \dots, z_n=x_n.
$$ 
This chart, however, does not contain any point of
the birational transform of $H$, so it does not matter.\qed
\medskip

Note that (\ref{bir.tr.ids.lem})  fails if $Z=H$.
In this case $I|_H$ is the zero ideal, $\pi_Z$ is an isomorphism,
and we have only the bad chart, which we did not need to consider 
in the proof above.
Because of this, 
we will have to consider codimension 1 subsets of $\cosupp I$
separately. 

\begin{warning}\label{bir.tr.marked.warning} 
\index{Birational!transform, warning for marked case}
\index{Ideal sheaf!marked, warning about birational transform}
\index{Transform!birational, warning about marked case}
Note that, while the  birational transform
of an ideal with $I=I_{{\rm cod}\geq 2}$ 
is  defined for an arbitrary birational morphism
(\ref{bir.transf.ideals.say}),
 we have defined the  birational transform of a marked ideal
only for a single smooth blow-up
(\ref{bir.tr.ids.defn}).
This can be extended to a sequence of smooth blow-ups,
but one has to be very careful.
Let 
$$
\Pi:X'=X_r\stackrel{\pi_{r-1}}{\longrightarrow} X_{r-1}
\stackrel{\pi_{r-2}}{\longrightarrow}\cdots
\stackrel{\pi_{1}}{\longrightarrow} X_1
\stackrel{\pi_{0}}{\longrightarrow} X_0=X
\eqno{(\ref{bir.tr.marked.warning}.1)}
$$
be a smooth  blow-up sequence.
We can  inductively define the birational transforms
of the marked ideal $(I,m)$ by
\begin{enumerate}
\item $(I_0,m):=(I,m)$, and
\item $(I_{j+1},m):=(\pi_j)^{-1}_*(I_j,m)$ as in (\ref{bir.tr.ids.defn}).
\end{enumerate}
At the end we get $(I_r,m)$, which I
rather sloppily also denote by $\Pi^{-1}_*(I,m)$.

It is very important to keep in mind
that this notation 
assumes that we have a particular blow-up sequence
in mind. That is, $\Pi^{-1}_*(I,m)$ 
depends not only on the morphism $\Pi$ but on the
actual sequence of blow-ups we use to get it.

Consider, for instance, the blow-ups
$$
\begin{array}{rcccc}
\Pi:X_2&\stackrel{\pi_{1}}{\longrightarrow}& X_{1}=B_pX_0&
\stackrel{\pi_{0}}{\longrightarrow} &X_0,\\
||\   &&&& ||\\
\Sigma:X'_2&\stackrel{\sigma_{1}}{\longrightarrow}& X'_{1}=B_CX_0&
\stackrel{\sigma_{0}}{\longrightarrow}& X_0
\end{array}
$$
introduced in (\ref{bu.dndet.ord.rem}).

Let us  compute the birational transforms of $(I,1)$, where $I:=I_C$.
The first blow-up sequence gives
$$
\begin{array}{rcl}
(\pi_0)^{-1}_*(I,1)&=&\o_{X_1}(E_0)\cdot \pi_0^*I
\qtq{and} \\
(\pi_1)^{-1}_*\bigl((\pi_0)^{-1}_*(I,1)\bigr)
&=&\o_{X_2}(E_1)\cdot \pi_1^*\bigl(
(\pi_0)^{-1}_*(I,1)\bigr)\\
&=&
\o_{X_2}(E_0+E_1)\cdot \Pi^*I.
\end{array}
$$
On the other hand, the second 
blow-up sequence gives
$$
\begin{array}{rcl}
(\sigma_0)^{-1}_*(I,1)&=&\o_{X'_1}(E'_0)\cdot \sigma_0^*I
\qtq{and} \\
(\sigma_1)^{-1}_*\bigl((\sigma_0)^{-1}_*(I,1)\bigr)&=&\o_{X'_2}(E'_1)\cdot 
\sigma_1^*\bigl((\sigma_0)^{-1}_*(I,1)\bigr)\\
&=&
\o_{X'_2}(E'_0+2E'_1)\cdot \Sigma^*I
\end{array}
$$
since $\sigma_1^*\o_{X'_1}(E'_0)=\o_{X'_2}(E'_0+E'_1)$.

Thus $\Pi^{-1}_*(I,1)\neq \Sigma^{-1}_*(I,1)$, although
$\Pi=\Sigma$.
\end{warning}

\section{The inductive setup of the proof}

In this section we set up the final notation and
state the main order reduction theorems.

\begin{notation}\label{triplet.notation}
For the rest of the chapter, $(X,I,E)$ 
or  $(X,I,m,E)$ denotes a
{\it triple}\footnote{I consider the pair $(I,m)$ as one item, so
$(X,I,m,E)$ is still a triple.},  where
\index{Triple}
\index{Xie@$(X,I,E)$}\index{Xime@$(X,I,m,E)$}
\begin{enumerate}
\item $X$ is a smooth, equidimensional (possibly reducible)
scheme of finite type over a field
 of characteristic zero,
\item $I\subset \o_X$ (resp., $(I,m)$) is a coherent  ideal sheaf
(resp.,  coherent  marked ideal sheaf),
which is nonzero on every irreducible component of $X$,  and
\item $E=(E^1,\dots,E^s)$ is an ordered set of
smooth divisors on $X$ such that $\sum E^i$ is a simple normal crossing
divisor. Each $E^i$ is allowed to be reducible or empty.
\end{enumerate}

The divisor $E$ plays an ancillary role
as a device that keeps track of the exceptional divisors that we
created and of the order in which we created them. As we saw in
(\ref{it.res.1bu.say}.3), one has to carry along
some information about the resolution process.

As we observed in (\ref{it.res.1bu.say}.2) and 
(\ref{nonproj.res.loc.glob.say}.1), it is
necessary to blow up disjoint subvarieties simultaneously.
Thus we usually do  get reducible smooth divisors  $E^j$.
\end{notation}

\begin{defn}\label{test.bu.defn}
 Given $(X,I,E)$ with $\mord I=m$, a {\it  smooth blow-up}
of order $ m$
\index{Blow-up!smooth of order $m$}
is a smooth blow-up $\pi:B_ZX\to X$ 
with {\it center}\index{Blow-up!center of} $Z$ such that
\begin{enumerate}
\item $Z\subset X$ has simple normal crossings only with $E$, and
\item $\ord_ZI= m$.
\end{enumerate}
The {\it birational transform} of  $(X,I,E)$
\index{Birational!transform of $(X,I,E)$}%
\index{Transform!birational, of $(X,I,E)$}%
under the above   blow-up
 is
$$
\pi^{-1}_*(X,I,E)=\bigl(B_ZX, \pi^{-1}_*I, \pi^{-1}_{\rm tot}(E)\bigr).
$$
Here $\pi^{-1}_*I$ is the birational transform of $I$
as defined in (\ref{bir.transf.ideals.say}),
and $\pi^{-1}_{\rm tot}(E)$ consists of the birational transform 
%%(\ref{bir.tr.var.defn})
 of
$E$ (with the same ordering as before) plus the
exceptional divisor $F\subset B_ZX$ added as the last divisor.
It is called the {\it total transform}\index{Total!transform of a divisor}
\index{Transform!total, for a divisor} of $E$.
(If $\pi$ is a trivial blow-up, then
 $\pi^{-1}_{\rm tot}(E)=E+Z$.)

 A {\it  smooth blow-up} of $(X,I,m,E)$
\index{Blow-up!smooth of order $\geq m$}
is a smooth blow-up $\pi:B_ZX\to X$ such that
\begin{enumerate}
\item[(1$'$)] $Z\subset X$ has simple normal crossings only with $E$, and
\item[(2$'$)] $\ord_ZI\geq  m$.
\end{enumerate} 
The {\it birational transform} of  $(X,I,m,E)$
\index{Birational!transform of $(X,I,m,E)$}%
\index{Transform!birational, of $(X,I,m,E)$}%
under the above   blow-up is defined as
$$
\pi^{-1}_*(X,I,m,E)=\bigl(B_ZX, \pi^{-1}_*(I,m), \pi^{-1}_{\rm tot}(E)\bigr).
$$
\end{defn}

\begin{defn}\label{seq.of.test.bu.defn}
\index{Blow-up!sequence of order $m$}
\index{Order!blow-up sequence of order $m$}
A {\it  smooth blow-up sequence} of order $ m$ and 
of length $r$ starting with $(X,I,E)$ such that $\mord I=m$ is a
smooth  blow-up sequence (\ref{bu.seq.functs})
$$
\begin{array}{l}
\Pi:(X_r,I_r,E_r)\stackrel{\pi_{r-1}}{\longrightarrow} 
(X_{r-1},I_{r-1},E_{r-1})
\stackrel{\pi_{r-2}}{\longrightarrow}\cdots\\
\hphantom{\Pi:(X_r,I_r,E_r)}
\stackrel{\pi_{1}}{\longrightarrow} (X_1,I_1,E_1)
\stackrel{\pi_{0}}{\longrightarrow} (X_0,I_0,E_0)=(X,I,E),
\end{array}
$$
where 
\begin{enumerate}
\item the $(X_i,I_i,E_i)$ are defined recursively by the
formula
$$
(X_{i+1},I_{i+1},E_{i+1}):=
\bigl(B_{Z_i}X_i, (\pi_i)^{-1}_*I_i, (\pi_i)^{-1}_{\rm tot}E_i\bigr),
$$
\item  each $\pi_i:X_{i+1}\to X_i$ is a  smooth blow-up with center
$Z_i\subset X_i$ and exceptional divisor $F_{i+1}\subset X_{i+1}$,
\item for every $i$, $Z_i\subset X_i$ has simple normal crossings
 with $E_i$, and
\item for every $i$, $\ord_{Z_i}I_i=m$.
\end{enumerate}

Similarly, a {\it  smooth blow-up sequence}
\index{Blow-up!sequence of order $\geq m$}%
\index{Order!blow-up sequence of order $\geq m$}%
 of order $\geq m$ and  
of length $r$ starting with $(X,I,m,E)$ is a
smooth  blow-up sequence
$$
\begin{array}{l}
\Pi:(X_r,I_r,m,E_r)\stackrel{\pi_{r-1}}{\longrightarrow} 
(X_{r-1},I_{r-1},m,E_{r-1})
\stackrel{\pi_{r-2}}{\longrightarrow}\cdots\\
\hphantom{\Pi:(X_r,I_r,m,E_r)}
\stackrel{\pi_{1}}{\longrightarrow} (X_1,I_1,m,E_1)
\stackrel{\pi_{0}}{\longrightarrow} (X_0,I_0,m,E_0)=(X,I,m,E),
\end{array}
$$
where
\begin{enumerate}
\item[(1$'$)] the $(X_i,I_i,m,E_i)$ are defined recursively by the
formula
$$
(X_{i+1},I_{i+1},m,E_{i+1}):=
\bigl(B_{Z_i}X_i, (\pi_i)^{-1}_*(I_i,m), (\pi_i)^{-1}_{\rm tot}E_i\bigr),
$$
\item[(2$'$--3$'$)] the sequence satisfies (2) and (3) above, and
\item[(4$'$)] for every $i$, $\ord_{Z_i}I_i\geq m$.
\end{enumerate}

As we noted in (\ref{bir.tr.ids.defn}), we allow the case where
$Z_i\subset X_i$ has codimension 1.
In this case $\pi_{i+1}$ is an isomorphism, but
$I_{i+1}\neq I_i$.

We also use the notation
$$
\begin{array}{rrl}
\Pi^{-1}_*(X,I,E)&:=&\bigl(X_r, \Pi^{-1}_*I, \Pi^{-1}_{\rm tot}(E)\bigr)\\
& := & (X_r,I_r,E_r),
\end{array}
$$
but keep in mind that, as we saw in (\ref{bir.tr.marked.warning}),
   this depends on the whole blow-up
sequence and not only on $\Pi$.

We also enrich the definition of {\it blow-up sequence functors}
considered in (\ref{busf.defnition}).
\index{Blow-up!sequence functor}%
From now on, we consider functors $\bb$ 
such that 
$\bb(X,I,E)$ (resp., $\bb(X,I,m,E)$)
is a  blow-up sequence
starting with $(X,I,E)$ (resp., $(X,I,m,E)$) as above.
That is, from now on we consider the sheaves $I_i$ and the divisors
$E_i$ as part of the functor. Since these are uniquely determined by
$(X,I,E)$ and the blow-ups $\pi_i$, this is a minor notational change.
\end{defn}

\begin{rem} The
 difference between the marked and unmarked versions
is significant, since the  birational transforms of the ideals are computed 
differently.

There is one  case, however, when one can freely pass between
the two versions. If $I$ is an ideal with $\mord I=m$, then
in any  blow-up sequence of order  $\geq m$  starting with $(X,I,m,E)$,
$\mord I_i\leq m$ by (\ref{ored.leq.maxord.bu}),
and so every  blow-up has order $=m$. Thus, by deleting $m$,
we automatically get a 
 blow-up sequence  of order $m$ starting with $(X,I,E)$.
The converse also holds.
\end{rem}

We can now state the two main
technical theorems that combine to give  an inductive proof of resolution.

\begin{thm}[Order reduction for ideals]\label{ord.red.I.thm}
\index{Order reduction!for ideal sheaves}
\index{Ideal sheaf!order reduction for}
For every $m$ there is a smooth blow-up sequence functor
$\bord_{m}$\index{BO@$\bord(X)$, order reduction functor} 
 of order $m$ (\ref{busf.defnition}) that is defined on triples 
 $(X,I,E)$  with  $\mord I\leq m$
such that if $\bord_{m}(X,I,E)=$ 
$$
\begin{array}{l}
\Pi:(X_r,I_r,E_r)\stackrel{\pi_{r-1}}{\longrightarrow} 
(X_{r-1},I_{r-1},E_{r-1})
\stackrel{\pi_{r-2}}{\longrightarrow}\cdots\\
\hphantom{\Pi:(X_r,I_r,E_r)}
\stackrel{\pi_{1}}{\longrightarrow} (X_1,I_1,E_1)
\stackrel{\pi_{0}}{\longrightarrow} (X_0,I_0,E_0)=(X,I,E),
\end{array}
$$
then
\begin{enumerate} 
\item  $\mord I_r<m$, and
\item  $\bord_{m}$ commutes with smooth morphisms
(\ref{funct.pack.proc.say}.1) 
and with change of fields (\ref{funct.pack.proc.say}.2).
\end{enumerate}
\end{thm}

In  our examples, the case  $\mord I<m$ 
is trivial, that is, 
$X_r=X$.

\begin{thm}[Order reduction for marked ideals]\label{ord.red.marked.thm}
\index{Order reduction!for marked ideal sheaves}
\index{Ideal sheaf!marked, order reduction for}
For every $m$ there is a smooth blow-up sequence functor
$\bmord_{m}$\index{BMO@$\bmord(X)$, order reduction functor}   
 of order $\geq m$ (\ref{busf.defnition}) that is defined on triples 
 $(X,I,m,E)$  
such that if $\bmord_{m}(X,I,m,E)=$ 
$$
\begin{array}{l}
\Pi:(X_r,I_r,m,E_r)\stackrel{\pi_{r-1}}{\longrightarrow} 
(X_{r-1},I_{r-1},m,E_{r-1})
\stackrel{\pi_{r-2}}{\longrightarrow}\cdots\\
\hphantom{\Pi:(X_r,I_r,m,E_r)}
\stackrel{\pi_{1}}{\longrightarrow} (X_1,I_1,m,E_1)
\stackrel{\pi_{0}}{\longrightarrow} (X_0,I_0,m,E_0),
\end{array}
$$
then
\begin{enumerate} 
\item  $\mord I_r<m$, and
\item  $\bmord_{m}$ commutes with smooth morphisms
(\ref{funct.pack.proc.say}.1) 
and also with change of fields (\ref{funct.pack.proc.say}.2).
\end{enumerate}
\end{thm}

\begin{say}[Main inductive steps of the proof]\label{main.red.2steps}

We prove (\ref{ord.red.I.thm})  and (\ref{ord.red.marked.thm}) 
together in two main reduction steps.
$$
\boxed{
\begin{array}{lcr}
&&\ \\*[-2ex]
&  \mbox{(\ref{ord.red.marked.thm}) 
in dimensions $\leq n-1$} &\\*[1ex]
\hphantom{\mbox{(\ref{main.red.2steps}.1)}}
&\Downarrow &\mbox{(\ref{main.red.2steps}.1)}\\*[1ex]
&\mbox{(\ref{ord.red.I.thm}) 
in dimension $n$}&\\*[1ex]
\hphantom{\mbox{(\ref{main.red.2steps}.2)}}
&\Downarrow &\mbox{(\ref{main.red.2steps}.2)}\\*[1ex]
&\mbox{(\ref{ord.red.marked.thm}) 
in dimension $n$}&\\*[-2ex]
&&\ 
\end{array}
}
$$
The easier part is (\ref{main.red.2steps}.2). Its proof is given in
 Section  3.13. Everything before that is devoted
to proving (\ref{main.red.2steps}.1).

We can  start the induction with the case $\dim X=0$.
Here $I=\o_X$ 
since $I$ is assumed nonzero on every irreducible component of $X$.
Everything is resolved without blow-ups.

The case $\dim X=1$ is also uninteresting.
The cosupport of an ideal sheaf is a Cartier divisor
and our algorithm tells us to blow up
$Z:=\cosupp (I,m)$. In the unmarked case $m=\mord I$. After
one blow-up $I$ is replaced by $I':=I\otimes \o_X(Z)$
which has order $<m$.
In the marked case
$\mord \bigl(I\otimes \o_X(Z)\bigr)<\mord I$. Thus, after finitely many steps,
the maximal order drops below $m$.

The 2-dimensional case is quite a bit more involved
since it includes the resolution of plane curve singularities
(essentially as in Section  1.10) and the
principalization of ideal sheaves studied in Section 1.9.
\end{say}

\begin{say}[From principalization to resolution]\label{pr.to.res.say}
As we saw in the proof of (\ref{res.sings.main.thm}),
one can prove the existence of resolutions
for quasi-projective schemes using (\ref{ord.red.I.thm}),
but it is not clear that  the 
resolution is independent of the projective embedding chosen.
In order to prove it,
we establish two additional properties of the functors
$\bord$ and $\bmord$.
\medskip

{\it Claim} \ref{pr.to.res.say}.1.  
$\bmord_{m}(X,I,m,\emptyset)=\bord_{m}(X,I,\emptyset)$
if $m=\mord I$.

\medskip

{\it Claim} \ref{pr.to.res.say}.2.  
Let  $\tau:Y\into X$ be a closed embedding of smooth schemes
 and  $J\subset \o_Y$ and $I\subset \o_X$ ideal sheaves
such that  $J$ is nonzero on every
irreducible component of $Y$ and $\tau_*(\o_Y/J)=\o_X/I$.
Then
$$
\bmord_{1}(X,I,1,\emptyset)=
\tau_*\bmord_{1}(Y,J,1,\emptyset).
$$

In both of these claims we assume that $E=\emptyset$.
One can easily extend (\ref{pr.to.res.say}.1)  to arbitrary $E$,
by slightly changing the definition (\ref{mon.nonmon.defn}). 
The situation with (\ref{pr.to.res.say}.2) is more problematic.
If $E\neq \emptyset$, then (\ref{pr.to.res.say}.2) fails 
in some cases when $\cosupp J$ contains some
irreducible components of $Y\cap E$. Most likely,
this can also be fixed with relatively minor changes, but
I do not know how.

Note also that (\ref{pr.to.res.say}.2) would not make sense
for any marking different from $m=1$.
Indeed, the ideal $I$ contains the local equations of $Y$,
thus it has order 1. Thus  $\bmord_{m}(X,I,m,\emptyset)$
is the identity for any $m\geq 2$.
\end{say}

\begin{say}[Proof of (\ref{ord.red.marked.thm}) \& (\ref{pr.to.res.say}.2)
 $\Rightarrow$ 
(\ref{funct.princ.thm.IV})]\label{pf.ored->unord}

 The only tricky point is that
in (\ref{funct.princ.thm.IV}) $E$ is a usual divisor
but  (\ref{ord.red.I.thm}) assumes that the index set of $E$ is ordered.
We can order the index set somehow, so the existence of a
principalization is not a problem. However, if we want
functoriality, then we should not introduce arbitrary choices in the
process.

If, by chance, the irreducible components of $E$ are disjoint,
then we can just declare that $E$ is a single divisor, since
in (\ref{ord.red.I.thm}) we allow the components of $E$
to be reducible. Next we show how to achieve this by
some preliminary blow-ups.

Let $E_1,\dots, E_k$ be the irreducible components of
$E$. We make the $E_i$ disjoint in $k-1$ steps.

First, let $Z_0\subset X_0=X$ be the  subset where
all of the $E_1,\dots, E_k$ intersect. 
Let $\pi_0:X_1\to X_0$ be the blow-up of $Z_0$ with exceptional
divisor $F^1$. Note that the $(\pi_0)^{-1}_*E_1,\dots, (\pi_0)^{-1}_*E_k$
do not have any $k$-fold intersections. 

Next let $Z_1\subset X_1$ be the  subset where
$k-1$ of the $(\pi_0)^{-1}_*E_i$  intersect. 
$Z_1$ is smooth since the  $(\pi_0)^{-1}_*E_i$
do not have any $k$-fold intersections. 
Let $\pi_1:X_2\to X_1$ be the blow-up of $Z_1$ with exceptional
divisor $F^2$. Note that the 
$(\pi_0\pi_1)^{-1}_*E_i$
do not have any $(k-1)$-fold intersections. 

  Next let $Z_2\subset X_2$ be the  subset where
$k-2$ of the $(\pi_0\pi_1)^{-1}_*E_i$
  intersect, and so on.

After $(k-1)$-steps we get rid of all pairwise intersections as well.
The end result is
$\pi:X'\to X$ such that $E^0:=\pi^{-1}_*(E_1+\cdots+E_k)$
is a smooth divisor. Let $E^1,\dots,E^{k-1}$ denote
the birational transforms of $F^1,\dots,F^{k-1}$.

Thus $(X', \pi^*I, \sum_{i=0}^{k-1}E^i)$ satisfies the assumptions
of (\ref{ord.red.I.thm}).

(It may seem  natural to start with $\dim X$-fold intersections
instead of $k$-fold intersections. We want 
functoriality with respect to all smooth morphisms, so we should not
use the dimension of $X$ in constructing the
resolution process. However, ultimately the difference is only in some
empty blow ups, and we can forget about those at the end.) 

The rest is straightforward.
Construct
$$
\Pi_{(X,I,E)}: X_r\to \dots \to 
X_s=X'\to \dots \to X
$$
 by composing
$\bmord_{1}(X', \pi^*I, 1,\sum_{i=0}^{k-1}E^i)$
with $X'\to X$.
By construction,
$\Pi_{(X,I,E)}^*I=I_{r}\cdot\o_{X_r}(F)$ 
for some effective divisor $F$ supported on 
the total transform of 
$\sum_{i=0}^{k-1}E^i$. Here $I_r=\o_{X_r}$ since $\mord I_r<1$ and
 $F$ is a simple normal crossing divisor. Therefore
 $\Pi_{(X,I,E)}^*I$ is a  monomial ideal which
 can be written down explicitly as follows.

Let $F_j\subset X_{j+1}$ denote the exceptional divisor
of the $j$th step in the above smooth blow-up sequence
for $\Pi_{(X,I,E)}:X_r\to X$. Then
$$
\Pi_{(X,I,E)}^*I=\o_{X_r}
\bigl(-\textstyle{\sum_{j=s}^r} \Pi_{r,j+1}^*F_j\bigr),
$$
where $\Pi_{r,j+1}:X_r\to X_{j+1}$ is the 
corresponding composite of blow-ups.

The functoriality properties required in
 (\ref{funct.princ.thm.IV}) 
follow from the corresponding
functoriality properties in (\ref{ord.red.marked.thm})
and from  (\ref{pr.to.res.say}.2).
\qed
\end{say}

\section{Birational transform of  derivatives}

\begin{defn}[Derivative of an ideal sheaf]
\label{der.of.ids.defn}

On a smooth scheme $X$ over a field $k$,  let $\der_X$ denote the sheaf
of derivations $\o_X\to \o_X$.
\index{Derivation}
\index{Derx@$\der_X$, sheaf of derivations}
If $x_1,\dots,x_n$ are local coordinates at a point $p\in X$, then
the derivations $\partial/\partial x_1, \dots, \partial/\partial x_n$
are local generators of $\der_X$.
Derivation gives a $k$-bilinear map
$$
\der_X\times\ \o_X\to \o_X.
$$

Let $I\subset \o_X$ be an ideal sheaf.
Its {\it first derivative} is the ideal sheaf
$D(I)$ generated by  all derivatives
\index{Derivative!ideal}\index{Ideal sheaf!derivative}%
\index{Di@$D(I)$, derivative ideal sheaf}%
of elements of $I$. That is,
$$
D(I):=\bigl( \im[\der_X\times\ I\to \o_X]\bigr).
\eqno{(\ref{der.of.ids.defn}.1)}
$$
Note that $I\subset D(I)$, as shown by the formula
$$
f=\frac{\partial (xf)}{\partial x}-x\frac{\partial f}{\partial x}.
$$
In terms of generators we can write $D(I)$ as
$$
D(f_1,\dots,f_s)=
\Bigl( f_i, \frac{\partial f_i}{\partial x_j}: 
1\leq i\leq s, 1\leq j\leq n\Bigr).
$$
Higher derivatives are defined inductively by
\index{Derivative!ideal, higher}
\index{Dri@$D^r(I)$, higher derivative ideal}
$$
D^{r+1}(I):=D\bigl(D^r(I)\bigr).\eqno{(\ref{der.of.ids.defn}.2)}
$$
(Note that  $D^r(I)$ contains all $r$th partial derivatives of
elements of $I$, but over  general rings  it is bigger;
try second derivatives over $\z[x]$.
Over characteristic zero fields, they are actually equal,
as one can see using formulas like
$$
\frac{\partial f}{\partial y}=
\frac{\partial^2 (xf)}{\partial y\partial x}-
x\frac{\partial^2 f}{\partial y\partial x}\qtq{and}
2\frac{\partial f}{\partial x}=
\frac{\partial^2 (xf)}{\partial x^2}-
x\frac{\partial^2 f}{\partial x^2}.
$$
The inductive definition is easier to work with.)

If $\mord I\leq m$, then $D^m(I)=\o_X$, and thus
 the $D^r(I)$ give an ascending chain of ideal sheaves
$$
I\subset D(I)\subset D^2(I)\subset \cdots \subset D^m(I)=\o_X.
$$
This is, however, not the right way to look at  derivatives.
Since differentiating a function $r$ times is expected to
reduce its order by $r$, we define the
derivative of a  marked ideal by
\index{Derivative!ideal, marked}
\index{Dim@$D(I,m)$, derivative ideal sheaf}
\index{Drim@$D^r(I,m)$, higher derivative ideal}
$$
D^r(I,m):=\bigl(D^r(I), m-r\bigr)\qtq{for $r\leq m$.}\eqno{(\ref{der.of.ids.defn}.3)}
$$
Before we can
usefully compare the ideal $I$ and its higher derivatives,
we have to correct for the difference in their markings.
\end{defn}

Higher derivatives have the usual properties.

\begin{lem}\label{der.of.ids.lem}
Let the notation be as above. Then
\begin{enumerate}
\item $D^r(D^s(I))=D^{r+s}(I)$,
\item
$D^r(I\cdot J)\subset
\sum_{i=0}^r D^i(I)\cdot D^{r-i}(J)$ (product rule),
\item $\cosupp(I,m)=\cosupp (D^r(I),m-r)$  for $r<m$ (char.\ $0$ only!),
\item if $f:Y\to X$ is smooth, then $D(f^*I)=f^*(D(I))$,
\item $D\bigl( \hat I\, \bigr)=\widehat{D(I)}$, 
where $\hat{\ }$ denotes completion (\ref{completions.def}).\qed
\end{enumerate}
\end{lem}

\ref{der.of.ids.lem}.6
(Aside about positive characteristic). 
The above  definition of higher derivatives is ``correct'' only
in characteristic zero. In general, one should use the
{\it Hasse-Dieudonn\'e derivatives,} which are essentially
\index{Derivative!Hasse-Dieudonn\'e}
given by
$$
\frac{1}{r_1!\cdots r_n!}\cdot \frac{\partial^{\sum r_i}}
{\partial x_1^{r_1}\cdots \partial x_n^{r_n}}.
$$
These operators then have other problems. One of the main difficulties
of resolution in positive characteristic is a lack of good
replacement for higher derivatives.
\medskip

\begin{say}[Birational transform of derivatives]
\label{bir.trans.ders.say}
  Let $X$ be a smooth variety, $Z\subset X$ a
smooth subvariety and $\pi:B_ZX\to X$ the blow-up
with exceptional divisor $F\subset B_ZX$.
Let $(I,m)$ be a marked ideal sheaf  on $X$ such that
 $m\leq \ord_ZI$. Choose
 local coordinates $(x_1,\dots,x_n)$ such that 
$Z=(x_1=\cdots=x_r=0)$.  Then
$$
y_1=\tfrac{x_1}{x_r},\dots, y_{r-1}=\tfrac{x_{r-1}}{x_r},\ 
y_r=x_r,\dots, y_n=x_n
$$
are  local coordinates on a chart of $B_ZX$.
Let us compute the derivatives of 
$$
\pi^{-1}_*\bigl(f(x_1,\dots,x_n),m\bigr)=
\bigl(y_r^{-m}f(y_1y_r,\dots,y_{r-1}y_r,y_r,\dots,y_n),m\bigr),
$$
defined in (\ref{bir.tr.ids.defn}.3).
The easy formulas are
$$
\begin{array}{lll}
\frac{\partial}{\partial y_j}\pi^{-1}_*(f,m)&=&
\pi^{-1}_*\bigl(\frac{\partial}{\partial x_j}f,m-1)
\qtq{for $j<r$,}\\ 
\frac{\partial}{\partial y_j}\pi^{-1}_*(f,m)&=&
\frac1{y_r}
\pi^{-1}_*\bigl(\frac{\partial}{\partial x_j}f,m-1)\qtq{for $j>r$,}
\end{array}
$$
and a more complicated one using the chain rule
for $j=r$:
$$
\begin{array}{lll}
\frac{\partial}{\partial y_r}\pi^{-1}_*(f,m)&=&
\tfrac{y_i}{y_r}\sum_{i<r} \pi^{-1}_*\bigl(\frac{\partial}{\partial x_i}f,m-1)+
\frac1{y_r}
\pi^{-1}_*\bigl(\frac{\partial}{\partial x_r}f,m-1)\\*[1ex]
&& +(\frac{-m}{y_r},-1)\cdot\pi^{-1}_*(f,m),
\end{array}
$$
where, as in (\ref{marked.id.defn}),
 multiplying by $(\frac{-m}{y_r},-1)$ means multiplying the
function by $\frac{-m}{y_r}$ and lowering the marking by $1$.

These can be rearranged to
$$
\begin{array}{lllr}
\pi^{-1}_*\bigl(\frac{\partial}{\partial x_j}f,m-1)&\!\!\!\!=\!\!\!\!&
\frac{\partial}{\partial y_j}\pi^{-1}_*(f,m)
\qtq{for $j<r$,}
&(\ref{bir.trans.ders.say}.1)\\*[1ex]
\pi^{-1}_*\bigl(\frac{\partial}{\partial x_j}f,m-1)
&\!\!\!\!=\!\!\!\!&y_r\frac{\partial}{\partial y_j}\pi^{-1}_*(f,m)
\qtq{for $j>r$,}
& (\ref{bir.trans.ders.say}.2)\\*[1ex]
\pi^{-1}_*\bigl(\tfrac{\partial}{\partial x_r}f,m-1)&\!\!\!\!=\!\!\!\!&
y_r\tfrac{\partial}{\partial y_r}\pi^{-1}_*(f,m)-
y_r\sum_{i<r} \tfrac{\partial}{\partial y_i}\pi^{-1}_*(f,m)\!\!\!\!&\\
&& +(m,-1)\cdot\pi^{-1}_*(f,m).
& (\ref{bir.trans.ders.say}.3)
\end{array}
$$
For later purposes, also note the following version of
(\ref{bir.trans.ders.say}.1):
$$
\begin{array}{lllr}
\pi^{-1}_*\bigl(x_j\frac{\partial}{\partial x_j}f,m-1)&=&
y_ry_j\frac{\partial}{\partial y_j}\pi^{-1}_*(f,m)
\qtq{for $j<r$.}
& (\ref{bir.trans.ders.say}.4)
\end{array}
$$
Observe that the right-hand sides of these equations are in 
$D(\pi^{-1}_*(f,m))$.  Thus we have proved
 the following elementary but important statement.
\end{say}

\begin{thm} \label{bir.trans.ders.thm}
\index{Birational!transform of a derivative ideal sheaf}
\index{Ideal sheaf!derivative, birational transform of}
\index{Transform!birational, of a derivative ideal sheaf}
Let $(I,m)$ be a marked ideal and
$\Pi:X_r\to X$ the composite of a  smooth blow-up sequence
of order $\geq m$  
starting with $(X,I,m)$. Then 
$$
\Pi^{-1}_*\bigl(D^j(I,m)\bigr)\subset D^j\bigl(\Pi^{-1}_*(I,m)\bigr)
\qtq{for every $j\geq 0$.}
$$
\end{thm}

Proof. For $j=1$ and for one blow-up this 
is what
the above formulas
(\ref{bir.trans.ders.say}.1--3) say. 
The rest follows by  induction on $j$ and on the number of
blow-ups.\qed

\begin{cor}\label{test.bu.seq.comp.cor}
Let 
$$
\begin{array}{l}
\Pi:(X_r,I_r,m)\stackrel{\pi_{r-1}}{\longrightarrow} 
(X_{r-1},I_{r-1},m)
\stackrel{\pi_{r-2}}{\longrightarrow}\cdots\\
\hphantom{\Pi:(X_r,I_r,m)}
\stackrel{\pi_{1}}{\longrightarrow} (X_1,I_1,m)
\stackrel{\pi_{0}}{\longrightarrow} (X_0,I_0,m)
\end{array}
$$
be a smooth  blow-up sequence of order $\geq m$ starting with $(X,I,m)$.

Fix  $j\leq m$, and define inductively the ideal sheaves
$J_i$ by 
$$
 J_0:=D^j(I)\qtq{and}
(J_{i+1},m-j):=(\pi_i)^{-1}_*(J_i,m-j).
$$
Then, 
 $J_i\subset D^j(I_i)$ for every $i$, and
 we get a  smooth blow-up sequence of order $\geq m-j$
starting with $(X,D^j(I),m-j)$
$$
\begin{array}{l}
\Pi:(X_r,J_r,m-j)\stackrel{\pi_{r-1}}{\longrightarrow} 
(X_{r-1},J_{r-1},m-j)
\stackrel{\pi_{r-2}}{\longrightarrow}\cdots\\
\hphantom{\Pi:(X_r,J_r,m-j)}
\stackrel{\pi_{1}}{\longrightarrow} (X_1,J_1,m-j)
\stackrel{\pi_{0}}{\longrightarrow} (X_0,J_0,m-j).
\end{array}
$$
\end{cor}

Proof. We need to check that for every $i<r$ the inequality
$\ord_{Z_i}J_i\geq m-j$ holds, where $Z_i\subset X_i$ is the center of
the blow-up $\pi_i:X_{i+1}\to X_i$.
If $\Pi_i:X_i\to X$ is the composition, then
$$
J_i=(\Pi_i)^{-1}_*(D^jI,m-j)\subset D^j\bigl( (\Pi_i)^{-1}_*(I,m)\bigr)
=D^j(I_i,m),
$$
where the containment in the middle follows from
(\ref{bir.trans.ders.thm}).
By assumption $\ord_{Z_i}I_i\geq m$, and thus
$\ord_{Z_i}D^j(I_i)\geq m-j$ by (\ref{der.of.ids.lem}.3).\qed
\medskip

\section{Maximal contact and going down}\label{max.cont.sect}
\label{prob3.sec}

\begin{defn}\label{max.cont.defn}
Let $X$ be a smooth variety, $I\subset \o_X$ an ideal sheaf
and $m=\mord I$.
A smooth hypersurface $H\subset X$ is called a hypersurface of
{\it maximal contact} if the following holds.

For 
every open set $X^0\subset X$ and for
every  smooth blow-up sequence of order $ m$ 
starting with $(X^0,I^0:=I|_{X^0})$,
$$
\Pi:(X^0_r,I^0_r)\stackrel{\pi_{r-1}}{\longrightarrow} 
(X^0_{r-1},I^0_{r-1})
\stackrel{\pi_{r-2}}{\longrightarrow}\cdots
\stackrel{\pi_{1}}{\longrightarrow} (X^0_1,I^0_1)
\stackrel{\pi_{0}}{\longrightarrow} (X^0_0,I^0_0),
$$
the center of every blow-up $Z^0_i\subset X^0_i$
is contained in the birational transform $H^0_i\subset X^0_i$ of
$H^0:=H\cap X^0$.
This implies that
$$
\begin{array}{l}
\Pi|_{H^0_r}:(H^0_r,I_r|_{H^0_r},m)\stackrel{\pi_{r-1}}{\longrightarrow} 
(H^0_{r-1},I_{r-1}|_{H^0_{r-1}},m)
\stackrel{\pi_{r-2}}{\longrightarrow}\cdots\\
\hphantom{\Pi|_{H^0_r}:(H^0_r,I_r|_{H^0_r},m)}
\stackrel{\pi_{1}}{\longrightarrow} (H^0_1,I_1|_{H^0_1},m)
\stackrel{\pi_{0}}{\longrightarrow} (H^0_0,I_0|_{H^0_0},m)
\end{array}
$$
is a  smooth blow-up sequence of order $\geq m$ 
starting with $(H^0,I|_{H^0},m)$.

Being a hypersurface of
 maximal contact is a local property.

For now we ignore the divisorial part  $E$ of a triple $(X,I,E)$
since we
cannot guarantee that $E|_H$ is also a
simple normal crossing divisor. 
\end{defn}

\begin{defn}\label{MC.id.defn}
Let $X$ be a smooth variety, $I\subset \o_X$ an ideal sheaf
and  $m=\mord I$. The {\it maximal contact ideal} of $I$
is
$$
MC(I):=D^{m-1}(I).
$$
\index{Maximal contact!ideal sheaf of}
\index{Ideal sheaf!maximal contact}
\index{MCI@$MC(I)$, ideal of maximal contacts}
Note that $MC(I)$ has order 1 at $x\in X$ 
if $\ord_xI=m$ and order 0 if $\ord_xI<m$. 
Thus 
$$
\cosupp MC(I)=\cosupp (I,m).
$$
\end{defn}

\begin{thm}[Maximal contact] \label{max.contact.thm.I}
\index{Maximal contact!local existence}
Let $X$ be a smooth variety, $I\subset \o_X$ an ideal sheaf
and  $m=\mord I$.  
Let $L$ be a line bundle on $X$ and $h\in H^0(X,L\otimes MC(I))$ 
  a section with zero divisor $H:=(h=0)$.
\begin{enumerate}
\item 
If $H$ is smooth and $I|_H\neq 0$, then $H$ is a hypersurface of
 maximal contact.
\item Every $x\in X$ has an open neighborhood $x\in U_x\subset X$
and  $h_x\in H^0(U_x,L\otimes MC(I))$ 
 such that  $H_x:=(h_x=0)\subset U_x$ is smooth.
\end{enumerate}
\end{thm}

Proof. Being a hypersurface of
 maximal contact is a local question, and thus 
we may assume that $L=\o_X$.
Let
$$
\Pi:(X_r,I_r)\stackrel{\pi_{r-1}}{\longrightarrow} 
(X_{r-1},I_{r-1})
\stackrel{\pi_{r-2}}{\longrightarrow}\cdots
\stackrel{\pi_{1}}{\longrightarrow} (X_1,I_1)
\stackrel{\pi_{0}}{\longrightarrow} (X_0,I_0)
$$
be a  smooth blow-up sequence of order $m$  starting with $(X,I)$,
where $\pi_i$ is the blow-up of $Z_i\subset X_i$.

Applying (\ref{test.bu.seq.comp.cor}) 
for $j=m-1$, we obtain a 
 smooth blow-up sequence of order $\geq 1$ starting with $(X,J_0:=MC(I),1)$:
$$
\Pi:(X_r,J_r,1)\stackrel{\pi_{r-1}}{\longrightarrow} 
(X_{r-1},J_{r-1},1)
\stackrel{\pi_{r-2}}{\longrightarrow}\cdots
\stackrel{\pi_{1}}{\longrightarrow} (X_1,J_1,1)
\stackrel{\pi_{0}}{\longrightarrow} (X_0,J_0,1).
$$

Let $H_i:=(\Pi_i)^{-1}_*H\subset X_i$ denote the birational
transform of $H\subset X$.
Since $\o_{X_0}(-H_0)\subset J_0$ and $H_0$ is smooth,
we see that $\o_{X_i}(-H_i)\subset J_i$ for every $i$.
By assumption $\ord_{Z_i}I_i\geq m$. Thus, using
(\ref{der.of.ids.lem}.3) and (\ref{test.bu.seq.comp.cor}) we get that
$$
\ord_{Z_i}J_i\geq \ord_{Z_i}MC(I_i)\geq 1
$$
and hence also $\ord_{Z_i}H_i\geq 1$.
 Thus 
 $Z_i\subset H_i$ for every $i$, and so
$H$ is a hypersurface of
 maximal contact.

To see the second claim, 
pick $x\in X$ such that $\ord_xI=m$. Then
 $\ord_xMC(I)=1$ by 
(\ref{der.of.ids.lem}.3).
 Thus there is a local section of $MC(I)$ that has
order 1 at $x$, and so its zero divisor is smooth in a neighborhood of $x$.
\qed

\begin{aside} A section  $h\in MC(I)$ such that
$H=(h=0)$ is smooth
  always exists locally
but usually not globally, not even if we tensor $I$ by a very ample line bundle
$L$. By the Bertini-type theorem of
\cite[4.4]{kol-pairs},  the best one can achieve globally is that
$H$ has $cA$-type singularities.
(These are given by local  equations 
$x_1x_2+(\mbox{other terms})=0$.)
\end{aside}

The above results say that every 
smooth  blow-up sequence of order $m$ starting with $(X,I)$
can be seen as a smooth   blow-up sequence starting with $(H,I|_H,m)$.

An important remaining  problem is that
not every  smooth blow-up sequence starting with $(H,I|_H,m)$
corresponds to a smooth 
 blow-up sequence of order $ m$ starting with $(X,I)$,
and thus we cannot yet construct an order reduction of
 $(X,I)$ from an  order reduction of $(H,I|_H,m)$.

Here are some examples that show what can go wrong.

\begin{exmp}\label{mc.badrestr.exmp}

Let $I=(xy-z^n)$. Then $\ord_0I=2$ and
$D(I)=(x,y,z^{n-1})$. $H:=(x=0)$ is a
surface of maximal contact, and
$$
(H, I|_{H})\cong \bigl(\a^2_{y,z},(z^n)\bigr).
$$
Thus  $ (H, I|_{H})$  shows a 1-dimensional singular
locus of order $n$, whereas we have an isolated singular point
of order 2. The same happens if we use $(y=0)$ as a
surface of maximal contact.

In this case we do better if we use a general surface of maximal contact.
Indeed, for
$H_g:=(x-y=0)$, 
$$
(H_g, I|_{H_g})\cong \bigl(\a^2_{x,z},(x^2-z^n)\bigr),
$$
and we get an equivalence between smooth
 blow-up sequences of order $2$ starting with $(\a^3,(xy-z^n) )$
and  smooth 
blow-up sequences  of order $\geq 2$ starting with $(\a^2,(x^2-z^n),2 )$.

In some cases, even the general hypersurface of maximal contact
fails to produce an equivalence.
There are no problems on $H$ itself, but difficulties appear after
blow-ups.

Let $I=(x^3+xy^5+z^4)$. 
A general surface of maximal contact
is 
$$
H:=(x+u_1xy^3+u_2y^4+u_3z^2=0),
\qtq{where the $u_i$ are units.}
$$
Let us compute two blow-ups given by
$x_1=x/y, y_1=y, z_1=z/y$ and $x_2=x_1/y_1, y_2=y_1, z_2=z_1/y_1$.
We get the birational transforms
$$
\begin{array}{ll}
x^3+xy^5\hphantom{_1}+z^4 & x\hphantom{_1}+u_1xy^3\hphantom{_1}+u_2y^4+u_3z^2\\
x_1^3+x_1y_1^3+y_1z_1^4 & x_1+u_1x_1y_1^3+u_2y_1^3+u_3y_1z_1^2\\
x_2^3+x_2y_2+y_2^2z_2^4\qquad & x_2+u_1x_2y_2^3+u_2y_2^2+u_3y_2^2z_2^2.
\end{array}
$$
The second birational transform of the ideal
has order 2 on this chart.
However, its restriction to the  birational transform $H_2$ of $H$ 
still has order 3
since we can use the equation of $H_2$ to eliminate $x_2$ by the substitution
$$
x_2=-y_2^2(u_2+u_3z_2^2)(1+u_1y_2^3)^{-1}
$$
to obtain that $I_2|_{H_2}\subset (y_2^3, y_2^2z_2^4)$. 
\end{exmp}

\section{Restriction of derivatives and going up}
\label{prob4.sec}

In general, neither the order of an ideal nor its derivative ideal
commute with restrictions to smooth hypersurfaces.
For instance, if $I=(x^2+xy+z^3)$  and $S=(x=0)$
then $\ord_0 I=2$ but $\ord_0 (I|_S)=3$
and $(DI)|_S=(y,z^2)$ but $D(I|_S)=(z^2)$.
It is easy to see that 
$$
\ord_p I\leq \ord_p (I|_S)\qtq{and}
(DI)|_S\supset D(I|_S),
$$
but neither is an equality.
The notion of $D$-balanced ideals provides a
solution to the first of these problems and a partial remedy to the second.

\begin{defn} \label{weak.D.bal} As in  (\ref{D-balance.say}), 
\index{D-bal@$D$-balanced ideal sheaf}
\index{Ideal sheaf!D-bal@$D$-balanced}
an ideal $I$ with $m=\mord I$ is called {\it $D$-balanced} if
$$
\bigl(D^iI\bigr)^m\subset I^{m-i}\quad \forall\ i<m.
$$
\end{defn}

If $I$ is  $D$-balanced, then at every point it has order either $m$
or $0$. Indeed, if $\ord_pI<m$ then $(D^{m-1}I)_p=\o_{p,X}$, thus
$I^{m-1}$ and  $I$ both contain a unit at $p$.
In particular, $\cosupp (I,m)=\cosupp I$, hence
 the maximal order commutes
with  restrictions.

We can reformulate this observation as follows.
If $I$ is $D$-balanced, then
 any smooth blow-up
of order $\geq m$ for $I|_S$ corresponds to a
smooth blow-up
of order $\geq m$ for $I$.  We would like a similar statement not
just for one blow-up, but for all blow-up sequences.

\begin{thm}[Going-up property of  $D$-balanced ideals]
\label{D-bal.goup.thm}
\index{Going!up}
\index{D-bal@$D$-balanced ideal sheaf!and going up}
Let $X$ be a smooth variety and $I$ a  $D$-balanced sheaf of ideals
with $m=\mord I$. Let $S\subset X$ be any smooth hypersurface
such that $S\not\subset\cosupp(I,m)$ and
$$
\begin{array}{l}
\Pi^S:(S_r,J_r,m)\stackrel{\pi^S_{r-1}}{\longrightarrow} 
(S_{r-1},J_{r-1},m)
\stackrel{\pi^S_{r-2}}{\longrightarrow}\cdots\\
\hphantom{\Pi^S:(S_r,J_r,m)}
\stackrel{\pi^S_{1}}{\longrightarrow} (S_1,J_1,m)
\stackrel{\pi^S_{0}}{\longrightarrow} (S_0,J_0,m)=(S,I|_S,m)
\end{array}
$$
be a smooth blow-up sequence of order $\geq m$,
where $\pi^S_i$ is the blow-up of
$Z_i\subset S_i$. 
Then the pushed-forward  sequence (\ref{bu.seq.functs})
$$
\begin{array}{l}
\Pi:(X_r,I_r)\stackrel{\pi_{r-1}}{\longrightarrow} 
(X_{r-1},I_{r-1})
\stackrel{\pi_{r-2}}{\longrightarrow}\cdots\\
\hphantom{\Pi:(X_r,I_r)}
\stackrel{\pi_{1}}{\longrightarrow} (X_1,I_1)
\stackrel{\pi_{0}}{\longrightarrow} (X_0,I_0)=(X,I)
\end{array}
$$
is a smooth blow-up sequence of order $m$, where $\pi_i$ is the blow-up of
$Z_i\subset S_i\subset X_i$.
\end{thm}

\begin{cor}[Going up and down] \label{D-bal.maxcont.cor}
\index{Going!up and down}
\index{D-bal@$D$-balanced ideal sheaf!and going up and down}
Let $X$ be a smooth variety, $I\subset \o_X$ a  $D$-balanced
 ideal sheaf with
$m=\mord I$ and $E$ a divisor with simple normal crossings.
Let $H\subset X$ be  a smooth hypersurface of maximal
contact  such that $E+H$ is also a 
divisor with simple normal crossings
and  no irreducible component of $H$ is contained in $\cosupp(I,m)$.

Then pushing forward   (\ref{bu.seq.functs}) from $H$ to $X$
is a one-to-one correspondence between
\begin{enumerate}
\item smooth blow-up sequences of order $\geq m$ 
starting with the triple $(H,I|_H,m,E|_H)$, and
\item smooth blow-up sequences of order $ m$ 
starting with $(X,I,E)$.
\end{enumerate}
\end{cor}

Proof. This follows from (\ref{D-bal.goup.thm})
   and (\ref{max.contact.thm.I}),
except for the role played by $E$. 

Adding $E$ to $(X,I)$ (resp., to $(H,m,I|_H)$) 
means that now we can use only blow-ups whose centers
are in simple normal crossing with $E$ (resp., $E|_H$)
and their total transforms.
Since $E|_H$ is again a 
divisor with simple normal crossings, this poses the same restriction on
order reduction for $(X,I,E)$ as on
order reduction for $(H,I|_H,m,E|_H)$.\qed
\medskip

\begin{say}[First attempt to prove (\ref{D-bal.goup.thm})]
\label{1st.attempt.D-bal.goup.thm}
We already noted that we are ok for the first blow-up.
Let us see what happens with pushing forward the second
blow-up $\pi^S_1:S_2\to S_1$.
By assumption $\ord_{Z_1}(\pi_0^S)^{-1}_*(I|_S,m)\geq m$.
Can we conclude from this that
$\ord_{Z_1}(\pi_0)^{-1}_*(I,m)\geq m$? In other words, is
$$
S_1\cap \cosupp \bigl((\pi_0)^{-1}_*(I,m)\bigr)=
\cosupp \bigl(\pi_0^S)^{-1}_*(I|_S,m)\bigr)?
$$
Since the birational transform commutes with restrictions,
this indeed holds if the birational transform  $(\pi_0)^{-1}_*(I,m)$
 is again $D$-balanced.
By assumption $(D^i I)^m\subset I^{m-i}$
and so 
$$
\bigl((\pi_0)^{-1}_*(D^i I,m-i)\bigr)^m\subset 
\bigl((\pi_0)^{-1}_*(I,m)\bigr)^{m-i}.
$$
Unfortunately, when we interchange
$(\pi_0)^{-1}_*$ and $D^i$ on the left-hand side,
the inequality in (\ref{bir.trans.ders.thm})
goes the wrong way and indeed, in general
the birational transform is not $D$-balanced. 

Looking at the formulas (\ref{bir.trans.ders.say}.1--3), we see that
taking birational transform commutes with some 
derivatives 
but not with  others.

In order to exploit this, we introduce logarithmic derivatives.
This notion enables us to separate the  ``good'' directions from the
``bad'' ones. 
\medskip

{\it Example} \ref{1st.attempt.D-bal.goup.thm}.1. 
 Check that $(x^2, xy^m, y^{m+1})$
is $D$-balanced. After blowing up the origin,
one of the charts gives  $(x_1^2, x_1y_1^{m-1}, y_1^{m-1})$,
which is not  $D$-balanced. 
\end{say}

\begin{say}[Logarithmic derivatives]\label{log-der.say}
Let $X$ be a smooth variety and $S\subset X$ a smooth subvariety.
For simplicity, we assume that $S$ is a hypersurface.
At a point $p\in S$ pick local coordinates
$x_1,\dots,x_n$ such that $S=(x_1=0)$.
If $f$ is any function, then 
$$
\frac{\partial f}{\partial x_i}|_S=
\frac{\partial \bigl(f|_S\bigr)}{\partial x_i}
\qtq{for $i>1$,}
$$
but $\partial (f|_S)/{\partial x_1}$
does not even make sense.
Therefore, we would like to decompose $D(f)$ 
into two parts:
\begin{enumerate}
\item[$\bullet$] ${\partial f}/{\partial x_i}$ for $i>1$
(these commute with restriction to $S$), and
\item[$\bullet$]
 $\partial f/{\partial x_1}$
(which does not).
\end{enumerate}
\noindent Such a decomposition is, however,  not coordinate invariant.
The best one can do is the following.

Let $\der_X(-\log S)\subset \der_X$
be the largest subsheaf that maps $\o_X(-S)$ into itself by
derivations.  It is called the
sheaf of {\it logarithmic derivations} along $S$.
\index{Derivative!logarithmic}
\index{Derxs@$\der_X(-\log S)$, logarithmic derivations}
 In the above  local coordinates we can write
$$
\der_X(-\log S)=\Bigl( x_1\frac{\partial}{\partial x_1},
\frac{\partial}{\partial x_2},\dots, 
\frac{\partial}{\partial x_n}\Bigr).
$$
For an ideal sheaf $I$ set
\index{Dersi@$D(-\log S)(I)$, logarithmic derivative of an ideal sheaf}
\index{Logarithmic derivative}
$$
\begin{array}{rcl}
D(-\log S)(I)&:=&\bigl( \im[\der_X(-\log S)\times\ I\to \o_X]\bigr)
\qtq{and}\\
D^{r+1}(-\log S)(I)&:=&D(-\log S)\bigl(D^r(-\log S)(I)\bigr)
\qtq{for $r\geq 1$.}
\end{array}
$$
We need three properties of log derivations.

First,  log derivations behave well with respect to restriction to $S$:
$$
\bigl(D^r(-\log S)(I)\bigr)|_S=D^r(I|_S).
\eqno{(\ref{log-der.say}.1)}
$$

Second,  one can  filter the sheaf $D^s(I)$ by subsheaves
$$
D^s(-\log S)(I)\subset D^{s-1}(-\log S)\bigl(D(I)\bigr)
\subset\cdots\subset D^s(I).
$$
There are no well-defined complements, but in
 local coordinates $x_1,\dots,x_n$ we can write
$$
D^s(I)=D^s(-\log S)(I)+
D^{s-1}(-\log S)\Bigl(\frac{\partial I}{\partial x_1}\Bigr)
+\cdots +
\Bigl(\frac{\partial^s I}{\partial x_1^s}\Bigr),
\eqno{(\ref{log-der.say}.2)}
$$
and the first $j+1$ summands span $D^{s-j}(-\log S)(D^j(I))$.

Third,  under the assumptions of 
(\ref{D-bal.goup.thm}), 
we get a logarithmic version of (\ref{bir.trans.ders.thm}):
$$
\Pi^{-1}_*\Bigl(D^j(-\log S_r)(I,m)\Bigr)\subset 
D^j(-\log S)\Bigl(\Pi^{-1}_*(I,m)\Bigr),
\eqno{(\ref{log-der.say}.3)}
$$
which is proved the same way using (\ref{bir.trans.ders.say}.4).

\end{say}

We can now formulate the next result, which
 can be viewed as a way to reverse the
inclusion in (\ref{bir.trans.ders.thm}).

\begin{thm}\label{D.and.bu.restr.prop}
\index{Birational!transform and log derivative}%
\index{Transform!birational, and log derivative}%
\index{Logarithmic derivative!and birational transform}%
Consider a smooth  blow-up sequence 
of order $\geq m$:  
$$
\begin{array}{l}
\Pi:(X_r,I_r,m)\stackrel{\pi_{r-1}}{\longrightarrow} 
(X_{r-1},I_{r-1},m)
\stackrel{\pi_{r-2}}{\longrightarrow}\cdots\\
\hphantom{\Pi:(X_r,I_r,m)}
\stackrel{\pi_{1}}{\longrightarrow} (X_1,I_1,m)
\stackrel{\pi_{0}}{\longrightarrow} (X_0,I_0,m)=(X,I,m).
\end{array}
$$
Let $S\subset X$ be a smooth hypersurface and $S_i\subset X_i$ its
birational transforms.  Assume that  each blow-up  center
$Z_i$ is contained in $S_i$.Then
$$
D^s\Pi^{-1}_*(I,m)=\sum_{j=0}^{s}
D^{s-j}(-\log S_r)\Pi^{-1}_*\bigl(D^jI,m-j\bigr).
\eqno{(\ref{D.and.bu.restr.prop}.1)}
$$
\end{thm}

Proof.  Using (\ref{bir.trans.ders.thm}) we obtain that
$$
\begin{array}{rcl}
D^{s-j}(-\log S_r)\Pi^{-1}_*\bigl(D^jI,m-j\bigr)\!\!
&\subset&\!\! D^{s-j}\Pi^{-1}_*\bigl(D^jI,m-j\bigr)\\
&\subset&\!\! D^{s-j} D^j (\Pi^{-1}_*I,m)=D^s\Pi^{-1}_*(I,m),
\end{array}
$$
and thus the right-hand side of (\ref{D.and.bu.restr.prop}.1)
is contained in the left-hand side.

Next let us check  the reverse inclusion in (\ref{D.and.bu.restr.prop}.1)
  for one blow-up.
The question is local on $X$, and so
choose coordinates
$x_1,\dots,x_n$   such that $S=(x_1=0)$ and
the center of the blow-up $\pi$ is $(x_1=\cdots=x_r=0)$.
We have a typical local chart
$$
y_1=\tfrac{x_1}{x_r},\dots, y_{r-1}=\tfrac{x_{r-1}}{x_r},\ y_r=x_r,\dots,
 y_n=x_n,
$$
and $S_1=(y_1=0)$ is the birational transform of $S$.
Note that the 
blow-up is covered by $r$ different charts, but only $r-1$ of these
can be written in the above forms, where $x_r$ is different from $x_1$.
These  $r-1$
charts, however,  completely  cover 
$S_1$.

Applying (\ref{log-der.say}.2)
 to $\pi^{-1}_*(I,m)$
we obtain that
$$
D^s \pi^{-1}_*(I,m)=
\sum_{j=0}^sD^{s-j}(-\log S_1)
\textstyle{\left(\frac{\partial^j \pi^{-1}_*(I,m)}{\partial y_1^j}
\right)}.
$$
Although usually differentiation does not commute with
birational transforms, by (\ref{bir.trans.ders.say}.1)
it does so for $\partial/ \partial x_1$ and $\partial/ \partial y_1$. So
we can rewrite the above formula as
$$
\begin{array}{rcl}
D^s \pi^{-1}_*(I,m)&=&
\sum_{j=0}^sD^{s-j}(-\log S_1)\pi^{-1}_*
\left(\frac{\partial^j (I,m)}{\partial x_1^j}\right)\\*[1ex]
&\subset & 
\sum_{j=0}^sD^{s-j}(-\log S_1)\pi^{-1}_*(D^jI,m-j),
\end{array}
\eqno{(\ref{D.and.bu.restr.prop}.2)}
$$
where the inclusion is clear. As noted above, the right-hand side of
(\ref{D.and.bu.restr.prop}.2) is contained in the left-hand side,
and hence they are equal. This proves
(\ref{D.and.bu.restr.prop}) for one blow-up.

In the general case, we use induction on the number of blow-ups.
We factor $\Pi:X_r\to X$ as the composite of $\pi_{r-1}:X_r\to X_{r-1}$
and $\Pi_{r-1}:X_{r-1}\to X$.
Use (\ref{D.and.bu.restr.prop}) for $\pi_{r-1}$
to get that
$$
\begin{array}{l}
D^s\Pi^{-1}_*(I,m)
 = D^s(\pi_{r-1})^{-1}_*(\Pi_{r-1})^{-1}_*(I,m)\\*[1ex]
\hphantom{D^s\Pi^{-1}_*(I,m)} 
= \sum_{j=0}^{s} D^{s-j}(-\log S_r)(\pi_{r-1})^{-1}_*
D^j(\Pi_{r-1})^{-1}_*(I,m).
\end{array}
\eqno{(\ref{D.and.bu.restr.prop}.3)}
$$
By induction (\ref{D.and.bu.restr.prop}) holds  for $\Pi_{r-1}$
 and $s=j$, thus
$$
D^j(\Pi_{r-1})^{-1}_*(I,m)
= 
\sum_{\ell=0}^{j} D^{j-\ell}(-\log S_{r-1})(\Pi_{r-1})^{-1}_*
\bigl(D^{\ell}I,m-\ell\bigr).
$$
By (\ref{log-der.say}.3),  we can  interchange $(\pi_{r-1})^{-1}_*$
and $D^{j-\ell}(-\log S_{r-1})$, and so
$$
\begin{array}{l}
(\pi_{r-1})^{-1}_*D^j(\Pi_{r-1})^{-1}_*(I,m)\\*[1ex]
\hphantom{D^s\Pi^{-1}_*(I,m)} =(\pi_{r-1})^{-1}_*
\sum_{\ell=0}^{j} D^{j-\ell}(-\log S_{r-1})(\Pi_{r-1})^{-1}_*
\bigl(D^{\ell}I,m-\ell\bigr)\\*[1ex]
\hphantom{D^s\Pi^{-1}_*(I,m)} \subset
\sum_{\ell=0}^{j}  D^{j-\ell}(-\log S_r)
(\pi_{r-1})^{-1}_*(\Pi_{r-1})^{-1}_*
\bigl(D^{\ell}I,m-\ell\bigr)\\*[1ex]
\hphantom{D^s\Pi^{-1}_*(I,m)} = 
\sum_{\ell=0}^{j}  D^{j-\ell}(-\log S_r)
\Pi^{-1}_*
\bigl(D^{\ell}I,m-\ell\bigr).
\end{array}
$$
Substituting into (\ref{D.and.bu.restr.prop}.3), we obtain the desired result:
$$
\begin{array}{l}
D^s\Pi^{-1}_*(I,m)\\
\quad\quad \subset
 \sum_{j=0}^{s}  D^{s-j}(-\log S_r) \sum_{\ell=0}^{j}D^{j-\ell}(-\log S_r)
\Pi^{-1}_*
\bigl(D^{\ell}I,m-\ell\bigr)\\*[1ex]
\quad\quad  = \sum_{\ell=0}^{s }D^{s-\ell}(-\log S_r) 
\Pi^{-1}_* \bigl(D^{\ell}I,m-\ell\bigr). \qed
\end{array}
$$

\begin{cor}\label{D.and.bu.restr.thm}
\index{Birational!transform, cosupport of}%
\index{Ideal sheaf!cosupport of birational transform}%
\index{Transform!birational, cosupport of}%
Let the notation and assumptions be as in (\ref{D.and.bu.restr.prop}).
 Then
$$
S_r\cap \cosupp\bigl(\Pi^{-1}_*(I,m)\bigr)=
\bigcap_{j=0}^{m-1}
\cosupp(\Pi|_{S_r})^{-1}_*\bigl((D^jI)|_S,m-j\bigr).
$$
\end{cor}

Proof.
Restricting (\ref{D.and.bu.restr.prop}.1) 
to $S_r$ and  using (\ref{bir.tr.ids.lem}) and (\ref{log-der.say}.1) 
we get that
$$
\bigl(D^s\Pi^{-1}_*(I,m)\bigr)|_{S_r}=\sum_{j=0}^{s}
D^{s-j}(\Pi|_{S_r})^{-1}_*\bigl((D^jI)|_S,m-j\bigr).
\eqno{(\ref{D.and.bu.restr.thm}.1)}
$$
Set $s=m-1$ and 
take cosupports.
Since $D^{m-1}\Pi^{-1}_*(I,m)$ has order 1, its cosupport
commutes with restrictions, so the left-hand side of
(\ref{D.and.bu.restr.thm}.1) becomes
$$
\begin{array}{rcl}
\cosupp\Bigl(\bigl(D^{m-1}\Pi^{-1}_*(I,m)\bigr)|_{S_r}\Bigr)&=&
\cosupp\bigl(D^{m-1}\Pi^{-1}_*(I,m)\bigr)\cap S_r\\
&=&\cosupp\bigl(\Pi^{-1}_*(I,m)\bigr)\cap S_r,
\end{array}
\eqno{(\ref{D.and.bu.restr.thm}.2)}
$$
where the second equality follows from 
 (\ref{der.of.ids.lem}.3).

On the right-hand side of (\ref{D.and.bu.restr.thm}.1)
  use (\ref{marked.id.defn}.4)
to obtain that 
$$
\begin{array}{l}
\cosupp\Bigl(\sum_{j=0}^{m-1}
D^{m-1-j}(\Pi|_{S_r})^{-1}_*\bigl((D^jI)|_S,m-j\bigr)\Bigr)\\
\quad =\bigcap_{j=0}^{m-1}
\cosupp\Bigl(D^{m-1-j}(\Pi|_{S_r})^{-1}_*\bigl((D^jI)|_S,m-j\bigr)\Bigr)\\
\quad =\bigcap_{j=0}^{m-1}
\cosupp(\Pi|_{S_r})^{-1}_*\bigl((D^jI)|_S,m-j\bigr).\\
\end{array}
\eqno{(\ref{D.and.bu.restr.thm}.3)}
$$
The last lines of (\ref{D.and.bu.restr.thm}.2) and
(\ref{D.and.bu.restr.thm}.3) are thus equal.
\qed

\begin{say}[Proof of (\ref{D-bal.goup.thm})]
By induction, assume that this already holds 
for blow-up sequences of length
$<r$. We need to show that the last blow-up
also has order $\geq m$, or, equivalently,
$ \cosupp (J_{r-1},m)\subset \cosupp (I_{r-1},m)$.

Using first (\ref{D.and.bu.restr.thm}) for $\Pi_{r-1}:X_{r-1}\to X$,
 then  the 
 $D$-balanced property in line 2, we obtain that
$$
\begin{array}{rcl}
S_{r-1}\cap \cosupp(I_{r-1},m)\!\!&=&\!\!
\bigcap_{j=0}^{m-1}
\cosupp(\Pi_{r-1}^S)^{-1}_*\bigl((D^jI)|_S,m-j\bigr)\\*[1ex]
&=&\!\!
\bigcap_{j=0}^{m-1}
\cosupp(\Pi_{r-1}^S)^{-1}_*\bigl((D^jI)^m|_S,m(m-j)\bigr)\\*[1ex]
&\supset &\!\! 
\bigcap_{j=0}^{m-1}
\cosupp(\Pi_{r-1}^S)^{-1}_*\bigl(I^{m-j}|_S,m(m-j)\bigr)\\*[1ex]
&= & \!\!
\bigcap_{j=0}^{m-1}
\cosupp(\Pi_{r-1}^S)^{-1}_*\bigl(I|_S,m\bigr)\\
&=&\!\!
\cosupp(\Pi_{r-1}^S)^{-1}_*\bigl(J_0,m\bigr)\\
&=&\!\!\cosupp (J_{r-1},m). \qed
\end{array}
$$
\end{say}

\section{Uniqueness of maximal contact}
\label{prob5.sec}

Given $(X,I,E)$, let $j:H\into X$ and $j':H'\into X$
 be two hypersurfaces of maximal contact.
By (\ref{D-bal.maxcont.cor})
 we can construct smooth blow-up sequences
for  $(X,I,E)$ from $(H,I|_H,m,E_H)$ and also from
$(H',I|_{H'},m,E_{H'})$. We need to guarantee that we get the same
  blow-up sequences.

Assume that there is  an automorphism $\phi$ of
$X$ such that $\phi^*I=I$ and $\phi^{-1}(E+H')=E+H$.
Then
$(H,I_H,m,E_H)=\phi^*(H',I|_{H'},m,E_{H'})$,
thus if $\bb(H',I|_{H'},m,E_{H'})$
is the smooth blow-up sequence constructed using
$H'$, then the ``same''
construction using  $H$ gives
$$
\bb(H,I_H,m,E_H)=\phi^*\bb(H',I|_{H'},m,E_{H'}).
$$
Pushing these forward as in (\ref{D-bal.maxcont.cor}),
we obtain that 
$$
j_*\bb(H,I_H,m,E_H)=\phi^*\bigl(j'_*\bb(H',I|_{H'},m,E_{H'})\bigr).
$$
That is, the blow-up sequences we get from 
$H$ and $H'$ are isomorphic, but we would like them
to be identical.

Let $Z_0$ (resp., $Z'_0$) be the center of the first blow-up
obtained using $H$ (resp., $H'$). 
As above, 
 $\phi^{-1}(Z'_0)=Z_0$. 
Both $Z'_0$ and $Z_0$ are contained in $\cosupp(I,m)$,
so  if
$\phi$ is the identity on $\cosupp(I,m)$ then
 $Z'_0=Z_0$.

The assumption $\phi^*I=I$ implies that $\phi$ maps $\cosupp(I,m)$
into itself, but it does not imply that 
$\phi$ is the identity on $\cosupp(I,m)$.
How can we  achieve the latter?
Let $R$ be a ring, $J\subset R$ an ideal and $\sigma$ an automorphism of $R$.
It is easy to see that  $\sigma(J)=J$ and $\sigma$ induces the identity 
automorphism on $R/J$
iff   $r-\sigma(r)\in J$ for every $r\in R$.

How should we choose the ideal $J$ in our situation?
It turns out that $J=I$ does not work and 
the ideal sheaf of $\cosupp(I,m)$ behaves badly for blow-ups.
An intermediate choice is given by $D^{m-1}(I)=MC(I)$,
which works well.

Another twist is  that usually $X$ itself has no automorphisms
(not even Zariski locally), so we have to work
in a formal or  \'etale neighborhood of  a point $x\in X$.
(See (\ref{completions.def}) for completions.)

\begin{defn}
\label{maxcont.formal.equiv.defn}
Let $X$ be a smooth variety, $p\in X$ a point,
$I$  an  ideal sheaf such that $\mord I=\ord_pI=m$
 and $E=E^1+\cdots +E^s$  a simple normal crossing divisor.
Let $H,H'\subset X$ be two hypersurfaces of
maximal contact.

We say that  $H$ and $H'$ are
{\it formally equivalent at $p$} with respect to $(X,I,E)$ if
 there is an automorphism 
 $\phi:\hat X\to \hat X$
which moves $(X,I,H+E)$ into $(X,I,H'+E)$ and $\phi$ is close to the identity.
That is, 
\begin{enumerate}
\item $\phi (\hat H)=\hat H'$,
\item $\phi^* (\hat I)=\hat I$, 
\item $\phi (\hat E^i)=\hat E^i$ for $i=1,\dots,s$, and
\item $h-\phi^*(h)\in MC(\hat I)$ for every $h\in \hat{\o}_{x,X}$.
\end{enumerate}

While this is the important concept, it is somewhat inconvenient
to use since we defined resolution, 
order reduction, and so on for schemes of finite type and not
for general schemes like $\hat X$.

Even very simple formal automorphisms cannot be realized as
algebraic automorphisms on some \'etale cover.
(Check this for the map $x\mapsto \sqrt{x}$, which is a 
formal automorphism of $(1\in \hat{\c})$.)
Thus we need a slightly modified definition.

We say that $H$ and $H'$ are 
 {\it \'etale equivalent}
\index{etale@\'Etale equivalent!hypersurfaces of maximal contact}%
with respect to $(X,I,E)$ if
 there are \'etale surjections
 $\psi,\psi':U\rightrightarrows X$
such that
\begin{enumerate}
\item[(1$'$)] $\psi^{-1} ( H)= \psi'^{-1}(H')$,
\item[(2$'$)] $\psi^*  (I)= \psi'^*(I)$, 
\item[(3$'$)]  
 $\psi^{-1} (E^i)= \psi'^{-1}(E^i)$ for $i=1,\dots,s$, and
\item[(4$'$)] $\psi^*(h)-\psi'^*(h)\in MC\bigl(\psi^*  (I)\bigr)$
 for every $h\in \o_X$.
\end{enumerate}

The connection with the formal case comes from noting
that $\psi$ is invertible after completion, and
then $\phi:=\hat{\psi'}\circ \hat{\psi}^{-1}:\hat X\to \hat X$
is the automorphism we seek.

\end{defn}

A key observation of \cite{wlo}
is that for certain ideals $I$ any two
smooth hypersurfaces of
maximal contact are formal and \'etale equivalent.
 Recall 
(\ref{mc-homog.say})
 that an ideal $I$ is 
 MC-invariant if
$$
MC(I)\cdot D(I)\subset I,
$$
where $MC(I)$ is the ideal of maximal contacts defined in
(\ref{max.cont.intr.say}.2).
Since taking derivatives commutes with completion
(\ref{der.of.ids.lem}.5), we see that
$\widehat{MC(I)}=MC(\hat I\,)$.

\begin{thm}[Uniqueness of maximal contact]
\label{maxcont.formal.automs.thm}
\index{Maximal contact!formal uniqueness}
\index{Maximal contact!etale@\'etale uniqueness}
Let $X$ be a smooth variety over a field of characteristic zero,
 $I$  an MC-invariant ideal sheaf,  $m=\mord I$ 
and $E$ a simple normal crossing divisor.
Let  $H,H'\subset X$ be two smooth hypersurfaces of
maximal contact for $I$ such that 
$H+E$ and $H'+E$ both have simple normal crossings.

Then  $ H$ and $H'$ are  \'etale  equivalent 
with respect to  $(X, I, E)$.
\end{thm}

We start with a general result relating automorphisms and
derivations of complete local rings.
 Since derivations are essentially the first order
automorphisms, it is reasonable to expect that an ideal is invariant under
a subgroup of automorphisms iff it is invariant to first order.
We are, however, in an infinite-dimensional setting, so it is safer
to work out the details.

\begin{notation} Let $k$ be a field of characteristic zero,
$K/k$  a finite field extension
and $R=K[[x_1,\dots,x_n]]$ the formal power series ring
in $n$ variables with maximal ideal $m$, viewed as a $k$-algebra.
 For $g_1,\dots, g_n\in m$ the
map  $g:x_i\mapsto g_i$ extends to an automorphism of $R$
$\Leftrightarrow$    $g:m/m^2\to m/m^2$ is an isomorphism
$\Leftrightarrow$  the linear parts of the $g_i$ are
linearly independent.

Let $B\subset m$ be an ideal. For $b_i\in B$ the map
  $g:x_i\mapsto x_i+b_i$ need not generate  an automorphism,
but
$g:x_i\mapsto x_i+\lambda_ib_i$
gives an automorphism for general $\lambda_i\in k$.
We call these {\it automorphisms of the form}
${\mathbf 1}+B$.\index{Automorphism of the form
${\mathbf 1}+B$}
\end{notation}

\begin{prop}\label{aut-inv=der-inv}
 Let the notation be as above, and let $I\subset R$ be an ideal.
The following are equivalent:
\begin{enumerate}
\item $I$ is invariant under every automorphism of the form
${\mathbf 1}+B$,
\item $B\cdot D(I)\subset I$,
\item $B^j\cdot D^j(I)\subset I$ for every $j\geq 1$.
\end{enumerate}
\end{prop}

Proof. Assume that $B^j\cdot D^j(I)\subset I$ for every $j\geq 1$.
Given any $f\in I$, we need to prove that
$f(x_1+b_1,\dots,x_n+b_n)\in I$.
Take the  Taylor expansion 
$$
f(x_1+b_1,\dots,x_n+b_n)=
f(x_1,\dots,x_n)+\sum_i b_i\frac{\partial f}{\partial x_i}
+\frac12\sum_{i,j} b_ib_j\frac{\partial^2 f}{\partial x_i\partial x_j}+\cdots. 
$$
For any $s\geq 1$, this gives that
$$
f(x_1+b_1,\dots,x_n+b_n)\in 
I+B\cdot D(I)+\cdots+ B^s\cdot D^s(I)+ m^{s+1}\subset I+m^{s+1}
$$
since $B^j\cdot D^j(I)\subset I$ by assumption.
Letting $s$ go to infinity, 
by Krull's intersection theorem (\ref{completions.def}) we conclude that
$f(x_1+b_1,\dots,x_n+b_n)\in I$.

Conversely, for any $b\in B$ and general $\lambda_i\in k$,
invariance under the automorphism
$(x_1,x_2,\dots,x_n)\mapsto (x_1+\lambda_ib,x_2,\dots,x_n)$
gives that
$$
f(x_1+\lambda_ib,x_2,\dots,x_n)=
f(x_1,\dots,x_n)+\lambda_ib\frac{\partial f}{\partial x_1}+
\cdots (\lambda_ib)^s \frac{\partial f^s}{\partial x_1^s}\in I+m^{s+1}.
$$
Use $s$ different values  $\lambda_1,\dots,\lambda_s$.
 Since the Vandermonde determinant $(\lambda_i^j)$ is invertible,
we conclude that
$$
b \frac{\partial f}{\partial x_1}\in I+m^{s+1}.
$$
Letting $s$ go to infinity, we obtain that  $B\cdot D(I)\subset I$.

Finally, we prove by induction that 
$B^j\cdot D^j(I)\subset I$ for every $j\geq 1$.
$B^{j+1}\cdot D^{j+1}(I)$ is generated by elements of the form
$b_0\cdots b_j\cdot D(g)$, where $g\in D^j(I)$.
The product  rule gives that
$$
\begin{array}{rcl}
b_0\cdots b_j\cdot D(g)&=&b_0\cdot D(b_1\cdots b_j\cdot g)
-\sum_{i\geq 1} D(b_i)\cdot(b_0\cdots \widehat{b_i}\cdots b_j\cdot g)\\
&\in & B\cdot D\bigl(B^j\cdot D^j(I)\bigr)
+B^j\cdot D^j(I)\\
&\subset& B\cdot D(I)+B^j\cdot D^j(I)\subset I,
\end{array}
$$
where the entry $\widehat{b_i}$ is omitted from the products.
\qed

\begin{say}[Proof of (\ref{maxcont.formal.automs.thm})]
Let us start with  formal equivalence.

Pick local sections $x_1,x'_1\in MC(I)$ such that
 $H=(x_1=0)$ and $H'=(x'_1=0)$. Choose other local coordinates
$x_2,\dots,x_{s+1}$ at $p$ such that
$E^i=(x_{i+1}=0)$ for $i=1,\dots,s$.  For a general choice
of $x_{s+2},\dots,x_n$, we see that 
$x_1,x_2,\dots,x_n$ and $x'_1,x_2,\dots,x_n$ are both local
coordinate systems.

If $X$ is a $k$-variety and the residue field
of $p\in X$ is $K$, then 
$\widehat{\o}_{p,X}\cong K[[x_1,\dots,x_n]]$ by (\ref{completions.def}),
so the computations of (\ref{aut-inv=der-inv}) apply.

Since $x_1-x'_1\in MC(I)$, the automorphism
$$
\phi^*(x'_1,x_2,\dots,x_n)= 
\bigl(x'_1+(x_1-x'_1),x_2,\dots,x_n\bigr)=(x_1,x_2,\dots,x_n)
$$
is of the form ${\mathbf 1}+MC(I)$. Hence by
(\ref{aut-inv=der-inv}) we conclude that
$\phi^*\bigl(\hat I\bigr)=\hat I$. By construction $\phi(\hat H)=\hat H'$,
 $\phi(\hat E^i)=\hat E^i$ and
(\ref{maxcont.formal.equiv.defn}.4) is also clear.

In order to go from the formal to the \'etale case,
the key point is to realize the automorphism $\phi$
on some \'etale neighborhood.
 Existence follows  from the general approximation theorems
of \cite{art}, but in our case the choice is clear.

Take $X\times X$, and for some $p\in X$ let
$x_{11},x_{12},\dots,x_{1n}$  be the
corresponding local coordinates on the first factor and
$x'_{21},x_{22},\dots,x_{2n}$   on the second factor.
Set
$$
U_1(p):=(x_{11}-x'_{21}=x_{12}-x_{22}=\cdots=x_{1n}-x_{2n}=0)
\subset X\times X.
$$
The completion of $U_1(p)$ at $(p,p)$
is the graph of $\phi_p$. By shrinking $U_1(p)$, we get 
$(p,p)\in U_2(p)\subset U_1(p)$
such that  both coordinate projections 
$\psi_p,\psi'_p:U_2(p)\rightrightarrows X$ are \'etale.

From our previous considerations, we know that 
(\ref{maxcont.formal.equiv.defn}.1--4)
 hold after taking completions
at $(p,p)$. Thus (\ref{maxcont.formal.equiv.defn}.1$'$--4$'$)
 also hold in an open neighborhood $U(p)\ni (p,p)$
 by (\ref{completions.def}).

The images of finitely many of the $U(p)$ cover $X$.
We can take $U$ to be their disjoint union.
\qed
\end{say}

In Section 3.12 we use the maximal contact hypersurfaces
$H,H'$ to construct blow-up sequences ${\mathbf B}$ and ${\mathbf B'}$
which become isomorphic after pulling back to $U$.
The next result shows that they are
the same already on $X$. That is, our blow-ups
do not depend on the choice of a hypersurface of maximal contact.

\begin{defn} \label{et.eq.bu.defn}
 Let $X$ be a smooth variety over a field of characteristic zero and let
$$
\begin{array}{lcl}
{\mathbf B}&:=&(X_r,I_r)\stackrel{\pi_{r-1}}{\longrightarrow} 
\cdots \stackrel{\pi_{0}}{\longrightarrow} (X_0,I_0)=(X,I),\qtq{and}\\
{\mathbf B'}&:=&(X'_r,I'_r)\stackrel{\pi'_{r-1}}{\longrightarrow} 
\cdots \stackrel{\pi'_{0}}{\longrightarrow} (X'_0,I'_0)=(X,I)
\end{array}
$$
be two blow-up sequences  of order $m=\mord I$.
We say that ${\mathbf B}$ and ${\mathbf B'}$
are {\it \'etale equivalent}\index{Blow-up!sequences, \'etale equivalent}%
\index{etale@\'Etale equivalent!blow-up sequences}
  if 
 there are \'etale surjections
$\psi,{\psi'}:U\rightrightarrows X$ such that
\begin{enumerate}
\item  $\psi^*  (I)= {\psi'}^*(I)$,
\item  $\psi^*(h)-\psi'^*(h)\in MC\bigl(\psi^*  (I)\bigr)$
 for every $h\in \o_X$, and
\item $\psi^*{\mathbf B}={\psi'}^*{\mathbf B'}$.
\end{enumerate}
\end{defn}

\begin{thm}\label{et.eq=eq}
 Let $X$ be a smooth variety over a field of characteristic zero and
 $I$ an MC-invariant  ideal sheaf.
Let ${\mathbf B}$ and ${\mathbf B'}$ be two
blow-up sequences  of order $m=\mord I$
which are \'etale equivalent.

Then ${\mathbf B}={\mathbf B'}$.
\end{thm}

Proof. By assumption there are \'etale surjections
$\psi,{\psi'}:U\rightrightarrows X$ such that
$\psi^*{\mathbf B}={\psi'}^*{\mathbf B'}$.
Let 
$$
{\mathbf B}^U:=(U_r,I^U_r)\stackrel{\pi^U_{r-1}}{\longrightarrow} 
\cdots \stackrel{\pi^U_{0}}{\longrightarrow} (U_0,I^U_0)=(U,\psi^*I={\psi'}^*I)
$$
be the common pullback.
We prove by induction on $i$ the following claims.
\begin{enumerate}
\item  $(X_i, I_i)=(X'_i, I'_i)$.
\item  $\psi, {\psi'}$ lift to \'etale surjections
$\psi_i,{\psi'}_i:U_i\rightrightarrows X_i$ such that
$$
\im(\psi_i^*-{\psi'_i}^*)\subset (\Pi_i^U)^{-1}_*\bigl(MC(I^U_0),1\bigr).
$$
\item $Z_{i-1}=Z'_{i-1}$.
\end{enumerate}

For $i=0$ there is nothing to prove. Let us see how to go from
$i$ to $i+1$.
Set $W_i=\cosupp (\Pi_i^U)^{-1}_*\bigl(MC(I^U_0),1\bigr)$ and note that
$Z^U_i\subset W_i$ by (\ref{test.bu.seq.comp.cor}).
By the inductive assumption (2)
$\psi_i|_{W_i}= {\psi'_i}|_{W_i}$, thus
$Z_i=\psi_i(Z^U_i)={\psi'_i}(Z^U_i)=Z'_i$.
This in turn implies that
$X_{i+1}=X'_{i+1}$.

In order to compute the lifting of $\psi_i$ and $\psi'_i$,
 pick local coordinates
 $x_1,\dots,x_n$ on $X_i=X'_i$ such that
$Z_i=Z'_i=(x_1=x_2=\cdots=x_k=0)$.
By induction, 
$$
{\psi'_i}^*(x_j)=\psi_i^*(x_j)-b(i,j)
\qtq{for some $b(i,j)\in (\Pi_i^U)^{-1}_*\bigl(MC(I^U),1\bigr)$.}
$$
On the blow-up $\pi_i:X_{i+1}\to X_i$ consider the local chart
$$
y_1=\tfrac{x_1}{x_r},\dots, y_{r-1}=\tfrac{x_{r-1}}{x_r},\ y_r=x_r,\dots,
 y_n=x_n.
$$
We need to prove that 
$$
\psi_{i+1}^*(y_j)-{\psi'_{i+1}}^*(y_j)
\in (\Pi_{i+1}^U)^{-1}_*\bigl(MC(I^U),1\bigr)
$$
for every $j$. This is clear if $y_j=x_j$, that is, for $j\geq r$.
Next we compute the case when $j<r$.

The $b(i,j)$ vanish along $Z^U_i$ and so
 $(\pi^U_i)^*b(i,j)=\psi_{i+1}^*(x_r)b(i+1,j)$ for some
$b(i+1,j)\in (\Pi_{i+1}^U)^{-1}_*\bigl(MC(I^U),1\bigr)$.
Hence, for $j<r$, we obtain that
$$
\begin{array}{rcl}
{\psi'_{i+1}}^*(y_j)&=&(\pi_i^U)^*
\displaystyle{\frac{{\psi'_i}^*(x_j)}{{\psi'_i}^*(x_r)}}=
(\pi_i^U)^*
\displaystyle{\frac{\psi_i^*(x_j)-b(i,j)}{\psi_i^*(x_r)-b(i,r)}}\\*[2ex]
&=&
\displaystyle{\frac{\psi_{i+1}^*(x_j)-\psi_{i+1}^*(x_r)b(i+1,j)}
{\psi_{i+1}^*(x_r)-\psi_{i+1}^*(x_r)b(i+1,r)}}\\*[2ex]
&=&
\displaystyle{\frac{\psi_{i+1}^*(y_j)-b(i+1,j)}{1-b(i+1,r)}}.
\end{array}
$$
This implies that
$$
\psi_{i+1}^*(y_j)-{\psi'_{i+1}}^*(y_j)=
\displaystyle{\frac{b(i+1,j)-b(i+1,r)\psi_{i+1}^*(y_j)}
{1-b(i+1,r)}}
$$
is in $ (\Pi_{i+1}^U)^{-1}_*\bigl(MC(I^U_0),1\bigr)$,
as required.\qed

\section{Tuning of ideals}
\label{tune.sec}

Following  (\ref{maxcont.key.exmp}.9) and 
(\ref{test.bu.seq.comp.cor}),
we are looking for  ideals that contain information about
all derivatives of $I$ with equalized markings.

\begin{defn}[Maximal coefficient ideals]\label{W_s(I).defn}
\index{Ws@$W_s(I)$, maximal coefficient ideal}
\index{Coefficient!ideal, maximal}
 Let $X$ be a smooth variety,
$I\subset \o_X$ an ideal sheaf and $m=\mord I$.
The   {\it maximal  coefficient ideal} of order $s$ of $I$ is
$$
W_s(I):=
\left(
\prod_{j=0}^{m}  \bigl(D^{j}(I)\bigr)^{c_j} : \sum (m-j)c_j\geq s
\right)\subset \o_X.
$$
\end{defn}

The ideals $W_s(I)$ satisfy a series of useful properties.

\begin{prop}\label{w_s.props}
 Let $X$ be a smooth variety,
$I\subset \o_X$ an ideal sheaf and $m=\mord I$.
Then 
\begin{enumerate}
\item $W_{s+1}(I)\subset W_{s}(I)$ for every $s$,
\item $W_s(I)\cdot W_t(I)\subset W_{s+t}(I)$, 
\item $D(W_{s+1}(I))= W_s(I)$,
\item $MC(W_s(I))= W_1(I)=MC(I)$, 
\item $W_s(I)$ is MC-invariant,
\item $W_{s}(I)\cdot W_t(I)=W_{s+t}(I)$ whenever 
 $t\geq (m-1)\cdot
\lcm(2,\dots,m)$ and $s=r\cdot\lcm(2,\dots,m)$,
\item $\bigl(W_{s}(I)\bigr)^j=W_{js}(I)$ whenever 
$s=r\cdot\lcm(2,\dots,m)$ for some $r\geq m-1$, and
\item $W_{s}(I)$ is $D$-balanced whenever 
$s=r\cdot\lcm(2,\dots,m)$ for some $r\geq m-1$. 
\end{enumerate}
\end{prop}

Proof. Assertions (1) and (2) are clear and
 the inclusion $D(W_{s+1}(I))\subset W_s(I)$ 
 follows from the product rule.
Pick $x_1\in MC(I)$ that has  order 1 at $p$.
Then $x_1^{s+1}\in W_{s+1}(I)$, implying
$$
x_1^{s}=(s+1)^{-1}\tfrac{\partial}{\partial x_1}x_1^{s+1}\in
 D\bigl(W_{s+1}(I)\bigr).
$$
Next we prove by induction on $t$ that
$x_1^{s-t}W_t\subset  D\bigl(W_{s+1}(I)\bigr)$, which gives (3) for $t=s$.

Note that $x_1^{s+1-t}f\in W_{s+1}(I)$ for any $f\in W_t(I)$. Thus
$$
(s+1-t)x_1^{s-t}f+x_1^{s+1-t}\bigl( \tfrac{\partial}{\partial x_1}f\bigr)=
\tfrac{\partial}{\partial x_1}\bigl(x_1^{s+1-t}f\bigr)\in
 D\bigl(W_{s+1}(I)\bigr).
$$
Since $\tfrac{\partial}{\partial x_1}f\in W_{t-1}(I)$, then  by induction, 
$x_1^{s+1-t}\bigl( \tfrac{\partial}{\partial x_1}f\bigr)\in 
 D\bigl(W_{s+1}(I)\bigr)$. Hence also
$x_1^{s-t}f\in 
 D\bigl(W_{s+1}(I)\bigr)$.

Applying (3) repeatedly gives that
$MC(W_s(I))= W_1(I)$,
which in turn contains $D^{m-1}(I)=MC(I)$ by definition.
Conversely,  $W_1(I)$ is generated by products of derivatives,
at least one of which is a derivative of order $<m$.
Thus 
$$
W_1(I)\subset \sum_{j<m} D^j(I)=D^{m-1}(I),
$$
proving (4). Together with (2) and (3), this implies (5).

Thinking of elements of $D^{m-j}(I)$ as variables
of degree $j$, (6) is implied by (\ref{w_s.props}.9),
and (7) is a special case of (6).

Finally, if $s=r\cdot\lcm(2,\dots,m)$ for some $r\geq m-1$,
then using (3) and (7) we get that 
$$
\Bigl(D^i\bigl(W_s(I)\bigr)\Bigr)^s=
\bigl(W_{s-i}(I)\bigr)^s\subset W_{s(s-i)}(I)=\bigl(W_s(I)\bigr)^{s-i}.\qed
$$
\medskip

{\it Claim} \ref{w_s.props}.9. 
Let $u_1,\dots, u_m$ be variables such that $\deg(u_i)=i$.
Then any monomial $U=\prod u_i^{c_i}$ with 
$\deg(U)\geq (r+m-1)\cdot \lcm(2,\dots,m)$ 
can be written as $U=U_1\cdot U_2$, where $\deg(U_1)=r\cdot \lcm(2,\dots,m)$.
\medskip

Proof.   Set $V_i=u_i^{\lcm(2,\dots,m)/i}$, and write
$u_i^{c_i}=V_i^{b_i}\cdot W_i$ for some $b_i$ 
such that $\deg W_i<\lcm(2,\dots,m)$.

If $\sum b_i\geq r$, then choose $0\leq d_i\leq b_i$ such that
$\sum d_i=r$, and take
$U_1=\prod V_i^{b_i}$.
Otherwise, $\deg U< (r-1)\cdot \lcm(2,\dots,m)+m\cdot \lcm(2,\dots,m)$,
a contradiction.\qed
\medskip

 {\it Aside}  \ref{w_s.props}.10. Note that one can think of
(\ref{w_s.props}.9) as a statement about certain multiplication
maps 
$$
H^0(X,\o_X(a))\times H^0(X,\o_X(b))\to H^0(X,\o_X(a+b)),
$$ 
where $X$ is the weighted projective space
$\p(1,2,\dots,m)$. The above claim
is a combinatorial version  of the Castelnuovo-Mumford regularity
theorem in this case (cf.\ \cite[Sec.1.8]{laz-book}).

It seems to me that  (\ref{w_s.props}.6)
should hold for $t\geq \lcm(2,\dots,m)$ and even for many
smaller values of $t$ as well.

It is easy to see that $(m-1)\cdot\lcm(2,\dots,m)\leq m!$ for
$m\geq 6$, and one can check by hand that
(\ref{w_s.props}.6)
 holds for $t\geq m!$ for $m=1,2,3,4,5$.
Thus we conclude that $W_{m!}(I)$ is $D$-balanced.
This is not important, but the traditional choice of
the coefficient ideal corresponds to $W_{m!}(I)$.  
\medskip

The following close analog of (\ref{test.bu.seq.comp.cor})
leads to  ideal sheaves
that behave the ``same'' as a given ideal $I$, as far
as order reduction is concerned.

\begin{thm}[Tuning of ideals, I] \label{getting.equiv.ideals.1}
Let $X$ be a smooth variety, $I\subset \o_X$ an ideal sheaf and
$m=\mord I$.  Let $s\geq 1$ be an integer and $J$  any ideal sheaf satisfying
\index{Ideal sheaf!tuning}
\index{Tuning of ideal sheaves}
$$
I^{s}\subset J\subset W_{ms}(I).
$$
Then 
a smooth blow-up sequence
$$
X_r\stackrel{\pi_{r-1}}{\longrightarrow} 
X_{r-1}
\stackrel{\pi_{r-2}}{\longrightarrow}\cdots
\stackrel{\pi_{1}}{\longrightarrow} X_1
\stackrel{\pi_{0}}{\longrightarrow} X_0=X
$$
is a smooth blow-up sequence of order $\geq  m$ starting with $(X,I,m)$ iff
it is a 
smooth blow-up sequence of order $\geq ms$ starting with $(X,J,ms)$.
\end{thm}

Proof.
Assume that we get  a  smooth blow-up sequence 
 starting with $(X,I,m)$: 
$$
\begin{array}{l}
(X_r,I_r,m)\stackrel{\pi_{r-1}}{\longrightarrow} 
(X_{r-1},I_{r-1},m)
\stackrel{\pi_{r-2}}{\longrightarrow}\cdots\\
\hphantom{(X_r,I_r,m)}
\stackrel{\pi_{1}}{\longrightarrow} (X_1,I_1,m)
\stackrel{\pi_{0}}{\longrightarrow} (X_0,I_0,m)=(X,I,m).
\end{array}
$$
We prove by induction on $r$ that we also get a
  smooth blow-up sequence 
 starting with $(X,J,ms)$: 
$$
\begin{array}{l}
(X_r,J_r,ms)\stackrel{\pi_{r-1}}{\longrightarrow} 
(X_{r-1},J_{r-1},ms)
\stackrel{\pi_{r-2}}{\longrightarrow}\cdots\\
\hphantom{(X_r,J_r,ms)}
\stackrel{\pi_{1}}{\longrightarrow} (X_1,J_1,ms)
\stackrel{\pi_{0}}{\longrightarrow} (X_0,J_0,ms)=(X,J,ms).
\end{array}
$$
Assume that this holds up to step $r-1$. 
We need to show that the last blow-up
$\pi_{r-1}: X_r\to X_{r-1}$ is a  blow-up for
$(X_{r-1},J_{r-1},ms)$. That is, we need to show that
$$
\ord_ZI_{r-1}\geq m\ \Rightarrow\  \ord_ZJ_{r-1}\geq ms
\qtq{for any $Z\subset X_{r-1}$.}
$$
Let $\Pi_{r-1}:X_{r-1}\to X_0$ denote the composite.
Since $J\subset W_{ms}(I)$, we know that
$$
\begin{array}{rcl}
J_{r-1}&=&\bigl(\Pi_{r-1}\bigr)^{-1}_*(J, ms)\\
&\subset & \bigl(\Pi_{r-1}\bigr)^{-1}_*\bigl(W_{ms}(I), ms\bigr)\\
&= &\bigl(\Pi_{r-1}\bigr)^{-1}_*
\left(
\prod_j  \bigl(D^jI, m-j\bigr)^{c_j} : \sum (m-j)c_j\geq ms
\right)\\
&= &
\left(
\prod_j  \bigl((\Pi_{r-1})^{-1}_*(D^{j}I,m-j)\bigr)^{c_j} : 
\sum (m-j)c_j\geq ms
\right)\\
&\subset & 
\left(
\prod_j  \left(D^{j}(\Pi_{r-1})^{-1}_*(I,m)\right)^{c_j} : 
\sum (m-j)c_j\geq ms
\right)\quad \mbox{by (\ref{bir.trans.ders.thm})}\\
&=& 
\left(
\prod_j  \bigl(D^{j}(I_{r-1},m)\bigr)^{c_j} : 
\sum (m-j)c_j\geq ms
\right).
\end{array}
$$
If $\ord_ZI_{r-1}\geq m$, then
$\ord_ZD^{j}(I_{r-1})\geq m-j$, and so
$$ 
\ord_Z\prod_j  \bigl(D^{j}(I_{r-1},m)\bigr)^{c_j}\geq \sum (m-j)c_j\geq ms,
$$
proving one direction.

In order to prove the converse, let 
$$
\begin{array}{l}
(X_r,J_r,ms)\stackrel{\pi_{r-1}}{\longrightarrow} 
(X_{r-1},J_{r-1},ms)
\stackrel{\pi_{r-2}}{\longrightarrow}\cdots\\
\hphantom{(X_r,J_r,ms)}
\stackrel{\pi_{1}}{\longrightarrow} (X_1,J_1,ms)
\stackrel{\pi_{0}}{\longrightarrow} (X_0,J_0,ms)=(X,J,ms)
\end{array}
$$
be a  smooth blow-up sequence 
 starting with $(X,J,ms)$. Again by induction we show that
it gives a   smooth blow-up sequence 
 starting with $(X,I,m)$.
Since $I^{s}\subset J$, we know that
$$
I_{r-1}^{s}=\bigl((\Pi_{r-1})^{-1}_*I\bigr)^{s}\subset
\bigl(\Pi_{r-1}\bigr)^{-1}_*(J,ms)=J_{r-1}.
$$
Thus if $\ord_ZJ_{r-1}\geq ms$, then
$\ord_ZI_{r-1}\geq m$, and so $\pi_{r-1}: X_r\to X_{r-1}$ is also 
a  blow-up for
$(X_{r-1},I_{r-1},m)$.
\qed
\medskip

\begin{cor}[Tuning of ideals, II] \label{getting.equiv.ideals.cor}
\index{Ideal sheaf!tuning}
\index{Tuning of ideal sheaves}
Let $X$ be a smooth variety, $I\subset \o_X$ an ideal sheaf with
$m=\mord I$ and $E$ a divisor with simple normal crossings.
  Let $s=r\cdot \lcm(2,\dots,m)$ for some $r\geq m-1$.
Then $W_s(I)$ is
  MC-invariant,  $D$-balanced, and
 a smooth blow-up sequence
$$
X_r\stackrel{\pi_{r-1}}{\longrightarrow} 
X_{r-1}
\stackrel{\pi_{r-2}}{\longrightarrow}\cdots
\stackrel{\pi_{1}}{\longrightarrow} X_1
\stackrel{\pi_{0}}{\longrightarrow} X_0=X
$$
is a  blow-up sequence of order $\geq  m$ starting with $(X,I,m,E)$ iff
it is a 
 blow-up sequence of order $\geq s$ starting with $(X,W_s(I),s, E)$.
\end{cor}

Proof. Everything follows from (\ref{w_s.props})
   and (\ref{getting.equiv.ideals.1}),
except for the role played by $E$. 

Adding $E$ to $(X,I)$ (resp., to $(X,W_s(I))$) 
means that now we can use only blow-ups whose centers
are in simple normal crossing with $E$ and its total transforms.
This poses the same restriction on
smooth blow-up sequences for $(X,I,E)$ as on
smooth blow-up sequences for $(X,W_s(I),E)$.\qed

\section{Order reduction for  ideals}

In this section we prove the first main implication
 (\ref{main.red.2steps}.1) of the inductive proof.
We start with a much weaker result.
Instead of getting rid of all points of
order $m$, we  prove only  that the set of points of
order $m$ moves away from the birational transform of a given divisor
$E^j$.

\begin{lem}\label{ormI=>disj.prop}
\index{Order reduction!for ideal sheaves along a divisor}
\index{Ideal sheaf!order reduction along a divisor}
Assume that 
(\ref{ord.red.marked.thm}) holds in dimensions $< n$. Then
for every $m, j$ there is a smooth blow-up sequence functor
$\bdd_{n,m,j}$ of order $m$ that is defined on triples 
 $(X,I,E)$  with $\dim X=n$,  $\mord I\leq m$
and $E=\sum_i E^i$
such that if $\bdd_{n,m,j}(X,I,E)=$ 
$$
\begin{array}{l}
\Pi:(X_r,I_r,E_r)\stackrel{\pi_{r-1}}{\longrightarrow} 
(X_{r-1},I_{r-1},E_{r-1})
\stackrel{\pi_{r-2}}{\longrightarrow}\cdots\\
\hphantom{\Pi:(X_r,I_r,E_r)}
\stackrel{\pi_{1}}{\longrightarrow} (X_1,I_1,E_1)
\stackrel{\pi_{0}}{\longrightarrow} (X_0,I_0,E_0)=(X,I,E),
\end{array}
$$
then
\begin{enumerate} 
\item  $\cosupp(I_r,m)\cap \Pi^{-1}_*E^j=\emptyset$, and
\item  $\bdd_{n,m,j}$  commutes with smooth morphisms
(\ref{funct.pack.proc.say}.1) 
and also with change of fields (\ref{funct.pack.proc.say}.2).
\end{enumerate}
Assume in addition that there is an ideal $J\subset \o_{E^j}$
such that  $J$ is nonzero on every
irreducible component of $E^j$ and $\tau_*\bigl(\o_{E^j}/J\bigr)=\o_X/I$, 
where $\tau:E^j\into X$ is the natural injection. Then
\begin{enumerate} \setcounter{enumi}{2}
\item $\bdd_{n,1,j}(X,I,E):=\tau_*\bmord_{n-1,1}(E^j, J,1, (E-E^j)|_{E^j})$.
\end{enumerate}
\end{lem}

Proof. By (\ref{getting.equiv.ideals.cor}), 
$W_{m!}(I)$  is  $D$-balanced and
order reduction for $(X,I,E)$ is equivalent to
 order reduction for $(X,W_{m!}(I),E)$.
Thus from now on we assume that $I$ is
 $D$-balanced.

Let $Z_{-1}$ be the union of those irreducible
components 
$E^{jk}\subset E^j$ that are contained in $\cosupp(I,m)$.
Let $\pi_{-1}:X_0\to X$ be the blow-up of $Z_{-1}$.
The blow-up is an isomorphism, but the
order of $I$ along  $E^{jk}$ is reduced 
by $m$ and we get a new
 ideal sheaf $I_0$. Since $\mord_{E^{jk}}  I\leq m$
to start with,  
$\mord_{E^{jk}}  I_0=\mord_{E^{jk}}  I-m\leq 0$. Thus
 $\cosupp (I_0,m)$ does not contain any
irreducible
component of $E^j$.

Next, set 
$S:=E^j$ with injection $\tau:S\into X$, $E_S:=(E-E^j)|_S$
  and consider the triple  $(S, I_0|_S, E_S)$.
By the going-up theorem (\ref{D-bal.goup.thm}),
 every blow-up  sequence of order $\geq m$ starting with 
 $(S, I_0|_S,m, (E-E^j)|_S)$ corresponds to a
 blow-up  sequence of order $m$ starting with 
 $(X_0,I_0, E-E^j)$.
 Since $S= E^j $, every blow-up center is
a smooth subvariety of the birational transform of $E^j$; thus we
in fact get  a blow-up  sequence of order $m$ starting with 
 $(X_0,I_0, E)$. 
Set 
$$
\bdd_{n,m,j}(X,I,E):=\tau_*\bmord_{n-1,m}(S,  I_0|_S, m, E_S)
\stackrel{\pi_{-1}}{\longrightarrow} X.
$$
That is, we take $\bmord_{n-1,m}(S,  I_0|_S, m, (E-E^j)|_S)$,
 push it forward (\ref{bu.seq.functs}.3)  
 and compose the resulting blow-up sequence
on the right with our first blow-up $\pi_{-1}$.
(This is the reason for the subscript $-1$.)
By (\ref{D-bal.maxcont.cor}) we obtain
$$
\Pi_r:X_{r}\to X\qtq{with} 
 I_{r}:=\bigl(\Pi_{r}\bigr)^{-1}_*I,\ 
E_{r}:=\bigl(\Pi_{r}\bigr)^{-1}_{\rm tot}(E)
$$
such  that $\bigl(\Pi_{r}\bigr)^{-1}_*(E^j)$
is disjoint from $\cosupp (I_{r},m)$.

The functoriality properties of
$\bdd_{n,m,j}(X,I,E)$ follow from the corresponding
functoriality properties of
$\bmord_{n-1,m}(S,  I_0|_S, E_S)$.
All the steps are obvious, but for the first time, let us go
through the details.

Let $h:Y\to X$ be a smooth surjection. Set $E^j_Y:=h^{-1}(E^j)$.
Then $h|_{E^j_Y}:E^j_Y\to E^j$ is also a smooth surjection and
we get the same result whether we first pull back by $h$ and then
restrict to $E^j_Y$ or we first restrict to $E^j$ and
then pull back by  $h|_{E^j_Y}$. That is, 
$$
\bigl(h|_{E^j_Y}\bigr)^*\bigl(E^j, I_0|_{E^j},m, (E-E^j)|_{E^j}\bigr)
=
\bigl(E^j_Y, (h^*I)_0|_{E^j_Y}, h^{-1}(E-E^j)|_{E^j_Y}\bigr).
$$
Therefore, 
$$
\begin{array}{l}
\bmord_{n-1,m}\bigl(E^j_Y, (h^*I)_0|_{E^j_Y}, h^{-1}(E-E^j)|_{E^j_Y}\bigr)
\\*[1ex]
\qquad\qquad\qquad =\bigl(h|_{E^j_Y}\bigr)^*\bmord_{n-1,m}
\bigl(E^j, I_0|_{E^j},m, (E-E^j)|_{E^j}\bigr),
\end{array}
$$
and hence
$$
h^*\bdd_{n,m,j}(X,I,E)=\bdd_{n,m,j}(Y, h^*I, h^{-1}(E)\bigr).
$$

If $h:Y\to X$ is any smooth morphism, we see similarly that the same
blow-ups end up with empty centers.

The functoriality property (\ref{funct.pack.proc.say}.2)
holds since  change of the base field commutes with restrictions.

Assume finally that  $\o_X/I=\tau_*(\o_{E^j}/J)$.
Every local equation of $E^j$ is an order 1 element
in $I$. Thus $m=\mord I=1$ and so $W_{m!}(I)=I$.
If $J$ is nonzero on every
irreducible component of $E^j$ then 
$Z_{-1}=\emptyset$ and so $I_0=I$, proving (3).
\qed
\medskip

The  main theorem of this section 
is the following.

\begin{thm}\label{ormI=>orI.prop}
\index{Order reduction!for ideal sheaves}
\index{Ideal sheaf!order reduction for}
Assume that 
(\ref{ord.red.marked.thm}) holds in dimensions $< n$. Then
for every $m$ there is a smooth blow-up sequence functor
$\bord_{n,m}$ of order $m$ that is defined on triples 
 $(X,I,E)$  with $\dim X=n$ and  $\mord I\leq m$
such that if $\bord_{n,m}(X,I,E)=$ 
$$
\begin{array}{l}
\Pi:(X_r,I_r,E_r)\stackrel{\pi_{r-1}}{\longrightarrow} 
(X_{r-1},I_{r-1},E_{r-1})
\stackrel{\pi_{r-2}}{\longrightarrow}\cdots\\
\hphantom{\Pi:(X_r,I_r,E_r)}
\stackrel{\pi_{1}}{\longrightarrow} (X_1,I_1,E_1)
\stackrel{\pi_{0}}{\longrightarrow} (X_0,I_0,E_0)=(X,I,E),
\end{array}
$$
then
\begin{enumerate} 
\item  $\mord I_r<m$, and
\item  $\bord_{n,m}$   commutes with smooth morphisms
(\ref{funct.pack.proc.say}.1) 
and also with change of fields (\ref{funct.pack.proc.say}.2).
\end{enumerate}
Assume in addition that there is a smooth hypersurface
$\tau:Y\into X$ and an ideal sheaf $J\subset \o_Y$
such that  $J$ is nonzero on every
irreducible component of $Y$ and $\tau_*(\o_Y/J)=\o_X/I$.
Then $\mord I=1$ and 
\begin{enumerate} \setcounter{enumi}{2}
\item $\bord_{n,1}(X,I,\emptyset)=
\tau_*\bmord_{n-1,1}(Y, J,1,\emptyset )$.
\end{enumerate}
\end{thm}

The proof is done in three steps.

{\it Step 1} (Tuning $I$).
\index{Ideal sheaf!tuning}%
\index{Tuning of ideal sheaves}% 
By (\ref{getting.equiv.ideals.cor}), there is an ideal
$W(I)=W_s(I)$ for suitable $s$,
 which is  $D$-balanced, MC-invariant and
order reduction for $(X,I,E)$ is equivalent to
 order reduction for $(X,W(I),E)$.
(Let us take $s=m!$ to avoid further choices.)
Thus from now on we assume that $I$ is
 $D$-balanced and MC-invariant.

{\it Step 2} (Maximal contact case).
\index{Maximal contact!and order reduction}%
 Here we assume that there is a smooth hypersurface of
maximal contact $H\subset X$. This is always satisfied in a
suitable open neighborhood of any point by
(\ref{max.contact.thm.I}.2), but it may hold globally as well.
This condition is also preserved under disjoint unions.

Under a smooth blow-up of order $m$, the birational
transform of a smooth hypersurface of
maximal contact is again a smooth hypersurface of
maximal contact; thus we stay in the
maximal contact case.

We intend to restrict everything to $H$, but
we run into the problem  that $E|_H$ need not be
a simple normal crossing divisor. We take care of this problem first.

{\it Step 2.1.}  If $E=\sum_{i=1}^s E^i$, we apply
(\ref{ormI=>disj.prop}) to each $E^i$.
At the end we get a blow-up sequence
$\Pi:X_r\to X$ such that
$\cosupp(I_r,m)$ is disjoint from $\Pi^{-1}_*E$.

 Note that  the new exceptional divisors
obtained in the process (and added to $E$) have simple normal
crossings with the birational transforms of $H$, so
$H_r+E_r$ is a simple normal
crossing divisor.
(I have used the ordering of the
index set of $E$. This is avoided  traditionally  by
restricting $(X,I,E)$ successively to the multiplicity
$n-j$ locus of  $E$, starting with the case $j=0$.
The use of the ordering cannot  be avoided
in (\ref{orI=>ormI.say}.3), so
there is not  much reason  to go around it here.)

{\it Step 2.2.}  Once $H+E$ is a simple normal crossing divisor,
 we restrict 
everything to the birational transform of $H$, and 
we obtain order reduction using dimension induction
and (\ref{ormI=>disj.prop}).

{\it Step 3} (Global case). There may not be 
 a global smooth hypersurface of
maximal contact $H\subset X$, but we can cover $X$ with
 open subsets  $X^{(j)}\subset X$ such that
on each $X^{(j)}$  there is a smooth hypersurface of
maximal contact $H^{(j)}\subset X^{(j)}$.
Thus the disjoint union
$$
H^*:=\textstyle{\coprod_j} H^{(j)}\subset 
\textstyle{\coprod_j} X^{(j)}=:X^*
$$
is a smooth hypersurface of
maximal contact.
Let $g:X^*\to X$ be the coproduct of the injections
$X^{(j)}\into X$.

By the previous step
$\bord_{n,m}$ is defined on
$(X^*,g^*I,g^{-1}E)$.
Then we argue as in (\ref{aff.to.glob.prop}) to prove that
$\bord_{n,m}(X^*,g^*I,g^{-1}E)$ descends to give
$\bord_{n,m}(X,I,E)$.
\medskip

Only Steps 2 and 3 need amplification.

\begin{say}[Step 2, Maximal contact  case]\label{local.oredred.say}
\index{Maximal contact!and order reduction}%
\index{Order reduction!maximal contact case}%
We start with a triple $(X,I,E)$, where
$I$ is $D$-balanced and MC-invariant, and assume that
there is a smooth hypersurface of maximal contact
$H\subset X$. Set $m=\mord I$.

{\it Warning.} As we blow up, we get 
birational transforms of $I$ which may be neither 
$D$-balanced nor MC-invariant.
We do not attempt to ``fix'' this problem, since the
 relevant consequences of these properties
(\ref{D-bal.goup.thm}) and (\ref{maxcont.formal.automs.thm})
 are established for any sequence of blow-ups
of order $m$. This also means that we should not  pick
new hypersurfaces of maximal contact after a blow-up but rather
stick with the birational transforms of the old ones.
\medskip

{\it Step 2.1} (Making $\cosupp (I_r,m)$
 and $\Pi^{-1}_* E$ disjoint).
To fix notation, write $E=\sum_{i=1}^s E^i$,
and set $(X_0,I_0,E_0):=(X,I,E)$ and $H_0:=H$.
The triple $(X_0,I_0,E_0)$ satisfies the assumptions of Step 2.1.1.
\medskip

{\it Step 2.1.j.}
Assume that we have already constructed a smooth blow-up sequence of order
$m$ starting with  $(X_0,I_0,E_0)$ whose end result is
$$
\begin{array}{l}
\Pi_{r(j-1)}:X_{r(j-1)}\to \cdots \to  X_0,  \qtq{where}\\
I_{r(j-1)}:=\bigl(\Pi_{r(j-1)}\bigr)^{-1}_*I \qtq{and} 
\ E_{r(j-1)}:=\bigl(\Pi_{r(j-1)}\bigr)^{-1}_{\rm tot}(E),
\end{array}
$$
such that
$$
\bigl(\Pi_{r(j-1)}\bigr)^{-1}_*(E^i)\cap \cosupp (I_{r(j-1)},m)
=\emptyset \qtq{for $i< j$.}
$$

Apply (\ref{ormI=>disj.prop})
 to $\bigl(X_{r(j-1)}, I_{r(j-1)}, E_{r(j-1)}\bigr)$
and the divisor $E_{r(j-1)}^j$
to obtain
$$
\Pi_{r(j)}:X_{r(j)}\to X_{r(j-1)}\cdots \to  X_0
$$
such that 
$$
\bigl(\Pi_{r(j)}\bigr)^{-1}_*(E^i)\cap \cosupp (I_{r(j)},m)
=\emptyset \qtq{for $i\leq  j$.}
$$

Note that the center of every blow-up is contained in
every hypersurface of maximal contact. Thus
 $H_{r(j)}:=\bigl(\Pi_{r(j)}\bigr)^{-1}_*H$ is a
smooth hypersurface of maximal contact,
and every new  divisor
in $\bigl(\Pi_{r(j)}\bigr)^{-1}_{\rm tot}E$ is transversal to $H_{r(j)}$.
If  $E=\sum_{i=1}^s E^i$, 
then after Step 2.1.s, we have achieved that 
\begin{enumerate}
\item[$\bullet$] $\bigl(\Pi_{r(s)}\bigr)^{-1}_*E$
is disjoint from $\cosupp (I_{r(s)},m)$, and
\item[$\bullet$]
for any hypersurface of maximal contact $H\subset X$, 
the divisor
$H_{r(s)}+E_{r(s)}$
has simple normal crossing along
$\cosupp (I_{r(s)},m)$.
\end{enumerate}
Note that we perform all these steps even if
$H+E$ is a simple normal crossing divisor to start with, though in
this case they do not seem  to be necessary. 
We would, however,  run into problems with the 
compatibility of the numbering
in the blow-up sequences otherwise.

{\it Step 2.2} (Restricting to $H$).
After dropping the  subscript $r(s)$
we have
a triple $(X,I,E)$ and a smooth hypersurface of maximal contact
$H\subset X$ such that $H+E$ is also a simple normal crossing divisor.
We can again replace $I$ by $W(I)$ and thus assume that
$I$ is MC-invariant.
Note that we do not pick a new hypersurface of maximal contact,
but  use only
the birational transforms $H_{r(s)}$ of the old
 hypersurfaces of maximal contact.

Declare $E^0:=H$ to be the first divisor in $H+E$
and apply   (\ref{ormI=>disj.prop}) to $(X,I,H+E)$
with $j=0$.
This gives a sequence of blow-ups
$\Pi:X_r\to X$ such that $\cosupp\Pi^{-1}_*(I,m)$ is disjoint from
$\Pi^{-1}_*H$. However, $H$ is a 
 smooth hypersurface of maximal contact, and hence, by definition,
$\cosupp\Pi^{-1}_*(I,m)\subset \Pi^{-1}_*H$.
Thus $\cosupp\Pi^{-1}_*(I,m)=\emptyset$, as we wanted.
\medskip

{\it Step 2.3} (Functoriality). Assuming  functoriality
 in dimension $< n$, we have
functoriality in Step 2.1 by the corresponding functoriality
in (\ref{ormI=>disj.prop}).

In Step 2.2 we rely on
the choice of a hypersurface of maximal contact $H$,
which is not unique.  Let $H,H'$ be two hypersurfaces of maximal contact
such that $H+E$ and $H'+E$ are both 
simple normal crossing divisors.
We can use either of the  two blow-up sequences
 $\bdd_{n,m,0}(X,I,H+E)$  and $\bdd_{n,m,0}(X,I,H'+E)$
to construct  $\bord_{n,m}(X,I,E)$.

Here we need that  
$I$ is MC-invariant. By (\ref{maxcont.formal.automs.thm}) this implies that 
$(X,I,H+E)$ and $(X,I,H'+E)$ are  \'etale equivalent.
Thus the two
blow-up sequences
 $\bdd_{n,m,0}(X,I,H+E)$  and $\bdd_{n,m,0}(X,I,H'+E)$ are 
also \'etale equivalent. By (\ref{et.eq=eq})
this implies that these 
blow-up sequences are identical.

As we noted in (\ref{funct.pack.proc.say}), the
functoriality package is local, so
we do not have to consider it separately
in the next  step.
\medskip

{\it Step 2.4} (Closed embeddings) 
Let  $\tau:Y\into X$ be a smooth hypersurface
 and  $J\subset \o_Y$ an ideal sheaf
such that  $J$ is nonzero on every
irreducible component of $Y$ and $\tau_*(\o_Y/J)=\o_X/I$.
Then $I$ contains the local equations of $Y$, and
so it has order 1. In particular, $I=W(I)$.
If $E=\emptyset$ then Step 2.1 does nothing, and in
Step 2.2 we can choose $H=Y$. 
Thus (\ref{ormI=>orI.prop}.3) follows from
(\ref{ormI=>disj.prop}.3).
\end{say}

As in (\ref{aff.to.glob.prop}), going from the local to the global case is
essentially automatic.
For ease of reference, let us axiomatize
the process.

\begin{thm}[Globalization of blow-up sequences]\label{glob.bus.prop}
\index{Order reduction!global case}%
Assume that we have the following:
\begin{enumerate}
\item a class of smooth morphisms ${\mathcal M}$
that is closed under fiber products and coproducts
(for instance, ${\mathcal M}$ could be all smooth morphisms,
all \'etale morphisms or all open immersions);
\item two classes of triples
${\mathcal {GT}}$ (global triples) and ${\mathcal {LT}}$
(local triples) such that 
\begin{enumerate}
\item[(i)] for every
$(X,I,E)\in {\mathcal {GT}}$ and every $x\in X$
there is an ${\mathcal M}$-morphism
$g_x:(x'\in U_x)\to (x\in X)$ such that
$(U_x, g^*I,g^{-1}E)$ is in ${\mathcal {LT}}$, and
\item[(ii)] ${\mathcal {LT}}$ is closed under disjoint unions;
\end{enumerate}
\item a blow-up sequence functor $\bb$ defined on
${\mathcal {LT}}$ that commutes with 
surjections in ${\mathcal M}$.
\end{enumerate}
Then $\bb$ has a unique extension to
a blow-up sequence functor $\overline{\bb}$,
which is defined on ${\mathcal {GT}}$ and
which commutes with  surjections in ${\mathcal M}$.
\end{thm}

Proof. For any $(X,I,E)\in {\mathcal {GT}}$ choose 
${\mathcal M}$-morphisms $g_{x_i}:U_{x_i}\to X$
such that the images cover $X$.

Let $X':=\coprod_i U_{x_i}$ be the disjoint union
and $g:X'\to X$ the induced ${\mathcal M}$-morphism.
By assumption $(X', g^*I,g^{-1}E)\in {\mathcal {LT}}$.

Set $X'':=X'\times_XX'$. By assumption
the  two coordinate projections 
$\tau_1,\tau_2:X''\to X'$  are in ${\mathcal M}$
and are surjective.

The blow-up sequence $\bb$ for $X'$ starts with blowing up
$Z'_0\subset X'$,
and the blow-up sequence $\bb$ for $X''$ starts with blowing up
$Z''_0\subset X''$.
Since $\bb$ commutes with the $\tau_i$, we conclude
that
$$
\tau_1^*(Z'_0)=Z''_0=\tau_2^*(Z'_0).
\eqno{(\ref{glob.bus.prop}.4)}
$$
If ${\mathcal M}= \{\mbox{open immersions}\}$,
we have proved in (\ref{aff.to.glob.prop})
that the subschemes $Z'_0\cap U_{x_i}\subset X$ 
 glue together
to a subscheme $Z_0\subset X$.
This is the only case  we need for the proof of (\ref{ormI=>orI.prop}).

The conclusion still holds for any  ${\mathcal M}$, but
we have to use the theory of faithfully flat descent; see
\cite{gro-ffd} or \cite[Ch.VII]{murre}.

This way we obtain $X_1:=B_{Z_0}X$
such that $X'_1=X'\times_XX_1$. We can repeat the above argument
to obtain the center $Z_1\subset X_1$ and eventually get the
whole blow-up sequence for $(X,I,E)$.\qed
\medskip

The following example, communicated to me by
Bierstone and Milman, shows that
while  principalization proceeds by
smooth blow-ups, the resolution of singularities also
involves blowing up singular centers.

\begin{exmp}\label{nsbu.exmp} 
Consider the subvariety  $X\subset \a^4$ defined by the ideal
	$I=( x^3-y^2, x^4 + xz^2 -w^3)$.
Let us see how the principalization proceeds.

Note that $\ord I=2$ and $H=(y=0)$ is a hypersurface of maximal
contact. $I|_H=(x^3, xz^2 -w^3)$ has order 3
and $MC(I|_H)=(x,z,w)$. Thus the first step
is to blow up the origin in $\a^4$.

Consider the chart $x_1=x, y_1=y/x, z_1=z/x, w_1=w/x$.
The birational transform of $I$ is
$I_1=(x_1 - y_1^2, x_1(x_1 + z_1^2 -w_1^3))$
and $E_1=(x_1=0)$. The order  has dropped to 1,
so we continue with $(I_1,1,E_1)$.

Since $\cosupp(I_1,1)$  is not disjoint from $E_1$,
  we proceed as in (\ref{local.oredred.say}.1).
The restriction is $I_1|_{E_1}=(x_1,y_1^2)$, and thus
 next we have to blow up
$(x_1=y_1=0)$.

On the other hand, the birational transform of
$X$ is 
$$
X_1=(x_1 - y_1^2= x_1 + z_1^2 -w_1^3=0).
$$
Its intersection with $(x_1=y_1=0)$
is the cuspidal  curve 
$(x_1=y_1=z_1^2 -w_1^3=0)$. Thus the resolution of $X$
first blows up the origin and then the new exceptional curve,
which is singular.
\end{exmp}

\section{Order reduction for marked ideals}\label{or.red.m.sec}

In this section we prove the second main implication
 (\ref{main.red.2steps}.2) of the inductive proof. That is,
we prove the following.

\begin{thm}\label{orI=>ormI.prop}
\index{Order reduction!for marked ideal sheaves}
\index{Ideal sheaf!marked, order reduction}
Assume that (\ref{ord.red.I.thm})
 holds in dimensions $\leq  n$. Then
for every $m$, there is a smooth blow-up sequence functor
$\bmord_{n,m}$ defined on triples 
 $(X,I,m,E)$  with $\dim X=n$ 
such that if $\bmord_{n,m}(X,I,m,E)=$ 
$$
\begin{array}{l}
\Pi:(X_r,I_r,m,E_r)\stackrel{\pi_{r-1}}{\longrightarrow} 
(X_{r-1},I_{r-1},m,E_{r-1})
\stackrel{\pi_{r-2}}{\longrightarrow}\cdots\\
\hphantom{\Pi:(X_r,I_r,m,E_r)}
\stackrel{\pi_{1}}{\longrightarrow} (X_1,I_1,m,E_1)
\stackrel{\pi_{0}}{\longrightarrow} (X_0,I_0,m,E_0)=(X,I,m,E),
\end{array}
$$
then
\begin{enumerate} 
\item  $\mord I_r<m$, 
\item  $\bmord_{n,m}$ commutes with smooth morphisms
(\ref{funct.pack.proc.say}.1) 
and with change of fields (\ref{funct.pack.proc.say}.2), and
\item if $m=\mord I$ then 
$\bmord_{n,m}(X,I,m,\emptyset)=\bord_{n,m}(X,I,\emptyset)$.
\end{enumerate}
\end{thm}

Before proving (\ref{orI=>ormI.prop}), we show that
it implies the two claims in (\ref{pr.to.res.say}).

\begin{say}[Proof of (\ref{pr.to.res.say}.1--2)]
First, (\ref{pr.to.res.say}.1) is the same as 
(\ref{orI=>ormI.prop}.3).

 The claimed identity in (\ref{pr.to.res.say}.2)
is a local question on $X$, thus we may assume
that there is a chain of smooth subvarieties
$Y=Y_0\subset Y_1\subset \cdots \subset Y_c=X$
such that each is a hypersurface in the next one.
Thus it is enough to prove the case when $Y$ is a
hypersurface in $X$.

Every local equation of $Y$ is in $I$, thus
$\mord I=1$. Therefore, 

\noindent $\bmord_{\dim X,1}(X,I,1,\emptyset)=\bord_{\dim X,1}(X,I,\emptyset)$
by (\ref{orI=>ormI.prop}.3)
and  (\ref{ormI=>orI.prop}.3) gives that 

\noindent $\bord_{\dim X,1}(X,I,\emptyset)=
\tau_*\bmord_{\dim Y,1}(Y,J,1,\emptyset)$.
Putting the two together gives  (\ref{pr.to.res.say}.2).\qed
\end{say}

\begin{say}[Plan of the proof of (\ref{orI=>ormI.prop})]

{\it Step 1.} We start with the unmarked triple $(X,I,E)$,
 and using (\ref{ord.red.I.thm})  in dimension $n$,
we reduce its order below $m$. That is, we get a
 composite of  smooth blow-ups
$\Pi_1:X^1\to X$ such that
$(\Pi_1)^{-1}_*I$ has order $<m$. The problem is that $(\Pi_1)^{-1}_*I$
 differs from $(\Pi_1)^{-1}_*(I,m)$
along the exceptional divisors of $\Pi_1$, and the latter
can have very high order.
We decide not to worry about it for now.

{\it Step 2.}  Continuing
 with $(X^1, (\Pi_1)^{-1}_*(I,m), (\Pi_1)^{-1}_{\rm tot}E)$,  we
blow up  subvarieties where 
 the birational transform of 
$(I,m)$ has order $\geq m$, and
   the birational transform of $I$ has order $\geq 1$.

 Eventually we get $\Pi_2:X^2\to X$
such that  $\cosupp(\Pi_2)^{-1}_*I$ is disjoint from the
locus where $(\Pi_2)^{-1}_*(I,m)$  has order $\geq m$.
 We can now completely ignore
$(\Pi_2)^{-1}_*I$.
Since $(\Pi_2)^{-1}_*I$ and $(\Pi_2)^{-1}_*(I,m)$
agree up to tensoring with the ideal sheaf of a divisor whose support is
in $E_r$, we can assume from now on that
 $(\Pi_2)^{-1}_*(I,m)$ is the ideal sheaf of a divisor with simple normal 
crossing.

{\it Step 3.} 
Order reduction for  the marked ideal sheaf of a divisor with simple normal 
crossing is rather easy.
\medskip

Instead of strictly following this plan, we divide the ideal into
a ``simple normal crossing part'' and the ``rest'' using all of $E$, 
instead of 
exceptional divisors only.
This is  solely a notational convenience.
\end{say}

\begin{defn-lem}\label{mon.nonmon.defn}
 Given $(X,I,E)$, we can write $I$ uniquely as
$I=M(I)\cdot N(I)$, where $M(I)=\o_X(-\sum c_iE^i)$ for some $c_i$ and
$\cosupp N(I)$ does not contain any of the $E^i$.
$M(I)$ is called the {\it monomial part} of $I$ 
\index{Monomial!part of an ideal sheaf}%
\index{Ideal sheaf!monomial part of}%
\index{Mi@$M(I)$,  monomial part of an ideal sheaf}%
and
$N(I)$  the {\it nonmonomial part} of $I$.
\index{Non-monomial part}%
\index{Ideal sheaf!nonmonomial part of}%
\index{Ni@$N(I)$,  nonmonomial part}%

Since the $E^i$ are not assumed irreducible,
 $\cosupp N(I)$  may contain  irreducible components
 of some of the $E^i$.
\end{defn-lem}

\begin{say}[Proof of (\ref{orI=>ormI.prop})]\label{orI=>ormI.say}

We write $I=M(I)\cdot N(I)$ and try to deal with the two parts
separately.
\medskip

{\it Step 1}
(Reduction to $\ord N(I)< m$).
 If $\ord N(I)\geq m$,
we can apply order reduction (\ref{ord.red.I.thm}) to $N(I)$,
until its order drops below $m$. This happens at some
$\Pi_1:X^1\to X$. Note that the two birational transforms
$$
(\Pi_1)^{-1}_*N(I)\qtq{and} (\Pi_1)^{-1}_*(I,m)
$$
differ only by tensoring with an ideal sheaf of  exceptional divisors
of $\Pi_1$, thus only in their monomial part.
Therefore,
$$
N\bigl((\Pi_1)^{-1}_*(I,m)\bigr)=(\Pi_1)^{-1}_*N(I),
$$
and so 
 we have reduced to the case where
the maximal order of the nonmonomial part is  $< m$.
\medskip

To  simplify notation, instead of
 $(X^1, (\Pi_1)^{-1}_*(I,m), (\Pi_1)^{-1}_{\rm tot}(E))$,
 write
$(X,I,m,E)$. From now on
we may assume that $\mord N(I)<m$.
\medskip

{\it Step 2} (Reduction to 
$\cosupp (I,m)\cap \cosupp N(I)=\emptyset$).
Our aim is  to continue with order reduction further
and get rid of $N(I)$ completely. The problem is that we are allowed
to blow up only subvarieties along which 
 $(I,m)$ has order at least $m$.
Thus we can blow up $Z\subset X$ with $\ord_ZN(I)<m$
only if $\ord_Z I\geq m$. 
We will be able to  guarantee this interplay by a simple trick.

Let $s$ be the maximum order of $N(I)$ along $\cosupp (I,m)$.
We reduce this order step-by-step, eventually ending up
with $s=0$, which is the same as 
$\cosupp (I,m)\cap \cosupp N(I)=\emptyset$.

It would not have been difficult to develop order reduction theory
for  several marked ideals and to apply it to the pair of marked ideals  
$(N(I),s)$ and $ (I,m)$, but the following simple observation
  reduces the general case to a single ideal:
$$
\mbox{ $\ord_ZJ_1\geq s$ and $\ord_ZJ_2\geq m$}
\Leftrightarrow
\ord_Z(J_1^{m}+J_2^s)\geq ms.
$$
Thus 
we apply order reduction to the ideal $N(I)^{m}+I^s$,
which has order  $\geq ms$. Every smooth blow-up sequence of order
$ms$ starting with $N(I)^{m}+I^s$ is also a
smooth blow-up sequence of order
$s$ starting with $N(I)$ and a smooth blow-up sequence of order
$m$ starting with $I$.  Thus we stop after $r=r(m,s)$ steps when
we have achieved
$\cosupp (I_r,m)\cap \cosupp (N(I_r),s)=\emptyset$.
We can continue with $s-1$ and so on.

Eventually we  achieve a situation where
(after dropping the  subscript) the cosupports of 
$N(I)$ and of $(I,m)$ are disjoint.
Since the center of any further  blow-up is contained in
$\cosupp (I,m)$, we can replace $X$ by $X\setminus \cosupp N(I)$
and thus assume that
$I=M(I)$.
The final step is now to deal with monomial ideals.
 
\medskip
{\it Step 3} (Order reduction for $M(I)$).
\index{Order reduction!monomial case}
Let $X$ be a smooth variety, $\cup_{j\in J} E^j$
a simple normal crossing divisor with ordered index set $J$
and $a_j$ natural numbers giving the monomial ideal
$I:=\o_X(-\sum a_jE^j)$.

The usual method  would be to
look for the highest multiplicity locus and blow it up.
This, however, does not work, not even for surfaces;
see (\ref{last.surfexmp}).

The only thing that saves us at this point is that the
divisors $E^i$ come with an ordered index set.
This allows us to specify in which order to blow up.
There are many possible choices. As far as I can tell,
there is no natural or best variant.

{\it Step 3.1.} Find the smallest $j$ such that
$a_j\geq m$ is maximal. If there is no such $j$, go to the next step.
Otherwise,  blow  up $E^j$.
Repeating this, we eventually get to the point where
 $a_j<m$ for every $j$.

{\it Step 3.2.} Find the lexicographically smallest $(j_1<j_2)$ such that
$E^{j_1}\cap E^{j_2}\neq \emptyset$ and $a_{j_1}+a_{j_2}\geq m$ is maximal.
If there is no such $(j_1<j_2)$, go to the next step.
Otherwise, 
blow  up  $E^{j_1}\cap E^{j_2}$. We get a new divisor, and
put it last as $E^{j_{\ell}}$. Its
 coefficient is $ a_{j_{\ell}}= a_{j_1}+a_{j_2}-m<m$. The new pairwise
intersections are
$E^i\cap E^{j_{\ell}}$ for certain values of $i$.
Note that 
$$
a_i+a_{j_{\ell}}=a_i+a_{j_1}+a_{j_2}-m<a_{j_1}+a_{j_2},
$$
since $a_i<m$ for every $i$ by Step 3.1.

At each repetition, the pair $(m_2(E),n_2(E))$ decreases lexicographically
where
$$
\begin{array}{lll}
m_2(E)&:=&\max\{a_{j_1}+a_{j_2}: E^{j_1}\cap E^{j_2}\neq \emptyset\},\\
n_2(E)&:=& \mbox{number of $(j_1<j_2)$ achieving the maximum.}
\end{array}
$$
Eventually we reach the stage where
$a_{j_1}+a_{j_2}< m$ whenever $E^{j_1}\cap E^{j_2}\neq \emptyset$.

 {\it Step 3.r.} Assume that  for every $s<r$ we already have the property
$$
a_{j_1}+\cdots +a_{j_{s}}< m\quad \mbox{if}\quad 
j_1<\cdots<j_{s}  \ \mbox{and}\ 
E^{j_1}\cap\cdots \cap  E^{j_{s}}\neq \emptyset.
\eqno{(*_{s})}
$$
   Find the  lexicographically smallest 
$(j_1<\cdots <j_r)$ such that
$E^{j_1}\cap\cdots \cap E^{j_{r}}\neq \emptyset$ and
$a_{j_1}+\cdots +a_{j_{r}}\geq m$ is maximal. 
If there is no such $(j_1<\cdots <j_r)$, go to the next step.
Otherwise, 
blow up $E^{j_1}\cap\cdots \cap E^{j_{r}}$, and put the new divisor 
$E^{j_{\ell}}$ last
with coefficient $a_{j_1}+\cdots +a_{j_{r}}-m$.
As before,  the new $r$-fold
intersections are of the form
$E^{i_1}\cap\cdots \cap E^{i_{r-1}}\cap E^{j_{\ell}}$,
where $E^{i_1}\cap\cdots \cap E^{i_{r-1}}\neq \emptyset$.
Moreover, 
$$
a_{i_1}+\cdots +a_{i_{r-1}}+a_{j_{\ell}}=
\bigl(a_{i_1}+\cdots +a_{i_{r-1}}-m\bigr)+a_{j_1}+\cdots +a_{j_{r}},
$$
which is less than
$a_{j_1}+\cdots +a_{j_{r}}$ since
$a_{i_1}+\cdots +a_{i_{r-1}}<m$ by Step 3.$r-1$.
 Thus the pair $(m_r(E),n_r(E))$ decreases lexicographically,
 where
$$
\begin{array}{lll}
m_r(E)&:=&\max\{a_{j_1}+\cdots +a_{j_{r}}:
 E^{j_1}\cap \cdots \cap E^{j_{r}}\neq \emptyset\},\\
n_r(E)&:=& \mbox{number of $(j_1<\cdots<j_r)$ achieving the maximum.}
\end{array}
$$
Eventually we  reach the stage where
the property ($*_r$) also holds.
We can now move to the next step.

 At the end of Step 3.n we are done, where $n=\dim X$.

The functoriality conditions are just as obvious as before.

The process greatly simplifies if  $m=\mord I$ 
and $E=\emptyset$. 
First, if $E=\emptyset$ then $N(I)=I$.
Thus in Step 1 we apply $\bord_{n,m}(X,I,\emptyset)$.
Each blow-up has order $m$, and thus the
birational transforms of $I$ agree with the 
birational transforms of $(I,m)$.
At the end of Step 1,
$(\Pi_1)^{-1}_*I = (\Pi_1)^{-1}_*(I,m)$.
Thus
$$
\cosupp\bigl((\Pi_1)^{-1}_*(I,m)\bigr)=\emptyset
\qtq{and}
M\bigl((\Pi_1)^{-1}_*(I,m)\bigr)=\o_{X^1}.
$$
Steps 2 and 3 do nothing, and so
$\bmord_{n,m}(X,I,m,\emptyset)=\bord_{n,m}(X,I,\emptyset)$.
\qed
\end{say}

\begin{exmp}\label{last.surfexmp}
Let $S$ be a smooth surface,  $E^1,E^2$ two 
 2 curves  intersecting at a point $p=E^1\cap E^2$ and $a_1=a_2=m+1$.
Let $\pi:S_3\to S$ be the blow-up of $p$ with exceptional
curve $E^3$. Then 
$$
\pi^{-1}_*\bigl(\o_S(-(m+1)(E^1+E^2)),m\bigr)=
\bigl(\o_{S_3}(-(m+1)(E^1+E^2)-(m+2)E^3),m\bigr).
$$
Next we blow up the intersection point  $E^2\cap E^3$
and so on.  After $r-2$ steps we get a birational transform
$$
\bigl(\o_{S_r}(-\textstyle{\sum_{i=1}^r}(m+p_i)E^i),m\bigr),
$$
where $p_i$ is the $i$th Fibonacci number.
Thus  we get  higher and higher
multiplicity  ideals.
\end{exmp}

\bibliography{refs}

\bigskip

\noindent Princeton University, Princeton NJ 08544-1000

\begin{verbatim}kollar@math.princeton.edu\end{verbatim}

\end{document}